\documentclass[12pt]{article}
\usepackage[margin=.88in]{geometry}
\usepackage[utf8]{inputenc}
\usepackage{amsfonts}

\usepackage{amsmath,amssymb,amsfonts,amscd,amsxtra,amsthm}
\usepackage{bbm}

\usepackage[colorlinks,
linkcolor=blue,
anchorcolor=blue,
citecolor=blue]{hyperref}

\usepackage{graphics,graphicx} 
\usepackage{subcaption,caption,cleveref}

 \usepackage{xcolor}

\usepackage[round,authoryear]{natbib}
 
\newtheorem*{theorem}{Theorem}
 
\def\eps{\varepsilon}

\newtheorem{corollary}{Corollary}
\newtheorem{lemma}{Lemma}

\newtheorem{definition}{Definition}

\newtheorem{remark}{Remark}

\newcommand{\RR}{\mathbb{R}}

\newcommand{\PP}{\mathbb{P}}
\newcommand{\EE}{\mathbb{E}}
\newcommand{\OO}{\mathbb{O}}

\newcommand{\bfj}{\mathbf{j}}

\newcommand{\bfs}{\mathbf{s}}

\newcommand{\mcA}{\mathcal{A}}

\newcommand{\mcE}{\mathcal{E}} 
\newcommand{\mcF}{\mathcal{F}}

\newcommand{\mcR}{\mathcal{R}} 
\newcommand{\mcM}{\mathcal{M}}
\newcommand{\mcN}{\mathcal{N}}
\newcommand{\mcP}{\mathcal{P}} 
\newcommand{\mcQ}{\mathcal{Q}}
\newcommand{\mcS}{\mathcal{S}}
\newcommand{\mcT}{\mathcal{T}}
\newcommand{\mcU}{\mathcal{U}}

\newcommand{\wtU}{\widetilde{U}}

\newcommand{\wtLambda}{\widetilde{\Lambda}} 
\newcommand{\wtSig}{\widetilde{\Sigma}}

\newcommand{\whu}{\widehat{u}} 
\newcommand{\whU}{\widehat{U}} 

\newcommand{\whSig}{\widehat{\Sigma}}

\newcommand{\whl}{\widehat{\lambda}}
\newcommand{\whL}{\widehat{\Lambda}}
\newcommand{\whP}{\widehat{P}}

\newcommand{\tdU}{\tilde{U}}

\newcommand{\tdSig}{\tilde{\Sigma}}

\newcommand{\tdY}{\tilde{Y}}

\newcommand{\ot}{\otimes} 

\usepackage{algorithmic}
\usepackage{algorithm}

\makeatletter
\newcommand*{\algotitle}[2]{%
	\stepcounter{algocf}%
	\hypertarget{algocf.title.\theHalgocf}{}%
	\NR@gettitle{#1}%
	\label{#2}%
	\addtocounter{algocf}{-1}%
}
\makeatother

\title{Optimal Differentially Private PCA and Estimation for Spiked Covariance Matrices}  

\author{T. Tony Cai\footnote{Department of Statistics and Data Science, the Wharton School, University of Pennsylvania,  tcai@wharton.upenn.edu. The research of Tony Cai was supported in part by NIH grants R01-GM123056 and R01-GM129781.} ,  Dong Xia  \footnote{Department of Mathematics, Hong Kong University of Science and Technology, madxia@ust.hk. Dong Xia's research was partially supported by Hong Kong RGC Grant GRF 16302323.}, and  Mengyue Zha  \footnote{Department of Mathematics, Hong Kong University of Science and Technology, mzha@connect.ust.hk. } }

\date{(\today)}

\begin{document}
	
	\maketitle
	
\begin{abstract}
 
Estimating a covariance matrix and its associated principal components is a fundamental problem in contemporary statistics. While optimal estimation procedures have been developed with well-understood properties, the increasing demand for privacy preservation introduces new complexities to this classical problem. In this paper, we study optimal differentially private Principal Component Analysis (PCA) and covariance estimation within the spiked covariance model.

We precisely characterize the sensitivity of eigenvalues and eigenvectors under this model and establish the minimax rates of convergence for estimating both the principal components and  covariance matrix. These rates hold up to logarithmic factors and encompass general Schatten norms, including spectral norm, Frobenius norm, and nuclear norm as special cases.

We propose computationally efficient differentially private estimators and prove their minimax optimality for sub-Gaussian distributions, up to logarithmic factors. Additionally, matching minimax lower bounds are established. Notably, compared to the existing literature, our results accommodate a diverging rank, a broader range of signal strengths, and remain valid even when the sample size is much smaller than the dimension, provided the signal strength is sufficiently strong. Both simulation studies and real data experiments demonstrate the merits of our method.

\end{abstract}

\section{Introduction}

The covariance structure plays a fundamental role in multivariate analysis, and Principal Component Analysis (PCA) is a widely recognized technique known for its efficacy in dimension reduction and feature extraction \citep{anderson2003introduction}. PCA is particularly adept in settings where the data is high-dimensional but the underlying signal displays a low-dimensional structure. The estimation of covariance matrices and principal components finds applications across a diverse spectrum, encompassing tasks such as image recognition, data compression, clustering, risk management, portfolio allocation, mean tests, independence tests, and correlation analysis.
Methodologies and theoretical advancements, including minimax optimality, for covariance matrix estimation and PCA, have been well-established in both low-dimensional and high-dimensional settings. See, for example,  \cite{koltchinskii2017concentration, vershynin2012close, srivastava2013covariance,bickel2008regularized, cai2010optimal,ravikumar2011high, johnstone2001, cai_ma_wu2012spars, cai2015optimal, zhang2022heteroskedastic}. For a survey on optimal estimation of high-dimensional covariance structures, see \cite{cai2016estimating}.

Amidst the increasing availability of large datasets containing sensitive personal information, privacy concerns in statistical data analysis have gained heightened prominence. The utilization of personal information in statistical analyses raises apprehensions about the potential compromise of individual privacy. Consequently, there is a growing emphasis on developing methodologies and techniques that offer robust privacy guarantees while still facilitating accurate statistical insights. This motivates a comprehensive exploration of the optimal tradeoff between privacy and accuracy in fundamental statistical problems, including PCA and covariance matrix estimation.

Differential privacy (DP), a concept introduced by \cite{dwork2006proposeDP}, provides a framework for safeguarding individual privacy in statistical analysis.  DP has become a commonly accepted standard in both industrial and governmental applications \citep{google_privacy,ding2017collecting,apple2017,abowd2016challenge,abowd2020modernization}.  The goal of the present paper is to develop methods and optimality results for PCA and covariance matrix estimation within the framework of the spiked covariance model under DP constraints. 

\subsection{Problem formulation}
We begin by formally introducing the spiked covariance model and general formulation of the privacy constrained estimation problems.  

The spiked covariance structure \citep{johnstone2001, johnstone_lu2009} naturally arises from factor models with homoscedastic noise and has found diverse applications in signal processing, chemometrics, econometrics, population genetics, and various other fields. See, for example,  \cite{Fan08, Kritchman08, onat09, patterson}. The spiked covariance model assumes that the population covariance matrix can be decomposed as
\begin{align}\label{eq:spiked-cov}
	\Sigma=U\Lambda U^{\top}+\sigma^2 I_p,
\end{align}
where $U\in\OO_{p,r}$ and $\Lambda={\rm diag}(\lambda_1,\cdots,\lambda_r)$ represent the leading eigenvectors and eigenvalues (excluding $\sigma^2$), respectively. Here, $\OO_{p,r}$ denotes the set of $p\times r$ matrices satisfying $U^{\top}U=I_r$.   The spiked covariance model is convenient for studying the distribution of sample eigenvalues and eigenvectors, which play a critical role in the statistical inference of $\Sigma$ and its eigenvectors. For instance, \cite{donoho2018optimal} studied the optimal shrinkage of sample eigenvalues in the spiked covariance model. In particular, \cite{cai2015optimal} and \cite{zhang2022heteroskedastic} established the minimax optimal rates
\begin{equation}\label{eq:opt-spiked-cov}
	\begin{aligned}
		\inf_{\whU} \sup_{\Sigma\in\Theta (\lambda,\sigma^2)} &\EE \|\whU\whU^{\top}-UU^{\top}\|\asymp \bigg(\frac{\sigma^2}{\lambda}+\sqrt{\frac{\sigma^2}{\lambda}}\bigg)\sqrt{\frac{p}{n}};\\
		\inf_{\whSig} \sup_{\Sigma\in\Theta (\lambda,\sigma^2)} &\EE \|\whSig-\Sigma\| \asymp \lambda\sqrt{\frac{r}{n}}+\sqrt{\sigma^2(\lambda+\sigma^2)}\sqrt{\frac{p}{n}},
	\end{aligned}
\end{equation}
where the infimum is taken over all possible estimators based on the data $X=(X_1,\cdots,X_n)$ consisting of $n$ observations independently sampled  from the spiked covariance model (\ref{eq:spiked-cov}), the parameter set $\Theta(\lambda,\sigma^2)$ is a collection of covariance matrices in the form (\ref{eq:spiked-cov}) with all spiked eigenvalues have magnitudes of order $\lambda$ (see formal definition in (\ref{eq:def-Theta})), and $\|\cdot\|$ denotes the matrix spectral norm.

The concept of differential privacy was first introduced in \cite{dwork2006proposeDP}. For a given dataset $X$ and any $\varepsilon > 0$ and $\delta \in [0, 1)$, a randomized algorithm $A$ that maps $X$ into $\RR^{d_1\times d_2}$ is called $(\varepsilon, \delta)$-differentially private ($(\varepsilon, \delta)$-DP) over the dataset $X$ if
$$
\PP\big(A(X)\in\mcQ\big)\leq e^{\varepsilon}\PP\big(A(X')\in\mcQ\big)+\delta, 
$$
for all measurable subset $\mcQ\subset \RR^{d_1\times d_2}$ and all neighboring data set $X'$.  In the standard definition, a dataset $X'$ is a neighbor of $X$ if they differ by only one datum, i.e., one observation in $X$ is replaced by some other, possibly arbitrary, datum. In the context of PCA and covariance matrix estimation, as observations in $X$ are independently sampled from a common distribution, a neighboring dataset $X'$ is obtained by replacing one datum in $X$ with an independent copy. This facilitates exploration of the statistical properties of the sample data.

Under the $(\varepsilon,\delta)$-DP constraint, our goal is to investigate the cost of privacy in PCA and covariance matrix estimation.  This includes designing minimax optimal $(\varepsilon,\delta)$-DP estimators of the principal components and covariance matrix and establishing the privacy-constrained minimax lower bounds.  

\subsection{Main contribution}

In this paper, we establish the minimax optimal rates for PCA and covariance matrix estimation in the spiked model under DP constraints. Over the collection of sub-Gaussian distributions, these rates, up to logarithmic factors, are given by:
\begin{equation}\label{eq:dp-opt-spiked-cov}
	\begin{aligned}
		\inf_{\whU\in\mcU_{\varepsilon,\delta}}\ \sup_{\Sigma\in\Theta (\lambda,\sigma^2)} &\frac{\EE \|\whU\whU^{\top}-UU^{\top}\|_q}{r^{1/q}}\asymp \bigg(\frac{\sigma^2}{\lambda}+\sqrt{\frac{\sigma^2}{\lambda}}\bigg) \left(\sqrt{\frac{p}{n}} +\frac{p\sqrt{r}}{n\epsilon}\right)\bigwedge 1;\\
		\inf_{\whSig\in\mcM_{\varepsilon,\delta}}\ \sup_{\Sigma\in\Theta (\lambda,\sigma^2)} &\frac{\EE \|\whSig-\Sigma\|_q}{r^{1/q}} \asymp  \left( \lambda \left(\sqrt{\frac{r}{n}} + \frac{r^{3/2}}{n\varepsilon}\right) + \sqrt{\sigma^2(\lambda+\sigma^2)} \left(\sqrt{\frac{p}{n}} + \frac{\sqrt{r}p}{n \varepsilon} \right) \right)\bigwedge \lambda,
	\end{aligned}
\end{equation} 
where the infimum is taken over all possible $(\varepsilon,\delta)$-DP algorithms denoted by $\mcU_{\varepsilon,\delta}$ for principal components and $\mcM_{\varepsilon,\delta}$ for the covariance matrix. The expectation is taken with respect to the randomness of both the data and the differentially private algorithm. These rates hold in Schatten-$q$ norms for all $q \in [1, \infty]$, including spectral norm ($q=\infty$), Frobenius norm ($q=2$), and nuclear norm ($q=1$) as special cases. The rank $r$ can grow with respect to $p$ as long as $r\leq p/2$, and the sample size can be much smaller than $p$ as long as the signal-to-noise ratio (SNR) satisfies $\lambda/\sigma^2\geq C_1(\sqrt{p/n}+p/n)$. This condition is minimal since no consistent estimation is possible when this condition does not hold. 
To our knowledge, this represents the first comprehensive presentation of minimax optimal rates for PCA and covariance matrix estimation under DP constraints. For technical convenience and theoretical clarity, we focus on sub-Gaussian distributions in this paper. However, we believe that our results can be extended beyond sub-Gaussian distributions. For further details, see the discussion in Section~\ref{sec:discuss}.

Our contributions are multifod. Methodologically, we introduce $(\varepsilon,\delta)$-DP estimators for PCA and covariance matrices that are computationally efficient. Specifically, we employ the Gaussian mechanism for the sample spectral projector in differentially private PCA. Notably, our DP estimator for the covariance matrix incorporates a novel design to handle unknown orthogonal rotations. These estimators are shown to achieve minimax optimality, up to logarithmic factors. Theoretically, we provide a comprehensive understanding of the minimax optimal rates for PCA and covariance estimation under privacy constraints, valid across all Schatten norms. The derivation of minimax lower bounds employs Fano's lemma with a differential privacy constraint and the construction of well-separated spectral projectors based on the packing complexity of Grassmannians \citep{koltchinskii_xia15, zhang2018tensor}.

Differentially private PCA and covariance estimation are challenging because it is difficult to characterize a sharp sensitivity bound for the eigenvectors.  Our main technical contribution lies in a precise characterization of the sensitivity of the sample spectral projector $\whU\whU^{\top}$, quantifying its deviation when one datum $X_i$ is replaced by an independent copy $X_i'$. A key technical tool is an explicit spectral representation formula for $\whU\whU^{\top}$ adapted from \cite{xia2021normal}. We derive a similar formula specifically for the spiked covariance model, which is of independent interest. Based on this sharp sensitivity analysis, we  apply the Gaussian mechanism to achieve the upper bounds in (\ref{eq:dp-opt-spiked-cov}), up to logarithmic terms.

\subsection{Related work}
  
Minimax optimal rates under $(\epsilon, \delta)$-DP guarantees have been established for several statistical problems, such as mean estimation, linear regression, pairwise comparisons, matrix completion, factorization, generalized linear models (GLMs), and sparse GLMs \citep{cai2021cost, cai2023PrivateRanking, chien2021private, pmlr-v206-wang23d, cai2023score}. 
Additionally, optimality results have also been developed under local privacy constraints.  For example, 
\cite{duchi2018minimax} established minimax rates for mean estimation, GLMs, and nonparametric density estimation, while \cite{rohde2020geometrizing} developed minimax theory for estimating linear functionals under local  privacy. It is worth noting that local privacy is a stronger notion of privacy compared to $(\varepsilon, \delta)$-DP, and it may not be compatible with high-dimensional problems \citep{duchi2018minimax}. {A refined fingerprint lower bound method was introduced by \cite{narayanan2023better}, allowing for a broader range of $\delta$ and establishing a minimax lower bound for covariance matrix estimation (see also \cite{dong2022differentially} and \cite{mangoubi2022re}). Both studies focused on general covariance matrix estimation, but their results become suboptimal in the case of spiked covariance matrices. The Johnson-Lindenstrauss mechanism was examined by \cite{nikolov2023private}, providing optimal sample complexity for differentially private covariance estimation of a bounded high-dimensional distribution. While these privacy-preserving methods are centered on covariance estimation, their applicability and performance for PCA remain largely unclear. Additionally, although enforcing  boundedness can guarantee worst-case privacy protection, it may result in a pessimistic estimator in certain settings. See Remarks~\ref{rmk:worst-case} and \ref{rmk:dp-gauss} for a detailed comparison with existing literature. 

Differentially private PCA algorithms were proposed in \cite{blum2005dp-pca,chaudhuri2011differentially,dwork2014analyze} based on the perturbation mechanism, treating each datum $X_i$ as a fixed vector and investigating the sensitivity of sample eigenvectors. However, their deterministic sensitivity analysis disregards the statistical properties of sample data, resulting in suboptimal error rates when $X_i$'s are i.i.d. sampled from a common distribution, such as the spiked covariance model. Differentially private methods that explore statistical properties have been studied in \cite{brown2021covariance, kamath19bridge} and related works. However, optimal differentially private PCA has received much less attention, and existing results for private covariance estimation are generally suboptimal under the spiked covariance model.  Recently, \cite{liu_xiyang2022dp-pca} introduced an online PCA algorithm with DP, providing a much sharper upper bound for differentially private PCA under the spiked covariance model. The online Oja's algorithm in \cite{liu_xiyang2022dp-pca} consumes one datum at a time, allowing for an explicit representation formula in the updated estimate of eigenvectors and enabling a study of their sensitivity. However, their bound is valid only for the rank-one case ($r=1$) and is minimax optimal only when $\lambda\leq \sigma^2$. The optimality of their algorithm for general rank $r$ or $\lambda\gg \sigma^2$ remains unclear. Moreover, the minimax optimal rates for estimating $\Sigma$ under privacy constraints are still unknown under the spiked covariance model. 

\subsection{Organization of the paper}

The rest of the paper is organized as follows. In Section~\ref{sec: methodology}, we introduce the Gaussian mechanism and study the sensitivity of the empirical spectral projector under the spiked covariance model.  We present a DP algorithm for estimating the spectral projector and spiked covariance matrix in the same section. The upper bounds for our proposed DP algorithms are proven in Section~\ref{sec:dp-upb}, where an explicit spectral representation formula under the spiked covariance model is also developed. Section~\ref{sec:dp-lwb} establishes a differentially private Fano's lemma and minimax lower bounds. Extensions to the settings with diverging conditioning number and  sub-Gaussian distributions are discussed in Section~\ref{sec:discuss}. 
Some of the key technical lemmas are presented in Section \ref{sec:tech-lem}.  All the proofs as well as additional simulation results are given in  the Appendix.

\section{Methodology: Gaussian Mechanism and Sensitivity}	
\label{sec: methodology}

Our differentially private PCA and covariance estimation method relies on a precise characterization of the sensitivity for both eigenvectors and eigenvalues under the spiked covariance model.  
For technical convenience, we first focus on Gaussian PCA and provide a broader extension to general sub-Gaussian PCA in Section~\ref{sec:discuss}. 

For brevity, let $X:=(X_1,\cdots,X_n)$ represent the $p\times n$ matrix collecting all i.i.d. observations $X_i$ sampled from a centered normal distribution $\mcN(0, \Sigma)$.  The sensitivity of eigenvectors and eigenvalues denotes their perturbation if an observation $X_i$ is replaced by an independent copy $X_i'$ expressed briefly as $X^{(i)}:=(X_1,\cdots,X_{i-1},X_i', X_{i+1},\cdots, X_n)$. Here, $X$ and $X^{(i)}$ form a pair of neighboring datasets \citep{dwork2006proposeDP}. Notably, the sensitivity is contingent on the covariance matrix $\Sigma$.

Through out this paper, we consider the spiked covariance matrix model where $\Sigma$ is from the following parameter space 
\begin{equation}\label{eq:def-Theta}
	\begin{aligned}
		\Theta (p, r, \lambda, \sigma^2) = \Big\{\Sigma =   U \Lambda U^{\top}+& \sigma^2 I_p:\\
		&  U\in\OO_{p,r},  \Lambda={\rm diag}(\lambda_1,\cdots,\lambda_r), c_0\lambda \leq \lambda_r\leq \lambda_1\leq C_0\lambda\Big\}, 
	\end{aligned} 
\end{equation}
where $I_p$ is the identity matrix and and $\OO_{p,r}$ refers to the set of matrices with orthonormal columns, i.e., matrices satisfying $U^{\top}U=I_r$. Thus, our focus is on spiked covariance matrices with a bounded condition number, a common assumption in existing literature \citep{cai_ma_wu2012spars,chaudhuri2011differentially,liu_xiyang2022dp-pca}. However, our methodology remains valid, and the theoretical framework can be extended to the case of an unbounded condition number, as discussed in Section~\ref{sec:discuss}. For simplicity, we use $\Theta(\lambda,\sigma^2)$ without explicitly stating the dimensions $p$ and rank $r$. Let $\mcP$ denote the family of normal distributions $\mcN(0, \Sigma)$ with the population covariance matrix $\Sigma\in\Theta(\lambda, \sigma^2)$. Without loss of generality, we assume that $\sigma^2$ is known.

Formally, the sensitivity and Gaussian mechanism are described as follows without proofs. See, for example, \cite[Proposition 1]{dwork2006proposeDP}  and \cite[Theorem A.1]{dwork2014algorithmic} for more details.  Here, $\|\cdot\|_{\rm F}$ stands for the matrix Frobenius norm. 


\begin{lemma}[sensitivity and Gaussian mechanism] \label{lem:gaussian_mechanism}
	Let $X$ be a given data set and $X'$ be any neighboring data set of $X$, i.e., $X$ and $X'$ differs by at most one observation.  The sensitivity of a function $f$ that maps $X$ into $\RR^{d_1\times d_2}$ is defined by 
	\begin{equation}\label{eq:sensitivity-def}
		\omega_f:=\sup_{\textrm{neighboring}(X, X')} \|f(X)-f(X')\|_{\rm F}. 
	\end{equation}
	Then, for any $\varepsilon > 0$ and $\delta \in [0, 1)$,  the randomized algorithm $A$ defined by $A(X)=f(X)+Z$ where $Z$ has i.i.d. $\mcN\big(0, 2\omega_f^2\varepsilon^{-2}\log(1.25/\delta)\big)$  entries is $(\varepsilon, \delta)$-DP over the dataset $X$. 
\end{lemma}


The definition of sensitivity in Lemma~\ref{lem:gaussian_mechanism} relies on the pair of neighboring data sets. Here, $X$ is simply the data matrix where each column represents one observation. While $X$ and $X'$ differ only by one observation, the sensitivity can still be unbounded if no restriction is posed on the difference, e.g., by replacing one observation of $X$ by infinite. Since $X$ consists of i.i.d. columns under the spiked covariance model,  we assume that  a neighboring data set $X'$ is obtained by replacing some column of $X$ by its i.i.d. copy throughout this paper.

\subsection{Differentially private estimation by Gaussian mechanism}
Our DP-estimators of principal components and covariance matrix are built on Gaussian mechanism. Here, we assume that the rank $r$ and nuisance variance $\sigma^2$ are known for simplicity. Let $\whU$ be the top-$r$ eigenvectors of the sample covariance matrix $\whSig:=n^{-1}\sum_{i=1}^n X_iX_i^{\top}$ and denote $\whU\whU^{\top}$ the sample spectral projector.   By Lemma~\ref{lem:gaussian_mechanism},  differentially private PCA can be obtained by adding Gaussian noise $Z$ to $\whU\whU^{\top}$ provided that the entrywise variance of $Z$ dominates the sensitivity of $\whU\whU^{\top}$.  While publishing $\whU\whU^{\top}+Z$ protects privacy,  it is certainly not a preferable estimator of principal components as it generally lacks validity as a spectral projector.  We therefore take the eigenvectors of $\whU\whU^{\top}+Z$ as the ultimate estimator. This choice maintains differential privacy, as the post-processing of a differentially private algorithm retains differential privacy according to well-established results, as discussed in \cite{dwork2006proposeDP}. 

\begin{algorithm}
	\caption{Differentially private PCA and covariance estimation}\label{alg:DP-PCA}
	\begin{algorithmic}
		\STATE{\textbf{Input}: data matrix $X= (X_1, \cdots, X_n)\in\RR^{n\times p}$;  eigenvectors and eigenvalues sensitivity $\Delta_1$ and  $\Delta_2>0$; rank $r$; nuisance variance $\sigma^2$;  privacy budget $\varepsilon>0,  \delta \in(0,1)$.}
		\STATE{\textbf{Output}: $(\varepsilon,\delta)$-DP estimate of $U$ and $\Sigma$}. 
		\STATE{Compute the sample covariance matrix and top-$r$ eigenvectors:
			$$
			\whSig \longleftarrow \frac{1}{n}\sum_{i=1}^n X_iX_i^{\top}\quad {\rm and}\quad \whU\longleftarrow {\rm SVD}_r(\whSig);
			$$
			Compute $(\varepsilon/2,  \delta/2)$-DP PCA by adding artificial Gaussian noise:
			$$
			\widetilde U\longleftarrow {\rm SVD}_r\Big(\whU\whU^{\top}+Z\Big)\quad {\rm where}\quad Z_{ij}=Z_{ji}\stackrel{{\rm i.i.d.}}{\sim} \mcN\Big(0,  \frac{8\Delta_1^2}{\varepsilon^2}\log \frac{2.5}{\delta}\Big),\quad \forall 1\leq i\leq j\leq p;
			$$
			Compute $(\varepsilon/2,  \delta/2)$-DP estimates of eigenvalues up to rotations:
			\begin{align*}
				\wtLambda&\quad \longleftarrow \wtU^{\top} \big(\whSig-\sigma^2I_p\big)\wtU + E\quad {\rm where}\quad E_{ij}=E_{ji}\stackrel{i.i.d.}{\sim} \mcN\Big(0,  \frac{8\Delta_2^2}{\varepsilon^2}\log \frac{2.5}{\delta}\Big),\quad \forall 1\leq i\leq j\leq r;
			\end{align*} 
			Compute $(\varepsilon,  \delta)$-DP covariance estimate by :
			$$
			\widetilde\Sigma \longleftarrow \wtU\wtLambda \wtU^{\top}+\sigma^2 I_p.
			$$
			\textbf{Return}: $\wtU$ and $\wtSig$}
	\end{algorithmic}
\end{algorithm}

Our proposed differentially private PCA and covariance estimation procedures are given in Algorithm~\ref{alg:DP-PCA}.  The proper choice of sensitivities $\Delta_1$ and $\Delta_2$ is determined by Lemma~\ref{lem:sense-U} and Lemma~\ref{lem:sense-lambda} in Section~\ref{sec:sense}, respectively.
However,  $\wtU$ and $\whU$ are close up to an orthogonal rotation.  As a result,  our algorithm chooses to add Gaussian noise to $\wtU^{\top} \whSig \wtU$ instead of the empirical eigenvalues $\whL:=(\widehat \lambda_1,\cdots,\widehat\lambda_r)^{\top}$.  The added noise level depends on the sensitivity of $\wtU^{\top} \whSig \wtU$,  within which $\wtU$ is already differentially private.  It thus suffices to study the upper bound of $\|\wtU^{\top}(\whSig-\whSig^{(i)})\wtU\|_{\rm F}\leq \|\whSig-\whSig^{(i)}\|_{\rm F}$,  which will be established in Lemma~\ref{lem:sense-lambda}.  

Our approach to differentially privately estimating the main covariance term involves separately privatizing the eigenvectors and eigenvalues. This separation is driven by the observation that the relative sensitivity of eigenvalues is significantly larger than that of eigenvectors.  Note that a natural estimator of $U(\Lambda+\sigma^2I_r)U^{\top}$ is $\whU\whU^{\top}\whSig \whU\whU^{\top}$. It is possible to characterize the sensitivity of this estimator by directly studying the bound $\|\whU\whU^{\top}\whSig \whU\whU^{\top}-\whU^{(i)}\whU^{(i)\top}\whSig^{(i)} \whU^{(i)}\whU^{(i)\top}\|_{\rm F}$. However, the sensitivity of eigenvalues will be the dominating factor and force us to add unnecessarily large noise to a $p\times p$ matrix.  This delivers a statistically sub-optimal estimator of the spiked covariance matrix. 

The estimated eigenvectors $\wtU$ is $(\varepsilon/2,  \delta/2)$-DP and eigenvalues $\wtLambda$ is $(\varepsilon/2, \delta/2)$-DP with high probability.   By the composition property of differentially private algorithm,  the estimator $\wtU \wtLambda \wtU^{\top}$ is $(\varepsilon, \delta)$-DP. The conclusion is formally stated in the following lemma. Recall that $\tilde{r}=(r\lambda+p\sigma^2)/(\lambda+\sigma^2)$ is the effective rank of $\Sigma$. Here, $\lambda$ is regarded as the signal strength. 



\begin{lemma}\label{lem:DP-alg}
	Let the data matrix $X=(X_1,\cdots,X_n)$ consists of i.i.d. columns sampled from $\mcN(0, \Sigma)$ with $\Sigma\in\Theta(\lambda,\sigma^2)$, $\varepsilon>0, \delta\in(0,1)$, and assume $n\geq C_1(r\log n+\log^2n),  2r\leq p$, and $\lambda/\sigma^2\geq C_1(p/n+\sqrt{p/n})$ for some large absolute constant $C_1>0$.  If we choose 
	$$
	\Delta_1:=C_2\bigg(\frac{\sigma^2}{\lambda}+\sqrt{\frac{\sigma^2}{\lambda}}\bigg)\frac{\sqrt{p(r+\log n)}}{n}\quad {\rm and}\quad \Delta_2:=C_3\frac{\lambda(r+\log n)+\sigma^2 (p+\log n)}{n}, 
	$$
	for some large absolute constants $C_2, C_3>0$, then Algorithm~\ref{alg:DP-PCA} is $(\varepsilon, \delta)$-DP with probability at least $1-4n^{-99}-e^{-c_1(n\wedge p)}$ for some absolute constant $c_1>0$. 
\end{lemma}

\begin{remark}[Worst-case and high-probability privacy guarantee]\label{rmk:worst-case}
Compared to existing literature \citep{chaudhuri2011differentially, liu_xiyang2022dp-pca, kamath2019privately, nikolov2023private, dong2022differentially}, our algorithm does not truncate the observations, allowing $\|X_i\|$ to remain unbounded.  As a result, our algorithm is differentially private with high probability. The \textsf{DP-Oja} algorithm proposed by \cite{liu_xiyang2022dp-pca} ensures worst-case privacy guarantees due to its online nature. However, it is limited to the rank-one case, performs poorly in both simulation and real data experiments (see Section~\ref{sec:numeric}),  and the established error rate is much larger than ours under spiked covariance model when signal strength $\lambda\gg \sigma^2$ (see Remark~\ref{rmk:dp-oja}).  The \textsf{DP-Gauss} method \citep{dwork2014analyze, mangoubi2022re} ensures worst-case privacy by applying a global scaling, limiting each observation to at most unit norm. While we could apply global scaling to our method to ensure worst-case privacy, as discussed in Remark~\ref{rmk:dp-gauss}, this approach would result in an overly pessimistic estimator with a significantly larger error rate under the spiked covariance model. Therefore, we do not pursue worst-case privacy guarantees in this paper.
 Moreover,  note that the probability terms $n^{-99}$ in Lemma~\ref{lem:DP-alg} can be replaced by $n^{-C_5}$ with any absolute constant $C_5>0$ (by adjusting the constants $C_2, C_3$ in the definitions of $\Delta_1$ and $\Delta_2$ accordingly).  The failure probability decreases polynomially fast with respect to sample size $n$. 
\end{remark}

The sensitivities $\Delta_1$ and $\Delta_2$ play a critical role in guaranteeing the differential privacy of Algorithm~\ref{alg:DP-PCA}, which shall be developed in next section. The conditions $r\log n+\log^2n=O(n)$ and $2r\leq p$ are mild. The SNR condition $\lambda/\sigma^2\geq C_1(p/n+\sqrt{p/n})$ is typical in the existing literature of spiked covariance matrix model. See, e.g., \cite{nadler2008finite, zhang2022heteroskedastic} and references therein.

\subsection{Sensitivity analysis}\label{sec:sense}
In this section, we analyze the sensitivities of sample eigenvectors and eigenvalues under the spiked covariance model. 
The data matrix $X=(X_1,\cdots,X_n)\sim \mcN(0,\Sigma)^{\otimes n}$ for some $\Sigma\in\Theta(\lambda,\sigma^2)$. Similarly, its neighboring data matrix $X^{(i)}=(X_1,\cdots,X_i', \allowbreak \cdots,X_n)\sim \mcN(0,\Sigma)^{\otimes n}$. Define the sample covariance matrices by 
$$
\whSig:=\frac{1}{n}\sum_{i=1}^n X_i X_i^{\top}\quad{\rm and}\quad \whSig^{(i)}: = \frac{1}{n} \Big( X_i^{\prime} X_i^{\prime\top} + \sum_{j\neq i} X_j X_j^{\top}  \Big), 
$$
Denote $\whU\in\OO_{p,r}$ and $\whU^{(i)}\in \OO_{p,r}$ the top-$r$ left eigenvectors of $\whSig$ and $\whSig^{(i)}$, respectively. The sensitivity of sample eigenvectors characterizes the deviation between $\whU$ and $\whU^{(i)}$ caused by replacing the $i$-th observation by its i.i.d. copy. Since eigenvectors are determined up to an orthogonal rotation (note that we allow the eigengap $|\lambda_i-\lambda_j|$ to be zero), a commonly used metric for measuring the distance between eigenvectors is the projection distance defined by $\|\whU\whU^{\top}-\whU^{(i)}\whU^{(i)\top}\|_{\rm F}$. 

The primary challenge in differentially private PCA  lies in characterizing a precise upper bound for $\|\whU\whU^{\top}-\whU^{(i)\top}\whU^{(i)}\|_{\rm F}$. In most existing works \citep{blum2005dp-pca,chaudhuri2011differentially,dwork2014analyze}, the data matrix $X$ is assumed to be fixed, and its columns are all bounded, denoted as $\|X_i\|\leq \gamma$, where we slightly abuse the notation by letting $\|\cdot \|$ denote the $\ell_2$-norm for vectors and $\gamma$ is a deterministic value. This immediately implies an upper bound $\|\whSig-\whSig^{(i)}\|\leq 2\gamma^2/n$ and the sensitivity of $\whU\whU^{\top}$ is guaranteed by the Davis-Kahan theorem. 

However, this approach becomes invalid when observations are unbounded and sub-optimal when observations are randomly sampled from a common distribution. A more recent work \cite{liu_xiyang2022dp-pca} aimed to exploit the statistical properties of i.i.d. samples to achieve a sharper bound for differentially private PCA. This work focused on the rank-one case ($r=1$) and the Oja's algorithm, well-known for online PCA, which iteratively updates the estimation with one additional observation. The online fashion of Oja's algorithm in the rank-one case allows for an explicit representation of the eigenvector estimator, enabling a sharp upper bound of the sensitivity to be derived. Consequently, nearly optimal differentially private PCA for the case $r=1$ was achieved. However, it remains unclear how this approach can be extended to the rank-$r$ case and what the minimax optimal convergence rates are.

We take a fundamentally different approach by directly focusing on $\|\whU\whU^{\top}-\whU^{(i)}\whU^{(i)\top}\|_{\rm F}$. This task presents two challenges: the spectral projector $\whU\whU^{\top}$ involves a complicated  function of the data matrix $X$, and a sharp perturbation analysis is required for a set of $r$ empirical eigenvectors. Fortunately, we leverage an explicit spectral representation formula adapted from \cite{xia2021normal} and successfully establish a precise upper bound for $\|\whU\whU^{\top}-\whU^{(i)}\whU^{(i)\top}\|_{\rm F}$.

\begin{lemma} 	\label{lem:sense-U} 
Suppose the conditions in Lemma~\ref{lem:DP-alg} hold and assume $n\geq C_1(r\log n+\log^2n)$ and $2r\leq p$. 
	There exist absolute constants $c_1, C_2>0$ such that  with probability at least $1-3n^{-99}-e^{-c_1(n\wedge p)}$, 
	\begin{equation}\label{eq:lem-sense-U-bd}
		\max_{i\in[n]}\|\widehat{U}\widehat{U}^{\top} -\widehat{U}^{(i)}\widehat{U}^{(i)\top}\|_{\rm F}\leq C_2 \bigg(\frac{\sigma^2}{\lambda}+\sqrt{\frac{\sigma^2}{\lambda}}\bigg) \frac{\sqrt{p(r+\log n)}}{n}. 
	\end{equation}
\end{lemma}

The $\log n$ term in upper bound (\ref{eq:lem-sense-U-bd}) is due to the maximization over $n$.  Nevertheless, the bound is much smaller than that achieved by the deterministic analysis in  \citep{blum2005dp-pca,chaudhuri2011differentially,dwork2014analyze}. Indeed, a direct application of Davis-Kahan theorem yields an upper bound $O\big(\|\whSig-\whSig^{(i)}\|\sqrt{r}/\lambda\big)$, which is at least in the order $O\big((r\lambda+p\sigma^2)\sqrt{r}/(n\lambda)\big)$, with high probability. The significant improvement is due to a sharp spectral characterization showing that the difference $\widehat{U}\widehat{U}^{\top} -\widehat{U}^{(i)}\widehat{U}^{(i)\top}$ is mainly contributed by the term $\|U^{\top}(X_iX_i^{\top}-X_i'X_i'^{\top})U_{\perp}\|_{\rm F}/(n\lambda)$. Here,  $U_{\perp}\in \OO_{p,p-r}$ denotes the orthogonal complement of $U$ such that $(U, U_{\perp})$ is an orthogonal matrix. The proof of Lemma~\ref{lem:sense-U} is technically involved and deferred to Section \ref{appendix_lem:sense-U}.  It is worth noting that the original  spectral representation formula developed in \cite{xia2021normal} is inapplicable here because $\Sigma$ is not exactly rank-$r$.   Interestingly,  we establish a similar spectral representation formula exclusively for spiked covariance matrix,  which may be of independent interest.  See Lemma~\ref{lem:spectral-formula} in Section~\ref{sec:spectral-formula}.

The sensitivity of eigenvalues is also necessary for constructing differentially private covariance estimation. Let $\lambda_k(\whSig)$ and $\lambda_k(\whSig^{(i)})$ denote the $k$-th largest eigenvalue of $\whSig$ and $\whSig^{(i)}$, respectively. Compared to the eigenvectors, the sensitivity of eigenvalues can be easily characterized by Hoffman-Weilandt's inequality.  The proof of Lemma \ref{lem:sense-lambda} is deferred to Section~\ref{sec: proofs_DP-SPIKED-COV_is_DP}. 

\begin{lemma}	\label{lem:sense-lambda} 
Suppose the conditions in Lemma~\ref{lem:DP-alg} hold. 
	There exists an absolute constant $C_2>0$ such that with probability at least $1-n^{-100}$, 
	\begin{equation}\label{eq:lem-sense-lambda-bd}
		\sum_{k=1}^p\Big| \lambda_k(\widehat{\Sigma}) - \lambda_k(\widehat{\Sigma}^{(i)}) \Big|^2 \leq C_2\bigg(\frac{\lambda (r+\log n) + \sigma^2 (p+\log n)}{n}\bigg)^2, 
	\end{equation}
	for all $i\in[n]$.

\end{lemma} 

We can regard $\Big(\sum_{k=1}^p\big(\lambda_k(\whSig)-\lambda_k(\whSig^{(i)})\big)^2\Big)^{1/2}/\lambda$ and $\|\whU\whU^{\top}-\whU^{(i)}\whU^{(i)\top}\|_{\rm F}/\sqrt{r}$ as the relative sensitivity of eigenvectors and eigenvalues, respectively. Lemmas \ref{lem:sense-U} and \ref{lem:sense-lambda} show that the relative sensitivity of eigenvalues can be considerably larger than that of eigenvectors. This insight implies that, when designing a differentially private optimal estimation procedure for the population covariance matrix, it is advisable to privatize the eigenvalues and eigenvectors separately, as elaborated in Algorithm~\ref{alg:DP-PCA}.

\section{Upper Bounds with Differential Privacy}	
\label{sec:dp-upb}

\subsection{Spectral representation formula} 
\label{sec:spectral-formula}

Our key technical tool is the following spectral representation formula.  Recall that $\whU$ and $U$ denote the top-$r$ eigenvectors of $\whSig$ and $\Sigma$, respectively.   Denote the deviation matrix by $\widehat{\Delta}:=\whSig-\Sigma$ so that $\whSig=\Sigma+\widehat{\Delta}$  is viewed as a perturbation of the ``signal" matrix $\Sigma$.  The spectral representation formula was first introduced in  \cite{xia2021normal},  which,  however,  requires the ``signal" matrix to be exactly rank-$r$.  This is certainly not the case here since $\Sigma$ is full-rank.   Here, we develop the spectral representation formula exclusively for the perturbation of a spiked covariance matrix.  

The spectral representation formula is actually deterministic.  Let the symmetric matrix $\Delta\in\RR^{p\times p}$ be an arbitrary perturbation.  Denote $\whU$ the top-$r$ eigenvectors of $\Sigma+\Delta$ where $\Sigma=U\Lambda U^{\top}+\sigma^2 I_p$ with $\Lambda={\rm diag}(\lambda_1,\cdots,\lambda_r)$.     
We are interested in developing an explicit representation formula for the spectral projector $\whU\whU^{\top}$ in terms of $\Delta$.  Let $Q^{\perp}:=U_{\perp}U_{\perp}^{\top}=I_p-UU^{\top}$ denotes the orthogonal projection.  For all $t\geq 1$,  we define $Q^{-t} := U\Lambda^{-t} U^{\top}$.  We slightly abuse the notation and denote 
$Q^{0} :=Q^{\perp}=U_{\perp}U_{\perp}^{\top}$. 

\begin{lemma}
	\label{lem:spectral-formula} 
	Suppose that $\Sigma$ is a spiked covariance matrix as in (\ref{eq:spiked-cov}) and $2\|\Delta\|\leq \lambda_r$, then 
	\begin{equation*}
		\whU \whU^{\top}-U U^{\top}=\sum_{k \geq 1} \mathcal{S}_{\Sigma, k}(\Delta), 
	\end{equation*}
	where the $k$-th order term $\mathcal{S}_{\Sigma, k}(\Delta)$ is a summation of $\binom{2k}{k}$ terms defined by
	\begin{equation*}
		\mathcal{S}_{\Sigma, k}(\Delta)=\sum_{\mathbf{s}: s_1+\ldots+s_{k+1}=k}(-1)^{1+\tau(\mathbf{s})} \cdot Q^{-s_1} \Delta Q^{-s_2} \ldots \Delta Q^{-s_{k+1}}, 
	\end{equation*}
	where $\mathbf{s}=\left(s_1, \ldots, s_{k+1}\right)$ contains non-negative indices and $\tau(\mathbf{s})=\sum_{j=1}^{k+1} \mathbb{I}\left(s_j>0\right).$  A simple upper bound of the $k$-th order term is 
	$$
	\big\|\mcS_{\Sigma, k}(\Delta) \big\|\leq {2k\choose k}\Big(\frac{\|\Delta\|}{\lambda_r}\Big)^k. 
	$$ 
\end{lemma}

Based on Lemma~\ref{lem:spectral-formula},  the leading term,  i.e., the $1$st-order term,  of $\whU\whU^{\top}-UU^{\top}$ is contributed by $\Lambda^{-1}U^{\top}\Delta U_{\perp}$ and $U_{\perp}^{\top}\Delta U\Lambda^{-1}$.  The latter terms can be sharply controlled by exploiting the statistical properties of $\Delta$ if observations are i.i.d.  sampled.  

\subsection{Upper bounds}	

In this section,  we present the upper bounds of our $(\varepsilon,\delta)$-DP estimator $\wtU\wtU^{\top}$ and $\wtSig$.  In this section, we focus on the Gaussian setting, with an extension to the sub-Gaussian case provided in Section~\ref{sec:sub-Gaussian}. Cases beyond sub-Gaussian distributions are discussed in Section~\ref{sec:discuss-beyond}.
Let $\|\cdot\|_q$ denotes the matrix Schatten-$q$ norm for any $q\in[1,\infty]$,  e.g.,  the spectral norm $\|\cdot\|$ if $q=\infty$,  the Frobenius norm $\|\cdot\|_{\rm F}$ if $q=2$,  and the nuclear norm $\|\cdot\|_{\ast}$ if $q=1$.   
A straightforward application of the triangle inequality
$$
\|\wtU\wtU^{\top}-UU^{\top}\|_{q}\leq \|\wtU\wtU^{\top}-\whU\whU^{\top}\|_{q}+\|\whU\whU^{\top}-UU^{\top}\|_{q}, 
$$
leads to the following theorem. 

\begin{theorem} \label{thm:dp-pca-upb}
	Suppose that $X_1,\cdots, X_n\stackrel{{\rm i.i.d.}}{\sim}\mcN(0, \Sigma)$, $n\geq C_1(r\log n+\log^2n),  2r\leq p$, and $\lambda/\sigma^2\geq C_1(p/n+\sqrt{p/n})$ for some large absolute constant $C_1>0$.  If we choose 
	$$
	\Delta_1:=C_2\bigg(\frac{\sigma^2}{\lambda}+\sqrt{\frac{\sigma^2}{\lambda}}\bigg)\frac{\sqrt{p(r+\log n)}}{n}, 
	$$
	then, there exist absolute constants $c_1,  C_4>0$ such that,  for any $\varepsilon>0, \delta\in(0,1)$,  Algorithm~\ref{alg:DP-PCA} outputs an $(\varepsilon,\delta)$-DP estimator $\wtU\wtU^{\top}$ satisfying 
	\begin{align*}
		\frac{\|\wtU\wtU^{\top}-UU^{\top}\|_q}{r^{1/q}} \leq C_4\bigg(\frac{\sigma^2}{\lambda}+\sqrt{\frac{\sigma^2}{\lambda}}\bigg)\bigg(\sqrt{\frac{p}{n}}+\frac{p\sqrt{r+\log n}}{n\varepsilon}\sqrt{\log\frac{2.5}{\delta}}\bigg),
	\end{align*} 
	with probability at least $1-e^{-c_1(n\wedge p)}$.  Moreover, if $\lambda/\sigma^2\leq (p/n)e^{c_2(p\wedge n)}$ for some small absolute constant $c_2>0$, then 
	\begin{align*}
		\frac{\EE\|\wtU\wtU^{\top}-UU^{\top}\|_q}{r^{1/q}} \leq C_4\bigg(\frac{\sigma^2}{\lambda}+\sqrt{\frac{\sigma^2}{\lambda}}\bigg)\bigg(\sqrt{\frac{p}{n}}+\frac{p\sqrt{r+\log n}}{n\varepsilon}\sqrt{\log\frac{2.5}{\delta}}\bigg).
	\end{align*} 
	Here, $q$ can be any number in $[1,\infty]$.   
\end{theorem} 

Basically,  the upper bounds consist of two parts: the first one represent the statistical error rate and the second one is the cost of privacy constraint.  It is well-known that the first term is minimax optimal \citep{nadler2008finite, cai2015optimal,  koltchinskii2017concentration}.  The second term decays at the  rate $O\big(p/(n\varepsilon)\log^{1/2}\delta^{-1}\big)$ with respect to  the sample size,  dimension and privacy-related parameters,  which is typical in differentially private algorithms \citep{cai2023score,liu_xiyang2022dp-pca}.  In Section~\ref{sec:dp-lwb}, we shall develop matching minimax lower bounds showing that the rates in Theorem~\ref{thm:dp-pca-upb} are minimax optimal up the $\log n$ and $\log(2.5/\delta)$ terms. 

It worth to mention that the $\log n$ term appearing in the privacy-related rate is due to the requirement of differential privacy that applies to each of the $n$ observations.  This $\log n$ term seems to be present in the upper bounds of most differentially private algorithms.  See,  e.g.,  \cite{cai2021cost,  cai2023score,  dwork2014analyze} and references therein.  A slight difference here is that the $\log n$ term appears not as an additional factor,  but as an additive term.  If $r\geq \log n$,  the logarithmic factor can be ignored and the rate becomes minimax optimal except for the $\log\delta^{-1}$ factor.  

\begin{remark}[Comparison with Oja's algorithm \cite{liu_xiyang2022dp-pca}]\label{rmk:dp-oja}
The \textsf{DP-Oja} algorithm introduced in \cite[Corollary 5.2]{liu_xiyang2022dp-pca} delivered a rank-one ($r=1$) PCA estimator achieving the following error bound,  with probability $0.99$: 
$$
\|\widehat u_{\textsf{\tiny oja}}\widehat u_{\textsf{\tiny oja}}^{\top}-uu^{\top}\|=\widetilde O\Bigg(\bigg(1+\frac{\sigma^2}{\lambda}\bigg)\cdot \bigg(\sqrt{\frac{p}{n}}+\frac{p\sqrt{\log 1/\delta}}{\eps n}\bigg)\Bigg),
$$
where $\widetilde O(\cdot)$ hides logarithmic factors in $n$ and $p$.  Their established upper bound is much larger than ours when the signal strength $\lambda\gg \sigma^2$,  and their failure probability is a constant while ours decay polynomially fast as sample size $n$ increases.  
\end{remark}


We now present the performance bound for the differentially private estimator $\wtSig$.  

\begin{theorem}  \label{thm:dp-Sigma-upb}
	Suppose that  $X_1,\cdots, X_n\stackrel{{\rm i.i.d.}}{\sim}\mcN(0, \Sigma)$, $n\geq C_1(r\log n+\log^2n),  2r\leq p$, and $\lambda/\sigma^2\geq C_1(p/n+\sqrt{p/n})$ for some large absolute constant $C_1>0$.  If we choose 
	$$
	\Delta_1:=C_2\bigg(\frac{\sigma^2}{\lambda}+\sqrt{\frac{\sigma^2}{\lambda}}\bigg)\frac{\sqrt{p(r+\log n)}}{n}\quad {\rm and}\quad \Delta_2:=C_2\frac{\lambda(r+\log n)+\sigma^2 (p+\log n)}{n}, 
	$$
	then, there exist absolute constants $c_1,  C_4>0$ such that,  for any $\varepsilon>0, \delta\in(0,1)$,  Algorithm~\ref{alg:DP-PCA} outputs an $(\varepsilon,\delta)$-DP estimator $\wtSig$ satisfying 
	\begin{align*}
		&\frac{\|\wtSig-\Sigma\|_q}{r^{1/q}} \\
		&\leq C_4\Bigg(\lambda\bigg(\sqrt{\frac{r}{n}}+\frac{\sqrt{r}(r+\log n)}{n\varepsilon}\cdot \sqrt{\log\frac{2.5}{\delta}}\bigg)+\sqrt{\sigma^2(\lambda+\sigma^2)} \bigg(\sqrt{\frac{p}{n}}+\frac{p\sqrt{(r+\log n)}}{n\varepsilon}\sqrt{\log\frac{2.5}{\delta}}\bigg)\Bigg),
	\end{align*}
	with probability at least $1-3n^{-99}-e^{-c_1(n\wedge p)}$.  Moreover, if $\lambda/\sigma^2\leq (p/n)e^{c_2(p\wedge n)}$ for some small absolute constant $c_2>0$, then
	\begin{align*}
		&\frac{\EE\|\wtSig-\Sigma\|_q}{r^{1/q}} \\
		&\leq C_4\Bigg(\lambda\bigg(\sqrt{\frac{r}{n}}+\frac{\sqrt{r}(r+\log n)}{n\varepsilon}\cdot \sqrt{\log\frac{2.5}{\delta}}\bigg)+\sqrt{\sigma^2(\lambda+\sigma^2)} \bigg(\sqrt{\frac{p}{n}}+\frac{p\sqrt{(r+\log n)}}{n\varepsilon}\sqrt{\log\frac{2.5}{\delta}}\bigg)\Bigg). 
	\end{align*}
\end{theorem} 

By Theorem~\ref{thm:dp-Sigma-upb}, the privacy-irrelevant error rate 
$$
\lambda\sqrt{\frac{r}{n}}+\sqrt{\sigma^2(\lambda+\sigma^2)}\sqrt{\frac{p}{n}}, 
$$
matches the minimax optimal rate of spiked covariance estimation in the existing literature \citep{cai2015optimal,cai2010optimal}. For ease of discussion,  let us focus on the error rate in spectral norm. There are two terms related to the cost of privacy:
$$
\lambda\cdot \frac{\sqrt{r}(r+\log n)}{n\varepsilon}\sqrt{\log\frac{2.5}{\delta}}\quad {\rm and}\quad \sqrt{\sigma^2(\lambda+\sigma^2)} \bigg(\sqrt{\frac{p}{n}}+\frac{p\sqrt{(r+\log n)}}{n\varepsilon}\sqrt{\log\frac{2.5}{\delta}}\bigg),
$$
where the second term is approximately of order $\lambda\|\wtU\wtU^{\top}-UU^{\top}\|$,  contributed by the cost of estimating the eigenvectors. The first term grows at the rate $O(r^{3/2})$ with respect to the rank, which is contributed by the cost of estimating the eigenvalues. Due to the unknown orthogonal rotation measuring the alignment between $\whU$ and $\wtU$, privacy cost is also paid for the $r\times r$ unknown rotation matrix. 
Minimax lower bounds are developed in Section~\ref{sec:dp-lwb} demonstrating the optimality of these bound up to the $\log n$ and $\log(2.5/\delta)$ related terms.

\begin{remark}[Comparison with \cite{dwork2014analyze} and \cite{mangoubi2022re}]\label{rmk:dp-gauss}
The \textsf{DP-Gauss} method is a privacy-preserving low-rank approximation method originally proposed by \cite{dwork2014analyze} and later improved by \cite{mangoubi2022re}. The method applies the Gaussian mechanism to find the rank-$r$ approximation of the sample covariance matrix $\widehat \Sigma$, denoted by $\widehat \Sigma_r$ hereafter. Under the spiked covariance model, the \textsf{DP-Gauss} method provides an $(\epsilon,\delta)$-DP estimator, denoted as $\widetilde \Sigma_r$, achieving the rate (Corollary 2.3 in \cite{mangoubi2022re})
\begin{align}\label{eq:DP-Gauss-upb}
\big\|\widetilde \Sigma_r-\widehat \Sigma_r \big\|_{\rm F}\leq C_3 \max_{i\in[ n]}\|X_i\|^2\bigg(1+\frac{\sigma^2}{\lambda}\bigg)\frac{\sqrt{rp}}{n\eps}\log^{1/2}(1/\delta),
\end{align}
which can be viewed as the cost of privacy in their method. Note that the term $\max_{i\in[n]}\|X_i\|_2$ appears here because \cite{dwork2014analyze} and \cite{mangoubi2022re} require that each observation has at most unit norm. Under the spiked model, we have $\max_{i}\|X_i\|^2\asymp r\lambda +p\sigma^2$ up to $\log n$ factors. Plugging this into Equation (\ref{eq:DP-Gauss-upb}), we can conclude that the bound attained by \textsf{DP-Gauss} in \cite{mangoubi2022re} is much larger than ours under the spiked covariance model. 
\end{remark}

Private covariance estimation for Gaussian distributions was studied by \cite{kamath2019privately} using the Gaussian mechanism. Their rate is optimal in the case $\lambda \asymp \sigma^2$ and $r = p$, but becomes suboptimal otherwise. In contrast, our rate is optimal, allowing a much more relaxed condition on $\lambda$ and $\sigma^2$. Moreover, their method cannot be applied to differentially private PCA.

\begin{remark}[High-dimensional data]
Our methods work as long as the signal-to-noise ratio satisfies $\lambda/\sigma^2\gtrsim p/n+\sqrt{p/n}$,  meaning that a strong signal is required when the dimension $p$ is much larger than the sample size $n$.  However,  we emphasize that such a signal strength condition is necessary for a non-trivial estimate of the population eigenvectors, even in the conventional non-private setting.  See, e.g. ,  \cite{koltchinskii2017concentration}.  Moreover,  simulation results in Section~\ref{sec:numeric} and Appendix~\ref{sec:app-simulation} demonstrate that our method is much more robust than other methods (\textsf{DP-Oja} \citep{liu_xiyang2022dp-pca} and \textsf{DP-Gauss} \citep{dwork2014analyze,  mangoubi2022re}) when dimension $p$ is relatively large compared to the sample size $n$.  
\end{remark}

\section{Minimax Lower Bounds with Differential Privacy}
\label{sec:dp-lwb}

In this section, we establish the minimax lower bound of  PCA and covariance matrix estimation under the constraint of differential privacy. Our main technical tool is a version of Fano's lemma with privacy constraint. 

\subsection{DP-constrained Fano's Lemma} 

Several techniques have been developed to establish minimax lower bounds under the constraint of differential privacy.  Notable examples include the fingerprint method \citep{kamath19bridge},  Le Cam's method under differential privacy \citep{barber_duchi_2014lower}, differentially private Fano's lemma \citep{acharya_sun_zhang2021lower}, and the recently introduced Score Attack method \citep{cai2023score}.  Le Cam's method and Fano's lemma construct a multitude of hypotheses that are difficult to distinguish, while the fingerprint method and Score Attack design a test statistic with a prior distribution.

For our convenience, we employ the differentially private Fano's lemma, as detailed in Lemma~\ref{lem:dp-fano},  whose proof is provided in Section~\ref{sec:dp-fano} of the Appendix.  Here, ${\rm KL}(\cdot, \cdot)$ and ${\rm TV}(\cdot, \cdot)$ denote the Kullback-Leibler divergence and total variation distance between two distributions.

\begin{lemma}\label{lem:dp-fano} 
	Let $\mcP := \{P: P = \mu^{(1)} \times \cdots \times \mu^{(n)} \}$ be a family of product measures indexed by a parameter from a pseudo-metric space $(\Theta, \rho)$. Denote $\theta(P)\in\Theta$ the parameter associated with the distribution $P$.  Let $\mcQ=\{P_1,\cdots, P_N\}\subset \mcP$ contain $N$ probability measures  and there exist constants $\rho_0, l_0, t_0>0$ such that  for all $i\neq i^{\prime}\in[N]$, 
	$$
	\rho\left(\theta(P_i), \theta(P_{i^{\prime}})\right) \geqslant \rho_0, \quad \operatorname{KL}\left(P_{i} \| P_{i^{\prime}}\right) \leq l_0, 
	$$ 
	and 
	$$
	\sum_{k\in[n]} \operatorname{TV}\left(\mu_{i}^{(k)} , \mu_{i^{\prime}}^{(k)}\right) \leq t_0, 
	$$
	where $P_i = \mu_i^{(1)} \times \cdots \times \mu_i^{(n)} $ and $P_{i^{\prime}} = \mu_{i^{\prime}}^{(1)} \times \cdots \times \mu_{i^{\prime}}^{(n)} $. Then,
	\begin{align}\label{eq:fano-bd1}
		\inf_{A \in \mcA_{\varepsilon,\delta}(\mcP)} & \sup_{P\in\mcP } \EE_{A} \; \rho( A , \theta(P)) \geqslant \max \left\{\frac{\rho_0}{2}\left(1-\frac{ l_0 +\log 2}{\log N}\right), \frac{\rho_0}{4}\left(1 \bigwedge \frac{N-1}{\exp \left( 4 \varepsilon t_0 \right)}\right)\left(1-\frac{2\delta e^{4\varepsilon t_0}}{e^{\varepsilon}-1}\right)\right\},
	\end{align}
	where the infimum is taken over all the $(\varepsilon,\delta)$-DP randomized algorithm defined by $\mcA_{\varepsilon,\delta}(\mcP) := \{ A:X\mapsto \Theta \ {\rm and}\   A  \text{ is } \; (\varepsilon, \delta)\text{-differentially private} \; \text{for all} \; X\sim P \in\mcP \ \} $ . 
\end{lemma} 

Lemma~\ref{lem:dp-fano} provides a powerful tool for developing a minimax lower bound in estimation problems under the constraint of differential privacy.  Basically, if one can construct a sufficiently large set of distributions which are pairwise close in both Kullback-Leibler divergence and total variation distance, then a minimax lower bound can be derived if the underlying parameters are well-separated. 
The first term in the RHS of (\ref{eq:fano-bd1}) is derived from the classic Fano's Lemma without privacy constraint and serves as a lower bound for the statistical error rate. This term is a well-established outcome in information theory by the framework of hypothesis testing and has been extensively employed in the statistics literature. The second term in the RHS of (\ref{eq:fano-bd1}) characterizes the price one needs to pay for differential privacy. It is noteworthy that the cost of privacy is determined by $t_0$, which is the summation of marginal total variances. Intuitively, if the marginal total variance distances between $P_i = \mu_i^{(1)} \times \cdots \times \mu_i^{(n)} $ and $P_{i^{\prime}} = \mu_{i^{\prime}}^{(1)} \times \cdots \times \mu_{i^{\prime}}^{(n)} $ are small , it becomes challenging to identify the distribution from which the dataset is drawn. Therefore, the cost of privacy is expected to be low when $t_0$ is small. Moreover, if we assume that $X = (X_1, \cdots, X_n)\sim P_i$, then the cost of privacy resulting from replacing $X_k\sim \mu_i^{(k)}$ by $X_k^{\prime}\sim \mu_{i^{\prime}}^{(k)}$ should be upper bounded in  terms of $\operatorname{TV}(\mu_i^{(k)}, \mu_{i^{\prime}}^{(k)})$. We remark that the bound given in (\ref{eq:fano-bd1}) is meaningful only when $\delta\leq (e^{\eps}-1)e^{-4\eps t_0}$.

\subsection{Minimax lower bounds} 
\label{sec:lwb-PCA} 

In this section,  we apply Lemma~\ref{lem:dp-fano} to establish the minimax lower bounds for differentially private PCA and covariance estimation under the spiked covariance model.  Denote the family of normal distribution with a spiked covariance matrix by 
$$
\mcP(\lambda,\sigma^2):=\Big\{\mcN(0,  \Sigma): \Sigma=U\Lambda U^{\top}+\sigma^2 I_p\in \Theta(\lambda,\sigma^2)\Big\}. 
$$
By definition,  each distribution $P\in\mcP(\lambda,\sigma^2)$ is indexed by the pair of eigenvalues $\Lambda$ and eigenvectors $U\in\OO_{p, r}$.  We first focus on the minimax lower bounds for estimating the spectral projector $UU^{\top}$.   Similarly, the minimax lower bounds are established in all Schatten-$q$ norms for $q\in[1,\infty]$. 

\begin{theorem}\label{thm:dp-pca-lwb}
	Let the $p\times n$ data matrix $X$ have i.i.d.  columns sampled from a distribution $P=\mcN(0,  U^{\top}\Lambda U^{\top}+\sigma^2 I_p)\in\mcP(\lambda,\sigma^2)$. Suppose $\delta\leq c_0'\exp\big\{2\eps-c_0\big(\eps\sqrt{npr}+pr\big)\big\}$ for some small constants $c_0,c_0'>0$.  Then,  there exists an absolute constant $c_1>0$ such that 
	\begin{align*}
		\inf_{\wtU\in\mcU_{\varepsilon, \delta}} \sup_{P \in\mcP(\lambda,\sigma^2) }& \frac{\EE \|\tdU\tdU^{\top} - UU^{\top} \|_q}{r^{1/q}} \geq c_1
		\Bigg(  \bigg(\frac{\sigma^2}{\lambda}+\sqrt{\frac{\sigma^2}{\lambda}}\bigg) \left(\sqrt{\frac{p}{n}} +\frac{p\sqrt{r}}{n\epsilon}\right) \Bigg) \bigwedge 1, 
	\end{align*}
	where the infimum is taken over all the possible $(\varepsilon,\delta)$-DP algorithms, denoted by $\mcU_{\varepsilon, \delta}$, and the expectation is taken with respect to both $\wtU$ and $P$. 
\end{theorem}

Theorem~\ref{thm:dp-pca-lwb} imposes a strong restriction on the parameter $\delta$. For most interesting cases, $\delta$ needs to be near zero. Therefore, the minimax lower bounds hold primarily for the pure differential privacy case, i.e., $\delta=0$. 
It is worth noting the two terms in the minimax lower bound of spectral norm ($q=\infty$):
\begin{equation}\label{eq:lwb-pca-cmp}
	\bigg(\frac{\sigma^2}{\lambda}+\sqrt{\frac{\sigma^2}{\lambda}}\bigg)\sqrt{\frac{p}{n}}\quad {\rm and}\quad \bigg(\frac{\sigma^2}{\lambda}+\sqrt{\frac{\sigma^2}{\lambda}}\bigg)\frac{p\sqrt{r}}{n\varepsilon}.
\end{equation}
The first term concerns the statistical error of PCA without privacy constraint. The error bound is free of the rank $r$, which is very typical in spectral norm error rate and the rate matches the existing minimax optimal rate of PCA for spiked covariance model. See, e.g., \cite{cai2015optimal, zhang2018tensor, yu2015useful}. The second term is the price paid for differential privacy. Interestingly, the second term is dependent on the rank $r$ even though spectral norm is considered here. The technical explanation is that the sensitivity of empirical spectral projector increases as the number of PC's grows. Comparing the two terms in (\ref{eq:lwb-pca-cmp}), we observe that if $\varepsilon\geq (rp/n)^{1/2}$, the cost of privacy is dominated by the statistical error. 

A minimax lower bound for \emph{rank-one} PCA has been established in \cite[Theorem~5.3]{liu_xiyang2022dp-pca}. Their developed rate in spectral norm also have two terms:
$$
\sqrt{\frac{\sigma^2}{\lambda+\sigma^2}}\cdot \sqrt{\frac{p}{n}}\quad {\rm and}\quad \sqrt{\frac{\sigma^2}{\lambda+\sigma^2}}\cdot \frac{p}{n\varepsilon}. 
$$
Their rate matches ours when $r=1$ and $\lambda\geq \sigma^2$.  On the other hand, if $\lambda\ll \sigma^2$, our minimax lower bound is much stronger. Moreover, our minimax lower bounds hold for a diverging rank as long as $2r\leq p$.

We now shift our focus to the minimax lower bound of differentially private estimation of the spiked covariance matrix. Here, we assume $\sigma^2$ is known and it suffices to estimate the signal part $U\Lambda U^{\top}$. As a result, the minimax lower bound is essentially determined jointly by the lower bounds in estimating eigenvalues and eigenvectors. 

\begin{theorem}\label{thm:dp-Sigma-lwb}
	Let the $p\times n$ data matrix $X$ have i.i.d.  columns sampled from a distribution $P=\mcN(0,  U^{\top}\Lambda U^{\top}+\sigma^2 I_p)\in\mcP(\lambda,\sigma^2)$. Suppose $\delta\leq c_0'\exp\big\{2\eps-c_0\big(\eps\sqrt{npr}+pr\big)\big\}$ for some small constants $c_0,c_0'>0$.  Then,  there exists an absolute constant $c_1>0$ such that 
	\begin{align*}
	\inf_{\wtSig\ \in\mcM_{\varepsilon, \delta}} \sup_{P\in\mcP(\lambda,\sigma^2)} \frac{\EE \big\| \wtSig - \Sigma \big\|_q}{r^{1/q}} \geq c_1
	\left( \lambda \left(\sqrt{\frac{r}{n}} + \frac{r^{3/2}}{n\varepsilon}\right) + \sqrt{\sigma^2(\lambda+\sigma^2)} \left(\sqrt{\frac{p}{n}} + \frac{\sqrt{r}p}{n \varepsilon} \right) \right)\bigwedge \lambda, 
\end{align*} 
	where the infimum is taken over all the possible $(\varepsilon,\delta)$-DP algorithms, denoted by $\mcM_{\varepsilon, \delta}$, and the expectation is taken with respect to both $\wtSig$ and $P$.  Here, $q$ can be any number in $[1,\infty]$.  
\end{theorem}

Without loss of generality, let us discuss the two terms in the spectral norm distance
\begin{equation}\label{eq:lwb-Sigma-cmp}
	\lambda\bigg(\sqrt{\frac{r}{n}}+\frac{r^{3/2}}{n\varepsilon}\bigg)\quad {\rm and}\quad \sqrt{\sigma^2(\lambda+\sigma^2)}\bigg(\sqrt{\frac{p}{n}}+\frac{\sqrt{r}p}{n\varepsilon}\bigg). 
\end{equation} 
The second term is contributed by the differentially private estimation error of PCA in the form of $\lambda\|\wtU\wtU^{\top}-UU^{\top}\|_{\rm F}^2$. 
The first term dominates if the signal strength is exceedingly large, or more precisely, when $\lambda/\sigma^2\gg p/r$. In this case, we can simply regard $\sigma=0$ and the stochastic error mainly comes from the randomness of a low-dimensional distribution. Basically, it suffices to consider the minimax optimal estimation under a smaller family of normal distributions $\{\mcN(0, \lambda UU^{\top}+\lambda I_r):\ U\in\OO_{r, r/4}\}$. By replacing $\sigma\leftarrow \lambda$, $r\leftarrow r/4$, and $p\leftarrow r$, the second term reduces to the first term in (\ref{eq:lwb-Sigma-cmp}). 
Without the privacy constraint, the first term also matches the existing optimal rate in covariance estimation under spiked covariance model \citep{cai2015optimal, cai2010optimal}.

\section{Extensions}
\label{sec:discuss}

For the sake of clarity, we have assumed uniformity in the order of spiked eigenvalues and Gaussian distributions. In this section, we extend our analysis to provide upper bounds for differentially private PCA and covariance estimation without requiring these specific conditions. 

\subsection{Diverging condition number}

Suppose that $X_1,\cdots, X_n\stackrel{{\rm i.i.d.}}{\sim} \mcN(0, \Sigma)$ where $\Sigma=U\Lambda U^{\top}+\sigma^2 I_p$ with spiked eigenvalues $\Lambda={\rm diag}(\lambda_1,\cdots,\lambda_r)$. Denote $\kappa_0:=\lambda_1/\lambda_r$, the ratio of the largest and smallest spiked eigenvalues. The  proof of Corollary~\ref{cor:kappa-dvg} is almost identical to that of Theorems~\ref{thm:dp-pca-upb} and ~\ref{thm:dp-Sigma-upb}, and thus omitted. We  only present the upper bounds of the expected error in Schatten norms, but high probability bounds hold similarly. 

\begin{corollary} \label{cor:kappa-dvg}
	Suppose that $n\geq C_1(\kappa_0^2r\log n+\log^2n),  2r\leq p$, $\lambda_1/\sigma^2\leq (p/n)e^{(p\wedge n)/C_1}$, and $\lambda_r/\sigma^2\geq C_1(\kappa_0p/n+\sqrt{p/n})$ for some large absolute constant $C_1>0$.  If we choose 
	$$
	\Delta_1:=C_2\bigg(\frac{\sigma^2}{\lambda_r}+\sqrt{\frac{\kappa_0\sigma^2}{\lambda_r}}\bigg)\frac{\sqrt{p(r+\log n)}}{n}\quad {\rm and}\quad \Delta_2:=C_2\frac{\lambda_1(r+\log n)+\sigma^2 (p+\log n)}{n}, 
	$$
	then, there exist absolute constants $C_4>0$ such that,  for any $\varepsilon>0, \delta\in(0,1)$,  Algorithm~\ref{alg:DP-PCA} outputs an $(\varepsilon,\delta)$-DP estimators $\wtU\wtU^{\top}$ and $\wtSig$ satisfying 
	\begin{align*}
		\frac{\EE \|\wtU\wtU^{\top}-UU^{\top}\|_q}{r^{1/q}} \leq C_4\bigg(\frac{\sigma^2}{\lambda_r}+\sqrt{\frac{\kappa_0\sigma^2}{\lambda_r}}\bigg)\bigg(\sqrt{\frac{p}{n}}+\frac{p\sqrt{r+\log n}}{n\varepsilon}\log^{1/2}\Big(\frac{2.5}{\delta}\Big)\bigg).
	\end{align*} 
	and
	\begin{align*}
		&\frac{\EE \|\wtSig-\Sigma\|_q}{r^{1/q}} \\
		&\quad\leq  C_4\Bigg(\lambda_1\bigg(\sqrt{\frac{r}{n}}+\frac{\sqrt{r}(r+\log n)}{n\varepsilon}\cdot \sqrt{\log\frac{2.5}{\delta}}\bigg)+\sqrt{\sigma^2(\lambda_1+\sigma^2)} \bigg(\sqrt{\frac{p}{n}}+\frac{p\sqrt{(r+\log n)}}{n\varepsilon}\sqrt{\log\frac{2.5}{\delta}}\bigg)\Bigg). 
	\end{align*} 
	for all $q\in[1,\infty]$. 
\end{corollary} 

\subsection{Sub-Gaussian}\label{sec:sub-Gaussian}
Suppose that $X$ follows a sub-Gaussian distribution satisfying that, for any $u\in\RR^p$,  the following bound holds
$$
\EE \exp\bigg\{\frac{\langle X,  u\rangle^2}{u^{\top}\Sigma u}\bigg\}\leq 2,
$$
where $\Sigma\in\Theta(\lambda,\sigma^2)$.  For ease of exposition,  we focus on the case of bounded condition number.  Interestingly,  the sensitivity of eigenvectors and eigenvalues is actually identical to that under Gaussian distributions.  

\begin{corollary} \label{cor:sub-Gaussian}
	Suppose that $n\geq C_1\big(r\log (p+n)\log^2 r+\log^2n\big),  2r\leq p$, $\lambda/\sigma^2\leq (p/n)e^{(p\wedge n)/C_1}$, and $\lambda/\sigma^2\geq C_1(p/n+\sqrt{p/n})\log(p+n)$ for some large absolute constant $C_1>0$.  If we choose 
	$$
	\Delta_1:=C_2\bigg(\frac{\sigma^2}{\lambda}+\sqrt{\frac{\sigma^2}{\lambda}}\bigg)\frac{\sqrt{p(r+\log n)}}{n}\quad {\rm and}\quad \Delta_2:=C_2\frac{\lambda(r+\log n)+\sigma^2 (p+\log n)}{n}, 
	$$
	then,  there exist absolute constant $C_4>0$ such that,  for any $\varepsilon>0, \delta\in(0,1)$,  Algorithm~\ref{alg:DP-PCA} outputs an $(\varepsilon,\delta)$-DP estimators $\wtU\wtU^{\top}$ and $\wtSig$ satisfying 
	\begin{align*}
		\frac{\EE \|\wtU\wtU^{\top}-UU^{\top}\|_q}{r^{1/q}} \leq C_4\bigg(\frac{\sigma^2}{\lambda}+\sqrt{\frac{\sigma^2}{\lambda}}\bigg)\bigg(\sqrt{\frac{p\log p}{n}}+\frac{p\sqrt{r+\log n}}{n\varepsilon}\log^{1/2}\Big(\frac{2.5}{\delta}\Big)\bigg), 
	\end{align*} 
	and
	\begin{align*}
		&\frac{\EE \|\wtSig-\Sigma\|_q}{r^{1/q}} \\
		&\quad\leq  C_4\Bigg(\lambda\bigg(\sqrt{\frac{r}{n}}+\frac{\sqrt{r}(r+\log n)}{n\varepsilon}\cdot \sqrt{\log\frac{2.5}{\delta}}\bigg)+\sqrt{\sigma^2(\lambda+\sigma^2)} \bigg(\sqrt{\frac{p\log p}{n}}+\frac{p\sqrt{r+\log n}}{n\varepsilon}\sqrt{\log\frac{2.5}{\delta}}\bigg)\Bigg). 
	\end{align*} 
	for all $q\in[1,\infty]$. 
\end{corollary} 

As shown by Corollary~\ref{cor:sub-Gaussian},  the upper bounds of differentially private sub-Gaussian PCA and covariance estimation are almost the same as those for Gaussian distributions,  implying that these bounds are minimax optimal.  However,  some additional logarithmic factors appear in the upper bound and signal-to-noise ratio condition when controlling the higher-order terms in spectral perturbation.  


\subsection{Private estimation of nuisance variance}

In this section, we provide a differentially private estimator of $\sigma^2$, demonstrating that minimax optimal estimation of the spiked covariance matrix is still achievable even if $\sigma^2$ is unknown.

The estimation of $\sigma^2$ in spiked covariance models has been studied by \cite{donoho2018optimal, cai2015optimal, shabalin2013reconstruction}. These methods exploit the eigenvalues of the sample covariance matrix or the properties of the empirical spectral distribution. 
We utilize the robust estimator of $\sigma^2$, originally proposed by \cite{ke2023estimation} for detecting the number of spikes. The basic idea is to average several bulk eigenvalues that are separated from the spike eigenvalues. We begin by reviewing the well-known Marchenko-Pastur (MP) law \citep{marchenko1967distribution, bai2010spectral}.

Let $Z$ be a $p\times n$ matrix whose entries are independent, identically distributed random variables with mean $0$ and variance $\sigma^2$. Let $Y_n:=ZZ^{\top}/n$ be the sample covariance matrix, and let $\lambda_1(Y_n)\geq\cdots\geq\lambda_{p\wedge n}(Y_n)$ be its non-zero eigenvalues. Define the empirical spectral distribution (ESD) of $Y_n$ by 
$
\mu_n(A):=(p\wedge n)^{-1}\sum_{j=1}^{p\wedge n}\mathbbm{1}\big(\lambda_j(Y_n)\in A\big), \ \forall A\subset \RR.
$
Assume that $p/n\to \gamma$ as $n\to\infty$ for some $\gamma>0$. It is known that the ESD converges in distribution to the MP distribution $\mu(\cdot)$, whose density function is given as follows. 
\begin{definition}
Given $\gamma>0$, the zero-excluded MP distribution is defined by the density 
$$
f_{\gamma,\sigma^2}(x)=\frac{1}{2\pi\sigma^2} \cdot \frac{1}{x(1\wedge \gamma)}\sqrt{(x-\sigma^2\gamma_{-})(\sigma^2\gamma_+-x)}\cdot\mathbbm{1}\big(x\in[\sigma^2\gamma_-,\ \sigma^2\gamma_+]\big),
$$
where $\gamma_{\pm}:=(1\pm \sqrt{\gamma})^2$. 
\end{definition}
Since $\EE Y_n=\sigma^2 I_p$, a natural estimator of $\sigma^2$ is by taking the average of several bulk eigenvalues of $Y_n$. Recall the sample covariance matrix $\widehat\Sigma$ under the spiked model (\ref{eq:spiked-cov}), and let $\lambda_1(\widehat\Sigma)\geq \cdots\geq \lambda_{p\wedge n}(\widehat\Sigma)$ denote the non-zero eigenvalues of $\widehat\Sigma$. The eigenvalue sticking property \citep[Theorem 2.7]{bloemendal2016principal} tells that $\lambda_{j+r}(\widehat\Sigma)\approx \lambda_j(Y_n)$ with high probability for all $(p\wedge n)/4 \leq j\leq 3(p\wedge n)/4$ if $r\ll (p\wedge n)$.  Denote $q_k$ the $k/(p\wedge n)$-upper quantile of the MP distribution with $\gamma_n:=p/n$ and $\sigma^2=1$, i.e., 
$
\int_{q_k}^{(1+\sqrt{\gamma_n})^2} f_{\gamma_n, 1}(x)dx=k/(p\wedge n).  
$

We define the non-private estimator of $\sigma^2$ by 
$$
\widehat\sigma^2:=\frac{\sum_{(p\wedge n)/4\leq k\leq 3(p\wedge n)/4}q_k \lambda_k(\widehat\Sigma)}{\sum_{(p\wedge n)/4\leq k\leq 3(p\wedge n)/4}q_k^2}. 
$$
The convergence rate of $\widehat\sigma^2$ was established by Theorem~1 of \cite{ke2023estimation}. Its sensitivity is characterized in Theorem~\ref{thm:unknown-sigma}. 
Let $\widehat\sigma^{(i)2}$ be defined as $\widehat\sigma^2$ using $\widehat\Sigma^{(i)}$ instead of $\widehat\Sigma$. 

\begin{theorem}\label{thm:unknown-sigma}
Suppose the conditions in Theorem~\ref{thm:dp-Sigma-upb} hold, $r\leq C_3$ for any large constant $C_3$,  and $p/n\to \gamma$ for a constant $\gamma>0$ .  Then, there exists an absolute constant $C_4>0$ such that, with probability at least $1-n^{-100}$,  
$$
\big|\widehat\sigma^2-\widehat\sigma^{(i)2}\big|\leq \Delta_3:=\frac{C_4}{\sqrt{p\wedge n}}\cdot \frac{\lambda(r+\log n)+\sigma^2(p+\log n)}{n}. 
$$
Consequently, the estimator $\widetilde\sigma^2:=\big|\widehat \sigma^2+\mcN\big(0, 18(\Delta_3/\eps)^2\log(3.75/\delta)\big)\big|$ is an $(\eps/3, \delta/3)$-DP estimate of $\sigma^2$ with probability at least $1-n^{-100}$. 
\end{theorem}

We now study the DP PCA and covariance estimator using the private estimate $\widetilde\sigma^2$. 

\begin{theorem}\label{thm:unknown-sigma-2}
Suppose the conditions in Theorem~\ref{thm:unknown-sigma},  $p/n\to \gamma$ for a constant $\gamma>0$,  and there exists a small constant $c_1>0$ such that 
\begin{align}\label{eq:unknown-sigma-cond}
\frac{\lambda}{\sigma^2}\leq c_1\frac{n\sqrt{(p\wedge n)/\log n}}{(r+\log n)}\cdot \frac{\eps}{\sqrt{\log(4/\delta)}}\quad {\rm and}\quad \frac{p}{n}\sqrt{\frac{\log n}{p\wedge n}}\cdot \frac{\eps}{\sqrt{\log(4/\delta)}}\leq c_1
\end{align}
 Let $\widetilde U, \widetilde\Lambda$, and $\widetilde \Sigma$ be defined as Algorithm~\ref{alg:DP-PCA} with replacing $(\eps, \delta, \sigma^2)$ by $(\eps/3, \delta/3, \widetilde\sigma^2)$,  respectively. Then, with probability at least $1-e^{-c_2(n\wedge p)}-4n^{-99}$,  $\widetilde U$ and $\widetilde\Lambda$ are $(\eps/3, \delta/3)$-DP,  $\widetilde\Sigma$ is $(\eps,\delta)$-DP, and the following bounds hold
 \begin{align*}
		\|\wtU\wtU^{\top}-UU^{\top}\|_{\rm F}\leq C_4\bigg(\frac{\sigma^2}{\lambda}+\sqrt{\frac{\sigma^2}{\lambda}}\bigg)\bigg(\sqrt{\frac{pr}{n}}+\frac{p\sqrt{r(r+\log n)}}{n\varepsilon}\sqrt{\log\frac{4}{\delta}}\bigg),
	\end{align*} 
	and
\begin{align*}
&\max\Big\{\big\|\widetilde\Sigma-\Sigma\big\|, \big\|\widetilde\Sigma-\Sigma\big\|_{\rm F}\Big\}\\
&\leq C_{5,\gamma} \Bigg(\lambda\bigg(\sqrt{\frac{r}{n}}+\frac{(r+\log n)^{3/2}}{n\varepsilon}\cdot \sqrt{\log\frac{4}{\delta}}\bigg)+\sqrt{\sigma^2(\lambda+\sigma^2)} \bigg(\sqrt{\frac{p}{n}}+\frac{p\sqrt{(r+\log n)}}{n\varepsilon}\sqrt{\log\frac{4}{\delta}}\bigg)\Bigg),
\end{align*}
where $c_2, C_{5,\gamma}>0$ are constants. 
\end{theorem}
These rates are nearly identical to those in Theorems~\ref{thm:dp-pca-upb} and \ref{thm:dp-Sigma-upb}, making them minimax optimal up to logarithmic factors. Condition (\ref{eq:unknown-sigma-cond}) is imposed to ensure that the outputs of Algorithm~\ref{alg:DP-PCA} remain differentially private when using the private estimate $\widetilde{\sigma}^2$ instead of the true $\sigma^2$. Essentially, this condition guarantees that $\widetilde{\sigma}^2 \asymp \sigma^2$ with high probability, ensuring that the artificial noise added in Algorithm~\ref{alg:DP-PCA} is sufficiently strong to maintain the privacy guarantee.

\section{Numerical Experiments}\label{sec:numeric}

\subsection{Simulations}

We present simulation results comparing the performance of our differentially private (DP) algorithms with existing methods in the literature. The \textsf{DP-Oja} algorithm, proposed by \cite{liu_xiyang2022dp-pca}, estimates the first principal component under privacy constraints by extending Oja's algorithm, originally introduced by \cite{oja1982simplified} for online PCA. We also compare it with the \textsf{DP-Gauss} algorithm proposed by \cite{dwork2014analyze}, which uses the Gaussian mechanism for privacy-preserving PCA. Both our method and \textsf{DP-Gauss} use the Gaussian mechanism to ensure privacy. However, our method is specifically motivated by the spiked covariance model, while \textsf{DP-Gauss} is designed for deterministic data, assuming each observation has at most unit norm. This distinction is also reflected in the artificial noise levels introduced by the two methods. For the spiked covariance model, we implement the \textsf{DP-Gauss} method after applying a universal scaling to ensure each observation has at most unit norm. The \textsf{DP-Oja} method, as discussed by \cite{liu_xiyang2022dp-pca}, is directly applicable to the spiked covariance model for Gaussian distributions.

Under the spiked covariance model (\ref{eq:spiked-cov}), we carefully and fairly choose the constants involved in the artificial noise level for the \textsf{DP-Oja}, \textsf{DP-Gauss}, and our methods. For simplicity, we set the nuisance variance parameter $\sigma^2=1$ in all experiments. 
 In our method, the constant $C_2$  in the definition of $\Delta_1$ in Lemma~\ref{lem:DP-alg} arises from the concentration of the sub-exponential random variable $\max_{i\in[n]}\|U^{\top}X_iX_i^{\top}U_{\perp}\|$ (see Lemma~\ref{lem:con-bounds}). In \textsf{DP-Gauss}, a scaling by $\big(\max_{i\in[n]}\|X_i\|^2\big)^{-1}$ is applied to the sample covariance matrix to enforce the unit norm upper bound assumption. By Lemma~\ref{lem:con-X-norm}, $\max_{i\in[n]}\|X_i\|^2\leq (r+C_5\log n)\lambda+p\sigma^2$. The scaling factor is set as $(r+C_5\log n)\lambda+p\sigma^2$ in \textsf{DP-Gauss}. For a fair comparison, we set both the $C_2$ in our method and the $C_5$ in the scaling factor in \textsf{DP-Gauss} to $4$. We can also set the scaling factor in \textsf{DP-Gauss} exactly as $\big(\max_{i\in[n]}\|X_i\|^2\big)^{-1}$. The resultant algorithm is denoted as \textsf{DP-Gauss*}. 
  The artificial noise level in \textsf{DP-Oja} involves several unsettled constant factors, with the dominating one determined by the concentration property of $\max_{i\in[n]}\|X_iX_i^{\top}\hat u_{i-1}\|$ (see Lemma~3.2 in \cite{liu_xiyang2022dp-pca}). For a fair comparison, we should set the constant $C=2$ in their Algorithm 2. However, since \textsf{DP-Oja} is an online PCA algorithm with random initialization, it is unsurprising that this method significantly underperforms our method and $\textsf{DP-Gauss}$. For clear illustration, we set $C=0.2$ in their Algorithm 2. Note that \textsf{DP-Oja} only works for rank-one PCA, and its stepsize is set as $0.5/n$ after fine tuning for the best performance.  
 
In the first simulation setting, we set $p=50$, $\Sigma=\lambda UU^{\top}+I_p$ with $\lambda=10$, and $U$ contains the left singular vectors of a $p\times r$ random matrix with i.i.d. $N(0,1)$ entries.  The privacy budget is set as $\eps=1$ and $\delta=0.1$. We first set $r=1$ and compare the utility performances of the three methods in terms of the error $\|\widehat U\widehat U^{\top}-UU^{\top}\|_{\rm F}$ as the sample size $n$ varies. For each $n$, the simulation is repeated for 40 times, and the average error  and standard deviation are recorded. The results are displayed in the left panel of Figure~\ref{fig:S1}. It shows that \textsf{DP-Oja} significantly underperforms compared to other methods, while our method slightly outperforms \textsf{DP-Gauss}. Our method and \textsf{DP-Gauss*} achieve similar performance. Moreover, we observe that \textsf{DP-Oja} runs much slower than the others. 
We then compare our method and \textsf{DP-Gauss} while varying the rank $r$. The boxplots shown in the right panel of Figure~\ref{fig:S1} are  based on 40 independent simulations, demonstrating that our method consistently outperforms \textsf{DP-Gauss} for different ranks. Interestingly, as the rank $r$ increases, our method becomes better than \textsf{DP-Gauss*}.

\begin{figure}
	\centering
	\begin{subfigure}[b]{.48\linewidth}
		\includegraphics[width=\linewidth]{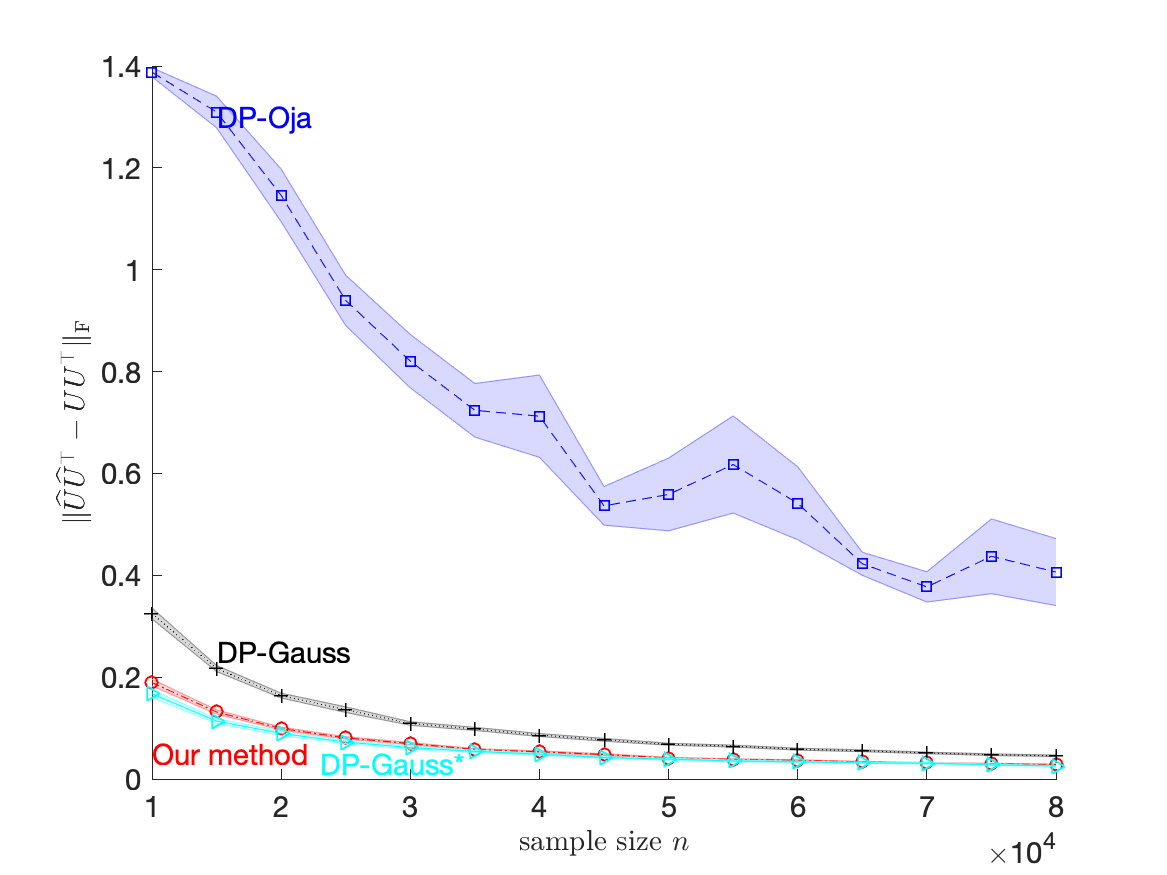}
		\caption{\tiny{Fixed rank $r=1$. }}
	\end{subfigure}
	\begin{subfigure}[b]{.48\linewidth}
		\includegraphics[width=\linewidth]{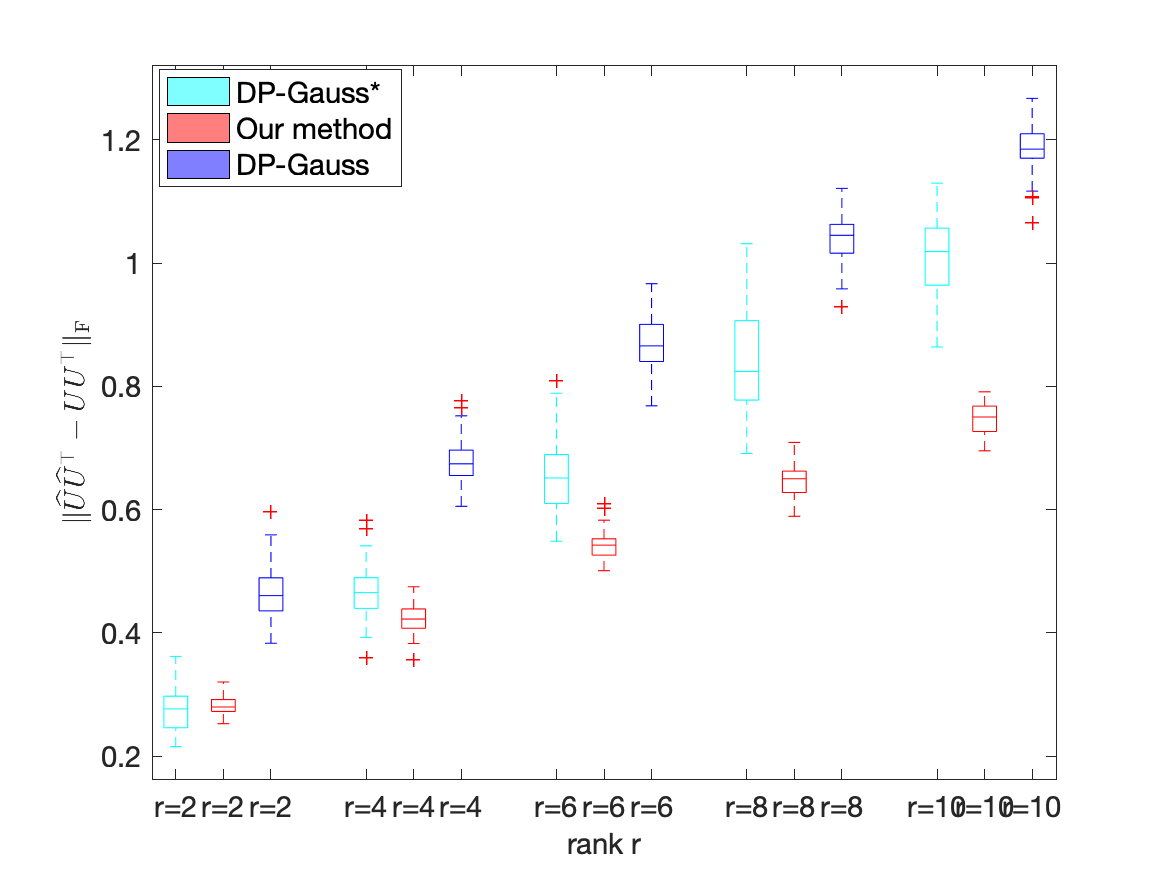}
		\caption{\tiny{Fixed sample size $n=10^4$.}}
	\end{subfigure}
	\caption{Comparison of our method, \textsf{DP-Oja} \cite{liu_xiyang2022dp-pca}, and \textsf{DP-Gauss}, \textsf{DP-Gauss*} \cite{dwork2014analyze} in differentially private PCA with varying $n$ and $r$.  The dimension $p=50$, $\lambda=10, \sigma^2=1$, and privacy constraints $\eps=1, \delta=0.1$. 
	}
	\label{fig:S1}
\end{figure}

The second simulation setting compares these methods with respect to the varying privacy parameter $\eps$ and signal strength $\lambda$. The results, presented in Figure \ref{fig:S2}, show that the error rates of all methods increase rapidly as $\eps$ decreases, which aligns with our theoretical predictions. When the signal strength $\lambda$ is small, \textsf{DP-Gauss*} performs the best. However, as $\lambda$ increases, our method outperforms the others.

\begin{figure}
	\centering
	\begin{subfigure}[b]{.48\linewidth}
		\includegraphics[width=\linewidth]{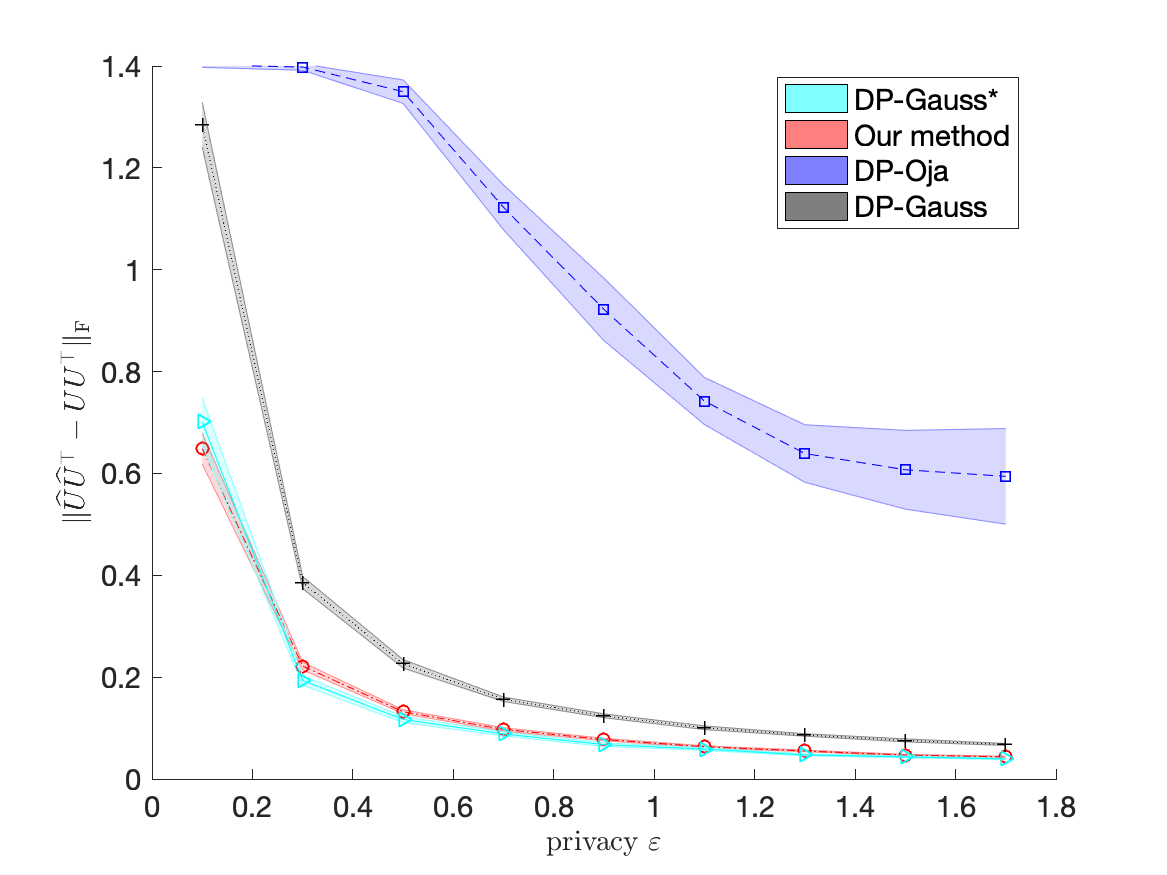}
		\caption{\tiny{Fixed rank $r=1$, $\lambda=10$, and sample size $n=3\times 10^4$. }}
	\end{subfigure}
	\begin{subfigure}[b]{.48\linewidth}
		\includegraphics[width=\linewidth]{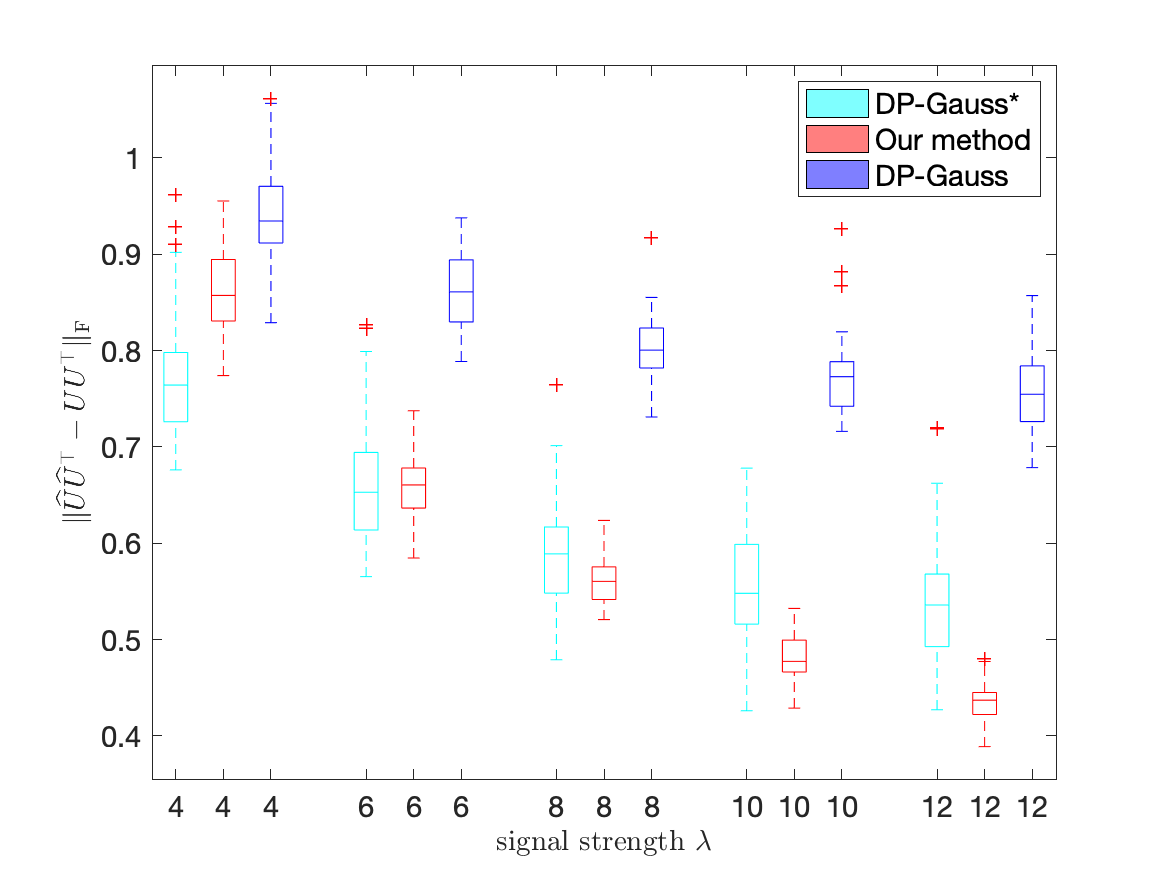}
		\caption{\tiny{Fixed rank $r=5$, $\eps=1$, and sample size $n=10^4$.}}
	\end{subfigure}
	\caption{Comparison of our method, \textsf{DP-Oja} \cite{liu_xiyang2022dp-pca}, and \textsf{DP-Gauss}, \textsf{DP-Gauss*} \cite{dwork2014analyze} in differentially private PCA with varying $\eps$ and $\lambda$.  The dimension $p=50$,  $\sigma^2=1$, and privacy constraint $\delta=0.1$. 
	}
	\label{fig:S2}
\end{figure}

The third simulation setting aims to test the performance of our method and others in the high-dimensional case where $p\geq n$. We demonstrate that our method is applicable as long as the signal strength condition $\lambda/\sigma^2\geq  C(p/n+\sqrt{p/n})$ hold.  The dimension, rank, and sample size are set to $p=50$, $r=3$, and $n=30$, respectively.  
The boxplots of error $\|\widehat U\widehat U^{\top}-UU^{\top}\|_{\rm F}$ based on $40$ simulations are displayed in Figure~\ref{fig:S3}. They show that our method significantly outperforms \textsf{DP-Gauss} and \textsf{DP-Gauss*} in this setting. The error rates of \textsf{DP-Gauss} and \textsf{DP-Gauss*} hardly improve as the signal strength increases.  One possible reason is that the universal scaling procedure used in \textsf{DP-Gauss} and \textsf{DP-Gauss*} leads to significant information loss, especially when sample size $n$ is small.

\begin{figure}
	\centering
	\begin{subfigure}[b]{.48\linewidth}
		\includegraphics[width=\linewidth]{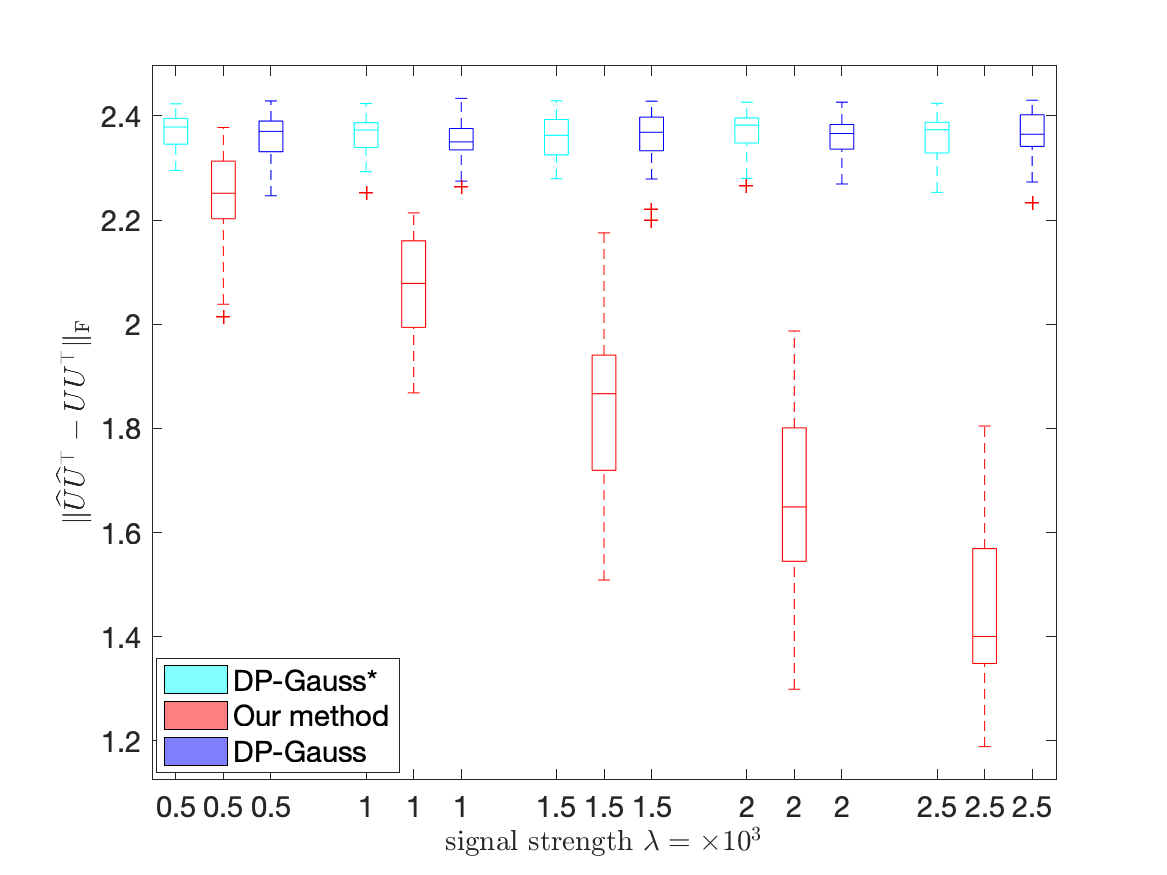}
		\caption{\tiny{Privacy constraint $\eps=3$. }}
	\end{subfigure}
	\begin{subfigure}[b]{.48\linewidth}
		\includegraphics[width=\linewidth]{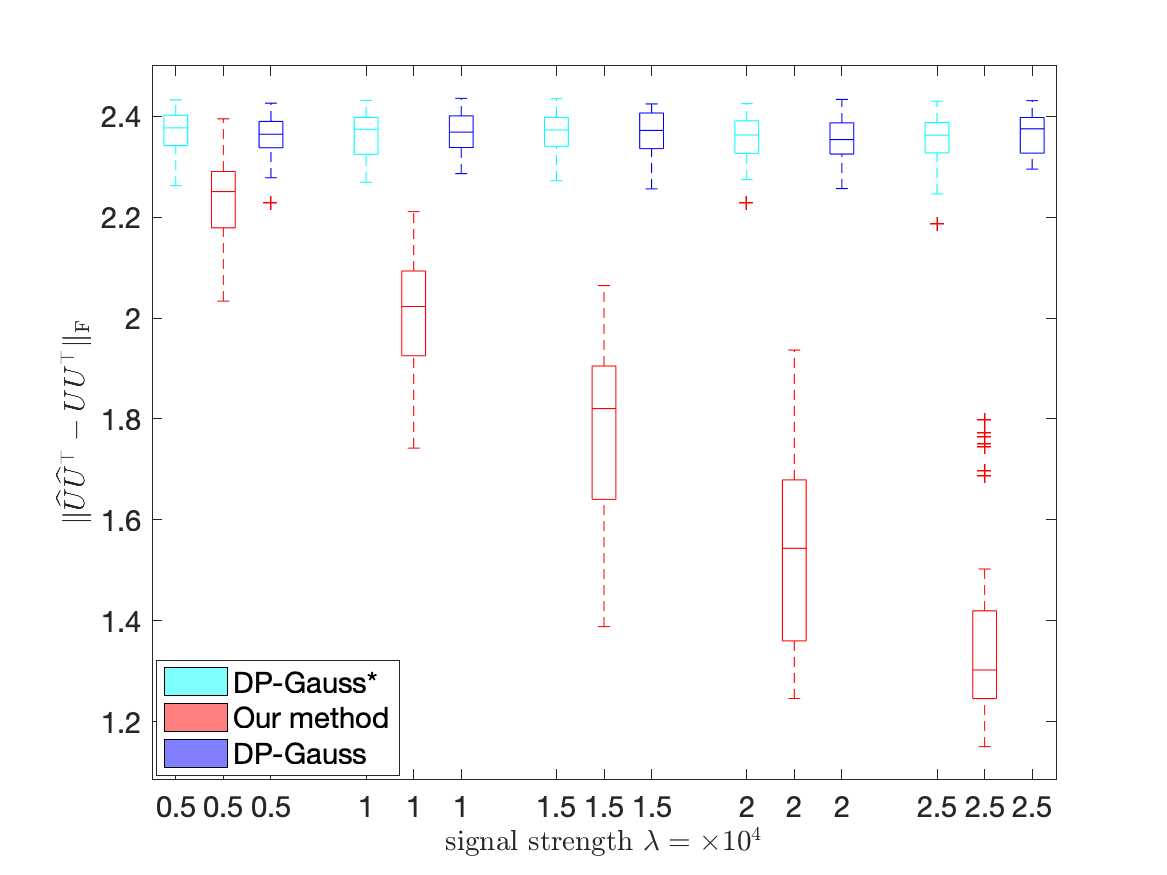}
		\caption{\tiny{Privacy constraint $\eps=1$.}}
	\end{subfigure}
	\caption{Comparison of our method, \textsf{DP-Gauss}, and \textsf{DP-Gauss*} \cite{dwork2014analyze} in differentially private PCA when $p\geq n$ and the signal strength $\lambda$ changes.  The dimension $p=50$, $n=30, r=3, \sigma^2=1$, and privacy constraints $\delta=0.1$. 
	}
	\label{fig:S3}
\end{figure}

In the fourth simulation settings, we compare the performance of our method, \textsf{DP-Gauss}, and \textsf{DP-Gauss*} for  covariance estimation.   Due to space constraints, the results are provided in Appendix~\ref{sec:app-simulation}.


\subsection{MNIST dataset}

We implemented our method, \textsf{DP-Gauss}, and \textsf{DP-Gauss*}  for differentially private PCA on the MNIST dataset, which contains grayscale images of handwritten digits from 0 to 9. The dataset includes 500 samples for each digit.} 
For a clear illustration, we used only the images corresponding to digits 1, 4, and 9, creating a sample of $n=1,500$ images. Each image is an observation of length $p=28\times 28=784$. The dimension is relatively large compared to the sample size.  We downscaled the original images to  $14\times 14$. The rank was set to $r=3$. We estimated $\lambda$ by averaging the first three eigenvectors of the sample covariance matrix and $\sigma^2$ by the mean the $50$-th to $140$-th sample eigenvalues. The privacy constraints were set as $\eps=2$ and $\delta=0.1$. After obtaining the DP estimates of eigenvectors, we applied dimension reduction to each observation, reported the variance explained, and visualized the scores corresponding to first and second components. The results are shown as in Figure~\ref{fig:MNIST}. As expected, the \textsf{DP-Gauss*} method performed poorly due to the relatively small sample size. These methods add too much noise, resulting in principal components that can explain only 2\% of the total variance.

\begin{figure}
	\centering
	\begin{subfigure}[b]{.48\linewidth}
		\includegraphics[width=\linewidth]{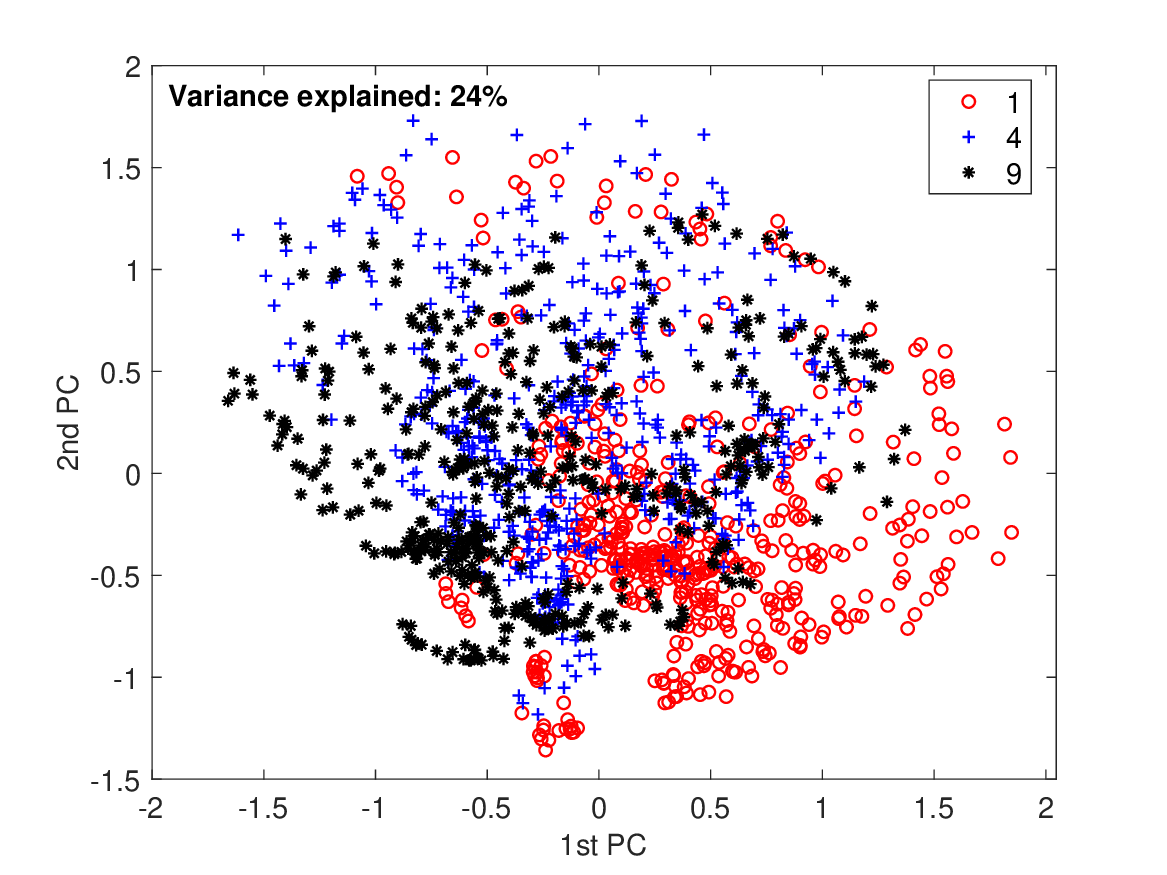}
		\caption{\tiny{Our method }}
	\end{subfigure}
	\begin{subfigure}[b]{.48\linewidth}
		\includegraphics[width=\linewidth]{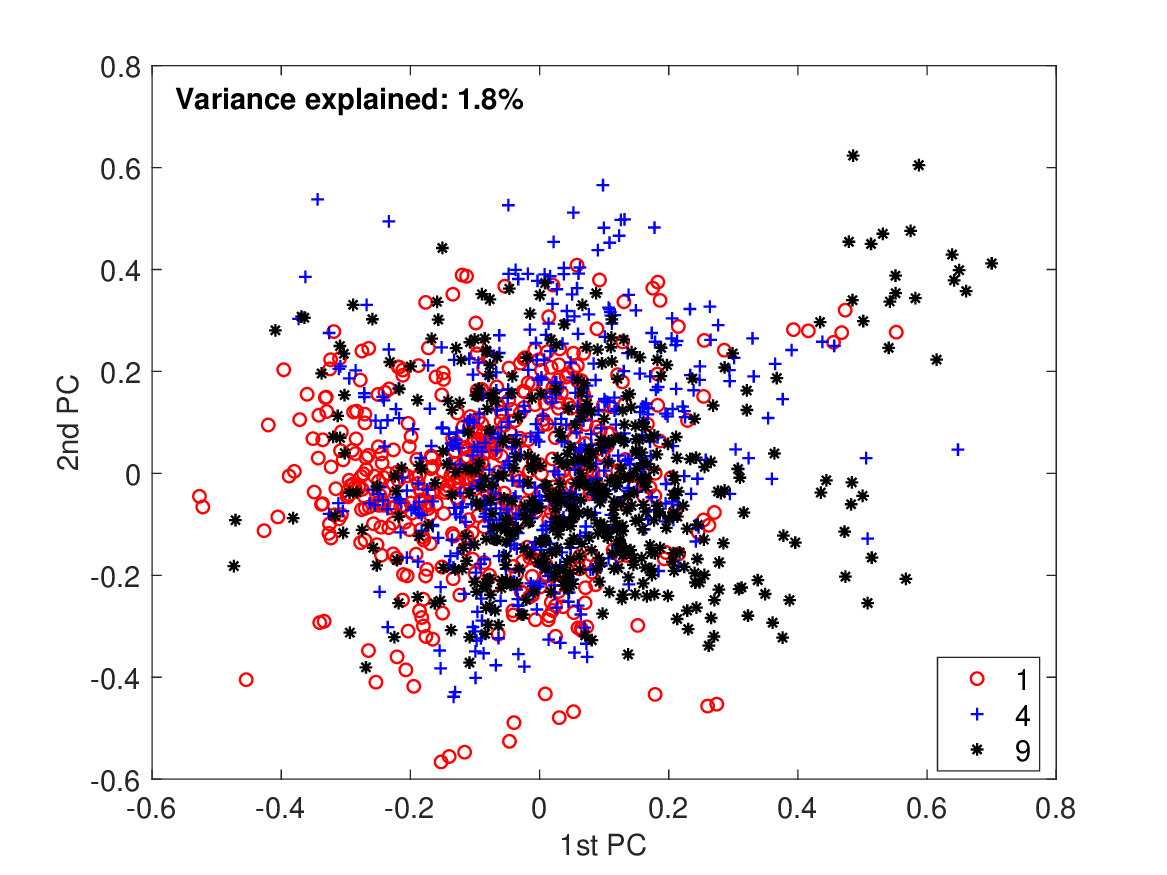}
		\caption{\tiny{\textsf{DP-Gauss*}}}
	\end{subfigure}
	\caption{Comparison of our method and \textsf{DP-Gauss*} \cite{dwork2014analyze} in differentially private PCA on MNIST dataset. The privacy constraints are $\eps=2$  and $\delta=0.1$. The total sample size is $n=1500$. All images are downscaled to a size $14\times 14$. 
	}
	\label{fig:MNIST}
\end{figure}

\section{Discussion}

In this paper, we establish optimal rates of convergence, up to logarithmic factors, for differentially private estimation of both the principal components and the covariance matrix under the spiked covariance model. We propose computationally efficient algorithms, and our results accommodate a diverging rank and a wider range of signal strengths.

\subsection{Private estimation of unknown rank} 
In the present paper, we assume that the number of components, $r$, in the spiked covariance model is known. However, in practice, $r$ is typically unknown.  The consistent estimation of the rank $r$ in such models has been extensively studied in the conventional setting (see, e.g., \cite{lam2012factor}, \cite{cai2015optimal}, \cite{ke2023estimation} and references therein). For instance, we can use the eigen-ratio estimator $\widehat r:=\arg\max_{1\leq k\leq R} \big(\lambda_k(\widehat \Sigma)+Z_k\big)\big(\lambda_{k+1}(\widehat\Sigma)+Z_{k+1}\big)^{-1}$, where $Z_1,\cdots, Z_R$ are i.i.d centered $\mcN\big(0, 8\Delta_2^2\eps^{-2}\log(2.5/\delta)\big)$ noise, and $R\leq (n\wedge p)$ is a postulated upper bound for $r$. Lemma~\ref{lem:sense-lambda} ensures that $\widehat r$ is $(\eps/2, \delta/2)$-DP with probability $1-4n^{-99}$. By the concentration property of edge eigenvalues and the sticking property of bulk eigenvalues of the sample covariance matrix \cite{bloemendal2016principal}, we can show that $\widehat r$ is a consistent estimate of $r$ as long as the signal-to-noise ratio satisfies $\lambda/\sigma^2\geq C_0p\sqrt{\log(R)\log(2.5/\delta)}/(n\eps)$ for some large enough constant $C_0>0$, in addition to the conditions required in Theorem~\ref{thm:dp-Sigma-upb}.

\subsection{Private estimation of unknown eigenvalue} 

We assume the signal strength $\lambda$ is known for simplicity. It can be estimated by averaging the first $r$ eigenvalues of the sample covariance matrix. This approach works well in our numerical experiment on the MNIST dataset. The properties of the sample eigenvalues under the spiked covariance model are well-understood (see, e.g., \cite{bloemendal2016principal} for a precise characterization). To protect privacy, one can also perturb sample eigenvalues with Gaussian noise, as discussed above using sensitivity $\Delta_2$ developed in Lemma~\ref{lem:sense-lambda}. The resultant estimator $\widetilde\lambda$  is consistent in the sense that $C_0^{-1}\leq|\widetilde{\lambda}/\lambda|\leq C_0$ with high probability under the conditions of Theorem~\ref{thm:dp-Sigma-upb} and if $n\eps \geq C_1(r+\log n)$ and $\lambda/\sigma^2\geq C_1p\sqrt{\log(2.5/\delta)}/(n\eps)$,  where $C_0>0$ is some absolute constant and $C_1>0$ is a large constant depending only on $C_0$.

\subsection{Beyond sub-Gaussian distribution} \label{sec:discuss-beyond}
 
Our differentially private PCA method can be extended to distributions beyond the sub-Gaussian case. The primary technical challenge lies in analyzing the sensitivity of the sample spectral projector. The leading term of sensitivity is primarily determined by the quantity 
$
\Delta_1\approx\max_{i\in[n]}\|U^{\top}X_i\|\|U_{\perp}^{\top}X_i\|/ (n\lambda),
$
 where $(\lambda+\sigma^2) I_r \preccurlyeq{\rm cov}\big(U^{\top}X_i\big)\preccurlyeq (\lambda+\sigma^2) I_r$ and $\sigma^2 I_{p-r}\preccurlyeq{\rm cov}(U_{\perp}^{\top}X_i)\preccurlyeq \sigma^2 I_{p-r}$. Assuming that $U^{\top}X_i/\lambda$ (also $U^{\top}_{\perp}X_i/\sigma^2$) has independent entries with a finite fourth moment, we can show that $\|U^{\top}X_i\|\|U_{\perp}^{\top}X_i\|/ (n\lambda)\asymp (\sigma^2/\lambda+\sqrt{\sigma^2/\lambda})\sqrt{p(r+\log n)}/n$ (same as Lemma~\ref{lem:sense-U}) up to additional logarithmic factors with probability $1-O\big(\log^{-1}(n)\big)$. Therefore, we believe that our established bounds for differentially private PCA and covariance matrix estimation are optimal for a wide range of distributions. However, completing the proof still requires establishing sharp upper bounds for the following terms: $\|\widehat\Sigma-\Sigma\|, \|U^{\top}(\widehat\Sigma-\Sigma)U_{\perp}\|$, $\|X_i\|$,  $\|U_{\perp}^{\top}X_i\|$, $\|U^{\top}X_i\|, \big\|h_1\big(\sum_{j\neq i}X_jX_j^{\top}\big)U_{\perp}^{\top}X_i\big\|$, and, $\big\|h_2\big(\sum_{j\neq i}X_jX_j^{\top}\big)U^{\top}X_i\big\|$, etc., which should hold for all $i\in[n]$,  where $h_1(\cdot)$ and $h_2(\cdot)$ are some deterministic functions. 
 It is possible to establish these bounds by imposing some moment condition on the distribution of $X$, such as requiring $\EE |\langle X, a\rangle|^m\leq L^m$ for some sufficiently large $m\geq 4$ and constant $L>0$. See, e.g., \cite{bai2010spectral} and \cite{bloemendal2016principal}. 
Developing these bounds would significantly complicate the current proof and presentation of theoretical results, introduce additional logarithmic factors, and weaken the privacy guarantee. We leave these extensions for future work.

\section{Technical lemmas}\label{sec:tech-lem}

In this section,  we provide some technical lemmas that will be frequently used in the subsequent proofs.  Due to page constraints,  all proofs are given in the Appendix.

Lemma~\ref{lem:con-hatSigma} is a well-known \emph{dimension-free} concentration inequality of sample covariance matrix developed by \cite{koltchinskii2017concentration}. Here, $\|\cdot\|$ denote the spectral norm of a matrix and $\ell_2$-norm of a vector. 

\begin{lemma}[\cite{koltchinskii2017concentration}]\label{lem:con-hatSigma}
	Suppose $X_1, \cdots, X_n$ are i.i.d. sampled from $ \mcN(0, \Sigma)$ and $\whSig := \sum_{i=1}^{n}X_i X_i^{\top}/n$. Then, 
	$$
	\mathbb{E}\|\whSig -\Sigma\| \asymp\Bigg(\sqrt{\frac{{\rm tr}(\Sigma)\|\Sigma\|}{n}} \bigvee \frac{{\rm tr}(\Sigma)}{n}\Bigg). 
	$$
	Moreover,  there exists an absolute constant $C_1>0$ such that,  for all $t \geq 1$,  with probability at least $1-e^{-t}$,
	$$
	\Big|\|\whSig-\Sigma\|-\EE\|\whSig-\Sigma\|\Big| \leq C_1\|\Sigma\|\left(\frac{t}{n}+\sqrt{\frac{t}{n}}\bigg(1+\sqrt{\frac{{\rm tr}(\Sigma)/\|\Sigma\|}{n}}\bigg) \right). 
	$$
	
\end{lemma}

The following lemma characterizes the concentration of the norm of a Gaussian random vector.  Recall that $X_i'$ is an independent copy of $X_i$.
\begin{lemma}\label{lem:con-X-norm}
	Let $X\sim \mcN(0, \Sigma)$ and the eigenvalues of $\Sigma$ are $\lambda_1\geq\cdots\geq \lambda_p\geq 0$.  Then,  there exist absolute constants $C_1,  C_2, c_1>0$ such that
	$$
	\PP\bigg(\Big|\|X\|^2-{\rm tr}(\Sigma) \Big|\leq C_1\Big(u\sum_{i=1}^p \lambda_i^2\Big)^{1/2}+C_2\lambda_1 u\bigg)\geq 1-e^{-c_1u}, 
	$$
	for any $u>0$.  Under the spiked covariance model $\Sigma\in\Theta(\lambda,\sigma^2)$ and the condition that $p\geq C_6\log n$ for some absolute constant $C_6>0$,  we have 
	\begin{align*}
		\PP\bigg(\Big\{\max_{i\in[n]} &\|X_i\|^2+\|X_i'\|^2\leq C_3(r\lambda+p\sigma^2)+C_4\sqrt{(r\lambda^2+p\sigma^4)\log n}+C_5(\lambda+\sigma^2)\log n\Big\}\bigg)\\
		\bigcap& \Big\{\max_{i\in[n]}\|U^{\top}X_i\|^2+\max_{i\in[n]}\|U^{\top}X_i'\|^2\leq C_3r(\lambda+\sigma^2)+C_4\sqrt{r(\lambda^2+\sigma^4)\log n}+C_5(\lambda+\sigma^2)\log n \Big\}\\
		\bigcap& \Big\{\max_{i\in[n]}\|U_{\perp}^{\top}X_i\|^2+\max_{i\in[n]}\|U_{\perp}^{\top}X_i'\|^2\leq C_3p\sigma^2\Big\}\geq 1-n^{-100},
	\end{align*}
	where $C_3, C_4,C_5>0$ are some absolute constants.   Let $\mcE_0$ denote the above event.  Moreover,  
	$$
	\EE \|X_i\|^2\leq C_3(r\lambda+p\sigma^2),\quad  \EE\|U^{\top}X_i\|^2\leq C_3r(\lambda+\sigma^2) \quad {\rm and}\quad \EE \|U_{\perp}^{\top}X_i\|^2\leq C_3p\sigma^2.  
	$$
\end{lemma}

Denote $\Delta:=\whSig-\Sigma$ and $\Delta^{(i)}:=\whSig^{(i)}-\Sigma$. We shall frequently use several concentration bounds related to $\Delta$ and $\Delta^{(i)}$ throughout the proof. For reader's convenience, these concentration bounds are collected in the following lemma.

\begin{lemma}\label{lem:con-bounds}
	Suppose that $\Sigma\in\Theta(\lambda, \sigma^2)$, $n\geq C_1(r\log n+\log^2n)$, $2r\leq p$, and $\lambda/\sigma^2\geq C_1p/n$ for some absolute constant $C_1>0$. There exist absolute constants $c_0>0$  and $C_2>0$ such that  the event
	\begin{align}\label{eq:tech-bounds-bd1}
		\mcE_{\Delta}:=\Bigg\{&\|\Delta\|+\max_{i\in[n]}\|\Delta^{(i)}\|\leq C_2\sqrt{\frac{(\lambda+\sigma^2)(r\lambda+p\sigma^2)}{n}}+\frac{\lambda}{10}\Bigg\} \notag\\
		&\bigcap \Bigg\{\big\|U^{\top}\Delta U_{\perp} \big\|+\max_{i\in[n]}\|U^{\top}\Delta^{(i)}U_{\perp}\|\leq C_3\sqrt{\frac{\sigma^2(\lambda+\sigma^2)p}{n}} \Bigg\}
	\end{align}
	holds with probability $\PP(\mcE_{\Delta})\geq 1-e^{-c_0(n\wedge p)}$.  Meanwhile, we have 
$$
	\EE\|\Delta\|\leq C_2\sqrt{\frac{(\lambda+\sigma^2)(r\lambda+p\sigma^2)}{n}}\quad {\rm and}\quad \EE\big\|U^{\top}\Delta U_{\perp} \big\|\leq C_3\sqrt{\frac{\sigma^2(\lambda+\sigma^2)p}{n}}.
$$

	There exist absolute constants $c_1>0$ and $C_3>0$ such that the event 
	\begin{align}
		\mcE_1:=\Bigg\{&\max_{i\in[n]}\big\|U^{\top}(X_iX_i^{\top}/n)U_{\perp}\big\|+\big\|U^{\top}(X_i'X_i'^{\top}/n)U_{\perp}\big\|\leq C_3\frac{\sqrt{\sigma^2(\lambda+\sigma^2)p(r+\log n)}}{n} \Bigg\}\label{eq:tech-bounds-bd2}\\
		&\bigcap \Bigg\{\max_{i\in[n]}\Big\|U^{\top}\Big(\frac{1}{n}\sum_{j\neq i}X_jX_j^{\top}\Big)U_{\perp}U_{\perp}^{\top}X_i \Big\|\leq C_3\sigma\cdot \sqrt{\frac{\sigma^2(\lambda+\sigma^2)p(r+\log n)}{n}} \Bigg\} \notag
	\end{align}
	holds with probability $\PP(\mcE_1)\geq 1-e^{-c_1(n\wedge p)}-2n^{-99}$.
\end{lemma}

The following perturbation bound of principal subspace will be useful.
\begin{lemma} \label{lem:hatU-bound}
	Suppose that $\lambda_r \geq(4+\delta)\|\Delta\|$ for some $\delta>0$, then
	\begin{equation*}
		\big\|\whU \whU^{\top}-U U^{\top}\big\| \leq 2\left\|\Lambda^{-1} U^{\top} \Delta U_{\perp}\right\| + 6(4+\delta) \frac{\|\Delta\|\left\|U^{\top} \Delta U_{\perp}\right\|}{\delta \lambda_r^2}. 
	\end{equation*}
\end{lemma}

\bibliographystyle{chicago}
\bibliography{reference}

\appendix

\section{More Simulation Results}\label{sec:app-simulation}

\subsection{Covariance Matrix Estimation}

\begin{figure}[H]
	\centering
	\begin{subfigure}[b]{.48\linewidth}
		\includegraphics[width=\linewidth]{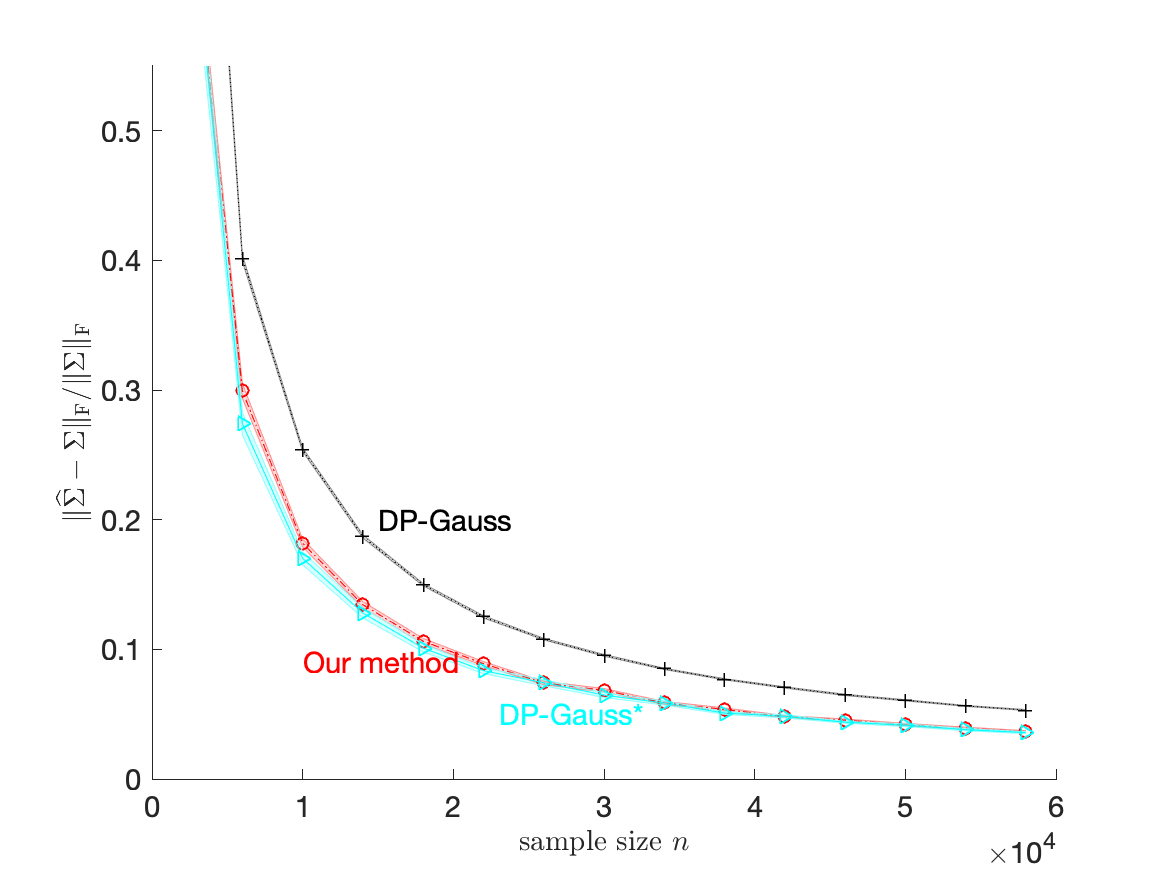}
		\caption{\tiny{Dimension $p=50$. }}
	\end{subfigure}
	\begin{subfigure}[b]{.48\linewidth}
		\includegraphics[width=\linewidth]{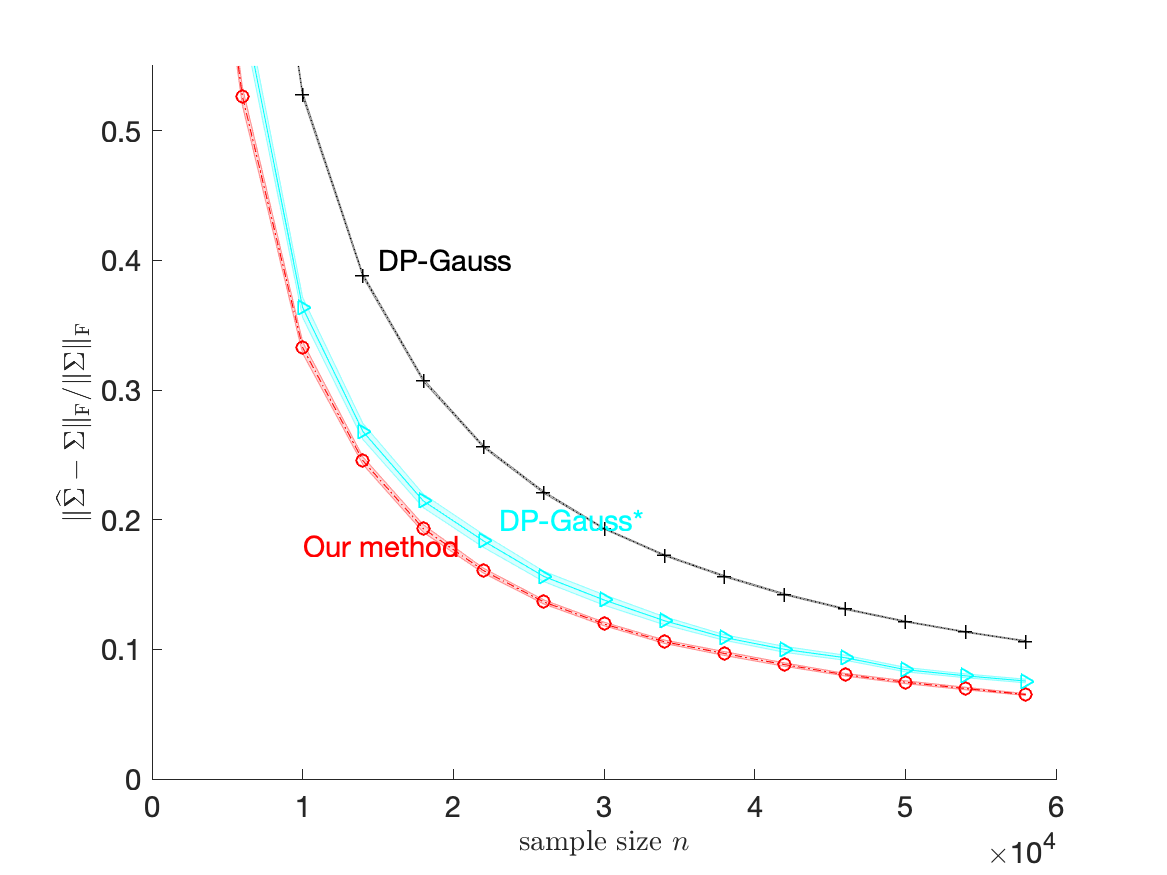}
		\caption{\tiny{Dimension $p=100$.}}
	\end{subfigure}
	\caption{Comparison of our method,  and \textsf{DP-Gauss}, \textsf{DP-Gauss*} \cite{dwork2014analyze} in differentially private covariance matrix estimation when sample size $n$ changes.  The rank  $r=3, \lambda=10$, $\sigma^2=1$, and privacy constraints $\eps=1, \delta=0.1$. 
	}
	\label{fig:S4}
\end{figure}

\begin{figure}[H]
	\centering
	\begin{subfigure}[b]{.48\linewidth}
		\includegraphics[width=\linewidth]{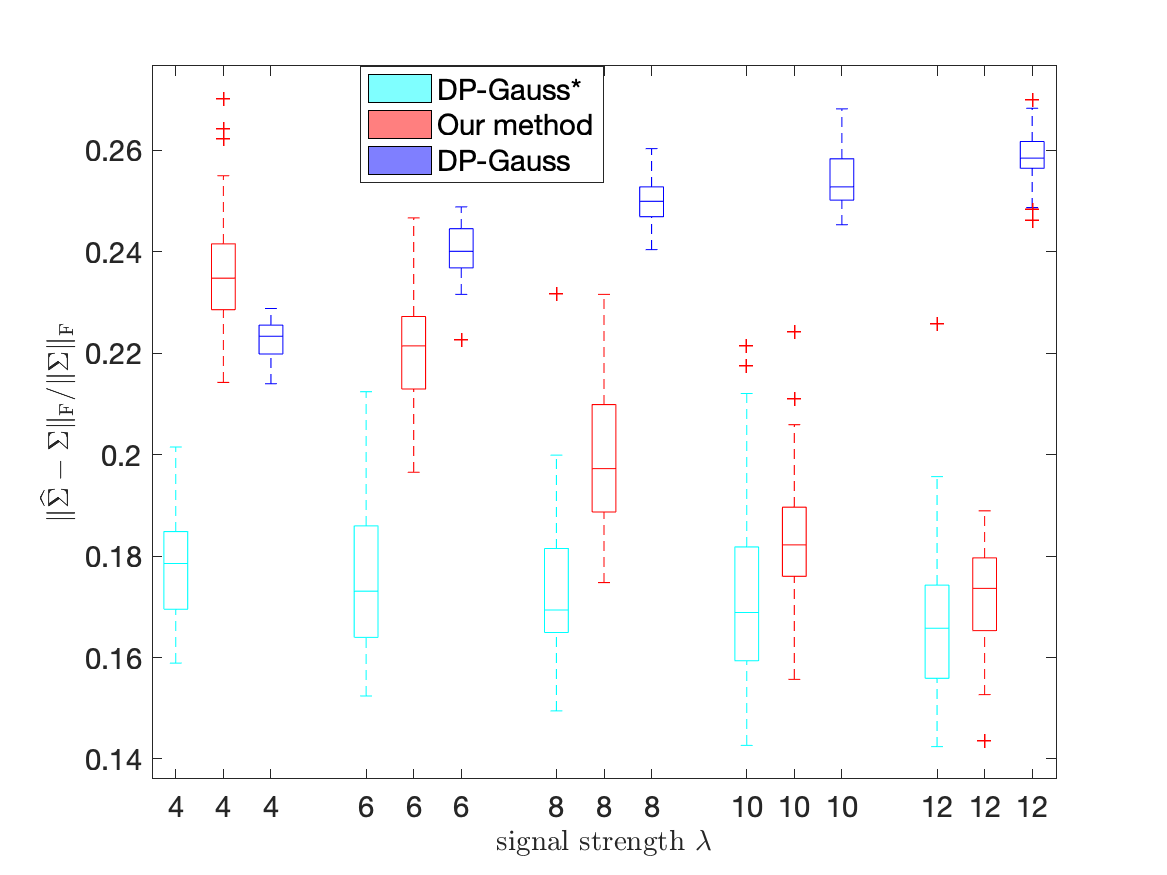}
		\caption{\tiny{Dimension $p=50$. }}
	\end{subfigure}
	\begin{subfigure}[b]{.48\linewidth}
		\includegraphics[width=\linewidth]{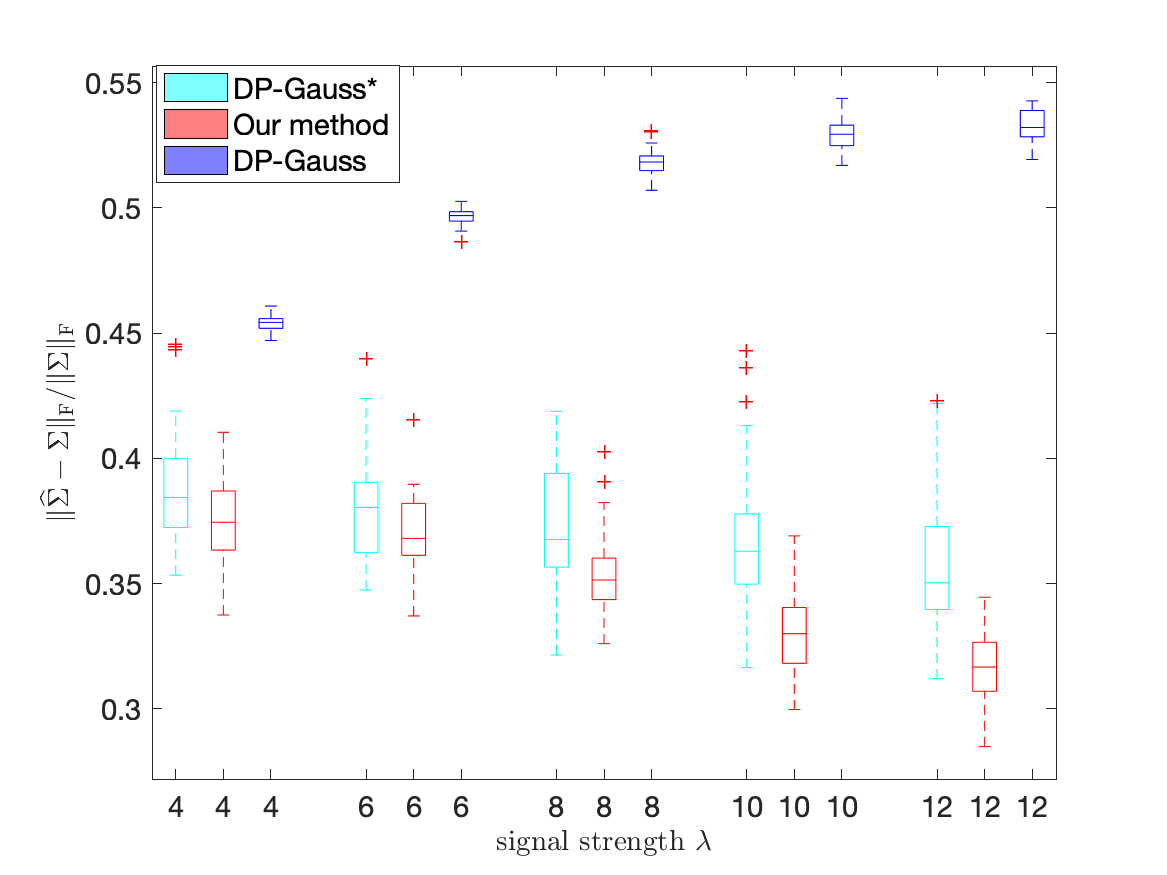}
		\caption{\tiny{Dimension $p=100$.}}
	\end{subfigure}
	\caption{Comparison of our method, and \textsf{DP-Gauss}, \textsf{DP-Gauss*} \cite{dwork2014analyze} in differentially private covariance matrix estimation when signal strength $\lambda$ changes.  The rank  $r=3$, sample size $n=10^4$, $\sigma^2=1$, and privacy constraints $\eps=1, \delta=0.1$. 
	}
	\label{fig:S4-1}
\end{figure}

The \textsf{DP-Gauss} method introduced by \cite{dwork2014analyze} can also be used for covariance matrix estimation. In the fourth simulation setting, we compare the performance of our method and \textsf{DP-Gauss} in spiked covariance matrix estimation. The rank is set to $r=3$, and the privacy constraints are fixed to $\eps=1$ and $\delta=0.1$. We observe that \textsf{DP-Gauss*} performs similarly as our method when $p=50$, but becomes less efficient than our method when $p=100$.  The results are displayed in Figure~\ref{fig:S4}. Similar phenomenon is also observed when signal strength changes.  See Figure~\ref{fig:S4-1}.

\section{Proofs of Theorems}

	\subsection{Proof of Theorem~\ref{thm:dp-pca-upb}}
	
	It suffices to bound $\|\wtU\wtU^{\top}-\widehat U\widehat U^{\top}\|_q$ and $\|\widehat U\widehat U^{\top}-UU^{\top}\|_q$.   Without loss of generality,  we start with $q=\infty$,  i.e.,  the upper bound of spectral norm. 
	
	Recall that $\wtU\wtU^{\top}$ denotes the spectral projector of the top-$r$ eigenvectors of $\Sigma+\Delta$.  Based on Lemma~\ref{lem:con-bounds},  we have
	$$
	\|\Delta\|\leq C_2\sqrt{\frac{(\lambda+\sigma^2)(r\lambda+p\sigma^2)}{n}}+\frac{\lambda}{10}
	$$
	with probability at least $1-e^{-c_0(n\wedge p)}$.   In this event,  we have $\lambda\geq 5\|\Delta\|$ since $n\geq C_1r$ and $\lambda/\sigma^2\geq C_1\big(p/n+\sqrt{p/n}\big)$.  Applying Lemma~\ref{lem:hatU-bound},  we get in this event,  
	$$
	\left\| \whU\whU^{\top}  -U U^{\top}\right\| \leq 2\left\|\Lambda^{-1} U^{\top} \Delta U_{\perp}\right\|+C_3 \frac{\|\Delta\|\left\|U^{\top} \Delta U_{\perp}\right\|}{\lambda^2}\leq C_4\bigg(\frac{\sigma^2}{\lambda}+\sqrt{\frac{\sigma^2}{\lambda}}\bigg)\sqrt{\frac{p}{n}},
	$$
	where the last inequality holds with probability at least $1-e^{-c_1(n\wedge p)}$ by Lemma~\ref{lem:con-bounds}. Denote this event as $\mcF_1$ so that $\PP(\mcF_1)\geq 1-e^{-c_1(n\wedge p)}$.
	
	Then, 
	\begin{align*}
	\EE \|\whU\whU^{\top}-UU^{\top}\|\leq& C_4\bigg(\frac{\sigma^2}{\lambda}+\sqrt{\frac{\sigma^2}{\lambda}}\bigg)\sqrt{\frac{p}{n}}+2\EE \mathbbm{1}(\mcF_1^{\rm c})\\
	\leq&C_4\bigg(\frac{\sigma^2}{\lambda}+\sqrt{\frac{\sigma^2}{\lambda}}\bigg)\sqrt{\frac{p}{n}},
	\end{align*}
	where we used the condition $\lambda/\sigma^2\leq (p/n)e^{c_2(n\wedge p)}$ for some small absolute constant $c_2>0$.

	It follows from the  Davis-Kahan Theorem that
	$
	\|\wtU\wtU^{\top}-\wtU\wtU^{\top}\|\lesssim \|Z\| \wedge 1,
	$
	where $Z$ is a symmetric matrix with i.i.d.  entries having distribution $\mcN\big(0,  8\Delta_1^2\varepsilon^{-2}\log\frac{2.5}{\delta}\big)$.  By  \cite[Theorem 4.4.5]{vershynin2020high}, 
	$$
	\|Z \|  \lesssim \bigg(\frac{\sigma^2}{\lambda}+\sqrt{\frac{\sigma^2}{\lambda}} \bigg)\cdot \frac{ p\sqrt{r+\log n}}{n\varepsilon} \log^{1/2}  (\frac{2.5}{\delta}), 
	$$  
	with probability at least $1-O\left( e^{-c_1(n\wedge p)}\right)$ and the same bound also holds for $\EE\|Z\|$.

	Combining the above two bounds,  we can conclude that with probability at least $1-e^{-c_1(n\wedge p)}$,
	\begin{align*}
		& \|\wtU\wtU^{\top}  - U U^{\top}\| \lesssim \bigg(\frac{\sigma^2}{\lambda}+\sqrt{\frac{\sigma^2}{\lambda}}\bigg)\cdot \bigg(\sqrt{\frac{p}{n}}+\frac{p\sqrt{(r+\log n)}}{n\varepsilon}\sqrt{\log\frac{2.5}{\delta}}\bigg),
	\end{align*}
	and the same bound also holds for $\EE \|\wtU\wtU^{\top}  - U U^{\top}\| $ as long as $e^{c_2(n\wedge p)}\geq n/p$ .
	In the same event,  we can also get the upper bound of nuclear norm distance:
	\begin{align*}
		& \|\wtU\wtU^{\top}  - U U^{\top}\|_{\ast} \lesssim \bigg(\frac{\sigma^2}{\lambda}+\sqrt{\frac{\sigma^2}{\lambda}}\bigg)\cdot r\bigg(\sqrt{\frac{p}{n}}+\frac{p\sqrt{(r+\log n)}}{n\varepsilon}\sqrt{\log\frac{2.5}{\delta}}\bigg)
	\end{align*}
	and $\EE \|\wtU\wtU^{\top}  - U U^{\top}\|_{\ast}$ as long as $\lambda/\sigma^2\leq (p/n)e^{c_2(n\wedge p)}$ . 
	For general $q\in[1,\infty]$,  we can apply the interpolation inequality:
	$$
	\|\wtU\wtU^{\top}-UU^{\top}\|_q\leq \|\wtU \wtU^{\top}-UU^{\top}\|_{\ast}^{1/q}\cdot \|\wtU\wtU^{\top}-UU^{\top}\|^{\frac{q-1}{q}},
	$$
	which completes the proof.

	\subsection{Proof of Theorem~\ref{thm:dp-Sigma-upb}}
	
	Recall that we assume $\sigma^2$ is known.  By the definition of $\wtSig$ stated in Algorithm~\ref{alg:DP-PCA},   for any $q\in[1, \infty]$, we have the Schatten-$q$ norm bounded as
	\begin{align*}
		\|\wtSig-\Sigma\|_q=&\|\wtU \wtLambda \wtU^{\top}-U\Lambda U^{\top}\|_q\\
		=&\big\|\wtU\big(\wtU^{\top}(\whSig-\sigma^2I_p+E) \wtU\big) \wtU^{\top}-UU^{\top}(\Sigma-\sigma^2 I_p) UU^{\top} \big\|_q\\
		\leq& \big\|\wtU \wtU^{\top}(\whSig-\sigma^2 I_p)\wtU\wtU^{\top}-UU^{\top}(\Sigma-\sigma^2I_p)UU^{\top}\big\|_q+\|E\|_q.
	\end{align*}
	Without loss of generality,  we begin with $q=\infty$ and bound the spectral norm.  Since $E$ is an $r\times r$ symmetric matrix with i.i.d.  entries $\mcN\big(0,  8\Delta_2^2\varepsilon^{-2}\log \left(\frac{2.5}{\delta}\right) \big)$,  by \cite[Theorem~4.4.5]{vershynin2020high},  we get 
	\begin{align}\label{eq:proof-dp-Sigma-upb-bd0}
		\|E\|\leq C_4 \frac{\lambda(r+\log n)^{3/2}+\sigma^2(p+\log n)\sqrt{r+\log n}}{n\varepsilon}\cdot \sqrt{\log\frac{2.5}{\delta}},
	\end{align}
	with probability at least $1-n^{-100}$ and for some absolute constant $C_4>0$. Moreover, the same bound holds for $\EE\|E\|$. 
	
	Observe that 
	\begin{align}
		\big\|\wtU \wtU^{\top}&(\whSig-\sigma^2 I_p)\wtU\wtU^{\top}-UU^{\top}(\Sigma-\sigma^2I_p)UU^{\top}\big\|\notag\\
		\leq& \big\|\big(\wtU \wtU^{\top}-UU^{\top}\big)\big\| \big\|(\whSig-\sigma^2 I_p)\big\|+\|UU^{\top}(\whSig -\Sigma)\|+ \|\Sigma-\sigma^2 I_p\|\|\wtU\wtU^{\top}-UU^{\top}\|.\label{eq:proof-dp-Sigma-upb-bd1}
	\end{align}
	Since
	$$
	\|\whSig-\sigma^2 I_p\|\leq \|\Sigma-\sigma^2 I_p\| + \|\whSig-\Sigma\|\lesssim \lambda,
	$$
	where the last inequality holds in event $\mcE_{\Delta}$ defined in Lemma~\ref{lem:con-bounds} and under the conditions $n\geq C_1r$ and $\lambda/\sigma^2\geq C_1\big(p/n+\sqrt{p/n}\big)$.   By (\ref{eq:proof-dp-Sigma-upb-bd1}),  we get 
	\begin{align}
		\big\|\wtU \wtU^{\top}&(\whSig-\sigma^2 I_p)\wtU\wtU^{\top}-UU^{\top}(\Sigma-\sigma^2I_p)UU^{\top}\big\|\notag\\
		\lesssim & \lambda \|\wtU\wtU^{\top}-UU^{\top}\|+\|U^{\top}(\whSig-\Sigma)\|\notag\\
		\lesssim& \sqrt{\sigma^2(\lambda+\sigma^2)} \bigg(\sqrt{\frac{p}{n}}+\frac{p\sqrt{(r+\log n)}}{n\varepsilon}\sqrt{\log\frac{2.5}{\delta}}\bigg)+\sqrt{\frac{(\lambda+\sigma^2)(r\lambda+p\sigma^2)}{n}}\notag\\
		\lesssim& \sqrt{\sigma^2(\lambda+\sigma^2)} \bigg(\sqrt{\frac{p}{n}}+\frac{p\sqrt{(r+\log n)}}{n\varepsilon}\sqrt{\log\frac{2.5}{\delta}}\bigg)+\lambda\sqrt{\frac{r}{n}},\label{eq:proof-dp-Sigma-upb-bd2}
	\end{align}
	where the last inequality is due to Theorem~\ref{thm:dp-pca-upb} and Lemma~\ref{lem:con-bounds}.    Combining (\ref{eq:proof-dp-Sigma-upb-bd0}) and (\ref{eq:proof-dp-Sigma-upb-bd2}),  we get 
	\begin{align*}
		\|\wtSig-\Sigma\|\lesssim \lambda\bigg(\sqrt{\frac{r}{n}}+\frac{(r+\log n)^{3/2}}{n\varepsilon}\cdot \sqrt{\log\frac{2.5}{\delta}}\bigg)+\sqrt{\sigma^2(\lambda+\sigma^2)} \bigg(\sqrt{\frac{p}{n}}+\frac{p\sqrt{(r+\log n)}}{n\varepsilon}\sqrt{\log\frac{2.5}{\delta}}\bigg),
	\end{align*}
	with probability at least $1-3n^{-99}-e^{-c_1(n\wedge p)}$.

	Moreover,
	\begin{align*}
		\EE \big\|\wtU \wtU^{\top}&(\whSig-\sigma^2 I_p)\wtU\wtU^{\top}-UU^{\top}(\Sigma-\sigma^2I_p)UU^{\top}\big\|\notag\\
		\lesssim & \EE \big\|\big(\wtU \wtU^{\top}-UU^{\top}\big)\big\| \big\|(\whSig-\sigma^2 I_p)\big\|+\EE\|UU^{\top}(\whSig -\Sigma)\|+ \EE\|\Sigma-\sigma^2 I_p\|\|\wtU\wtU^{\top}-UU^{\top}\|\\
		\lesssim &\EE^{1/2} \big\|\big(\wtU \wtU^{\top}-UU^{\top}\big)\big\|^2\EE^{1/2}\big\|(\whSig-\sigma^2 I_p)\big\|^2+\sqrt{\frac{(\lambda+\sigma^2)(r\lambda+p\sigma^2)}{n}}\\
		&\quad+\sqrt{\sigma^2(\lambda+\sigma^2)} \bigg(\sqrt{\frac{p}{n}}+\frac{p\sqrt{(r+\log n)}}{n\varepsilon}\sqrt{\log\frac{2.5}{\delta}}\bigg).
	\end{align*}
	Note that 
	\begin{align*}
		\EE\|\wtU\wtU^{\top}-UU^{\top}\|^2\leq& 2\EE\|\whU\whU-\wtU\wtU^{\top}\|^2+2\EE\|\wtU\wtU^{\top}-UU^{\top}\|^2\\
		\lesssim&\bigg(\frac{\sigma^2}{\lambda}+\sqrt{\frac{\sigma^2}{\lambda}}\bigg)\cdot \bigg(\sqrt{\frac{p}{n}}+\frac{p\sqrt{r+\log n}}{n\varepsilon}\sqrt{\log\frac{2.5}{\delta}}\bigg),
	\end{align*}
	where the last inequality is due to the classical concentration of operator norm of a Gaussian random matrix,  e.g., \cite{koltchinskii2016perturbation}, the condition that $\lambda/\sigma^2\leq (p/n)e^{c_2(n\wedge p)}$ as in the proof of Theorem~\ref{thm:dp-pca-upb}, 
	and 
	\begin{align*}
		\|\whSig-\sigma^2 I_p\|^2\leq& 2\|\whSig-\Sigma\|^2+2\|\Sigma-\sigma^2 I_p\|^2\\
		\lesssim& \lambda+\sqrt{\frac{(\lambda+\sigma^2)(r\lambda+p\sigma^2)}{n}}\lesssim \lambda,
	\end{align*}
	where the first inequality can be obtained by integrating the probability bound in Lemmas~\ref{lem:con-hatSigma}  and \ref{lem:con-bounds}. Therefore, we conclude that 
	$$
	\EE \|\wtSig-\Sigma\|\lesssim \lambda\bigg(\sqrt{\frac{r}{n}}+\frac{(r+\log n)^{3/2}}{n\varepsilon}\cdot \sqrt{\log\frac{2.5}{\delta}}\bigg)+\sqrt{\sigma^2(\lambda+\sigma^2)} \bigg(\sqrt{\frac{p}{n}}+\frac{p\sqrt{(r+\log n)}}{n\varepsilon}\sqrt{\log\frac{2.5}{\delta}}\bigg). 
	$$

	Similar bounds can also be derived for the nuclear norm distance $\|\wtSig-\Sigma\|_{\ast}$ and the general Schatten-$q$ norm $\|\wtSig-\Sigma\|_q$.  The detailed proof is skipped.

\subsection{Proof of Theorem~\ref{thm:dp-pca-lwb}}

	
	Some preliminary results on the KL-divergence and total variation distance between Gaussian distributions are required.	Let $\mcN(\mu_1, \Sigma_1)$ and $\mcN(\mu_2, \Sigma_2)$  be two $p$-dimensional multivariate Gaussians, then 
	\begin{align*}
		& \operatorname{KL}\left(\mcN\left(\mu_1, \Sigma_1\right) \| \mcN\left(\mu_2, \Sigma_2\right)\right) \\ & = \frac{1}{2}\left(\operatorname{Tr}\left(\Sigma_2^{-1} \Sigma_1 - I_p\right) +\left(\mu_2-\mu_1\right)^{\top} \Sigma_2^{-1}\left(\mu_2-\mu_1\right)+\log \left(\frac{\operatorname{det} \Sigma_2}{\operatorname{det} \Sigma_1}\right)\right). 
	\end{align*} 
	Suppose $U_i, U_j\in \OO_{p,r}$ satisfying $\big\|U_i U_i^\top - U_j U_j^{\top}\big\|_{\rm F} \leq \varepsilon_0$. Let $\lambda, \sigma^2 \geq 0$ be constants and define $\Sigma_i = \lambda U_i U_i^{\top} + \sigma^2 I_p$ and $\Sigma_j = \lambda U_j U_j^{\top}+\sigma^2 I_p$, respectively.  Then it is easy to check that 
	$$
	\begin{aligned}
		\operatorname{KL} \big(\mcN(0, \Sigma_1) \|& \mcN(0, \Sigma_j) \big) 
		= \frac{1}{2}\left(\operatorname{Tr}\left(\Sigma_j^{-1} \Sigma_i - I_p\right)+\log \left(\frac{\operatorname{det} \Sigma_j}{\operatorname{det} \Sigma_i}\right)\right) \\ 
		& =\frac{1}{2} \left (\frac{\lambda}{\sigma^2} - \frac{\lambda}{\lambda+\sigma^2} \right)\left[ r - \operatorname{Tr} \left(U_j U_j^{\top}U_i U_i^{\top}\right)  \right] 
		\leqslant \frac{1}{2} \frac{\lambda^2}{\sigma^2(\sigma^2+\lambda)}  \varepsilon_0^2,  
	\end{aligned}
	$$ 
	and further by Pinsker's inequality, we have 
	\begin{align*}
		& \operatorname{\operatorname{TV}}\left( \mcN(0, \Sigma_i),  \mcN(0, \Sigma_j) \right) \leq \sqrt{\frac{1}{2} \operatorname{KL} \left(\mcN(0, \Sigma_i) \| \mcN(0, \Sigma_j) \right)} \leq \frac{1}{2} \varepsilon_0 \sqrt{\frac{\lambda^2}{\sigma^2(\sigma^2+\lambda)}}. 
	\end{align*}

	In order to apply Fano's lemma,  we need to construct a large subset of $\OO_{p, r}$ within which the elements are well-separated.  Towards that end,  we apply existing results of the packing number of Grassmannians.    Indeed,  by \cite[Proposition 8]{pajor1998metric} and \cite[Lemma 5]{koltchinskii_xia15},  for any $q\in[1, \infty]$,  there exists an absolute constant $c'>0$ and a subset $\mcS_q^{(p-r)}\subset \OO_{p-r, r}$ such that for any $V_i\neq V_j\in \mcS_q^{(p-r)}$,
	$$
	\|V_iV_i^{\top}-V_jV_j^{\top}\|_{q}\geq c'r^{1/q}
	$$
	and the cardinality of $\mcS_q^{(p-r)}$ is at least $2^{r(p-r)}$.  Here, $\|\cdot\|_q$ denotes the Schatten-$q$ norm of a matrix.  In particular,  spectral norm is Schatten-$\infty$ norm,  Frobenius norm is Schatten-$2$ norm,  and nuclear norm is Schatten-$1$ norm.  Let $\varepsilon_0>0$ be a small number to be decided later.  Now,  for each $V\in\mcS_q^{(p-r)}$,  we define
	$$
	U=\left(
	\begin{array}{c}
		\sqrt{1-\varepsilon_0^2}I_{r}\\
		\sqrt{\varepsilon_0^2}V
	\end{array}
	\right)
	$$
	such that $U\in\RR^{p\times r}$ and $U^{\top}U=I_r$.  This means that,  for any $V\in\mcS_q^{(p-r)}$,  we can construct a $U\in\OO_{p, r}$.  This defines a subset $\mcS_q^{(p)}\subset \OO_{p,r}$ with ${\rm Card}\big(\mcS_q^{(p)}\big)\geq 2^{r(p-r)}$ such that for any $U_i\neq U_j\in \mcS_q^{(p)}$,
	\begin{align*}
		\|U_iU_i^{\top}-U_jU_j^{\top}\|_q\geq \sqrt{\varepsilon_0^2(1-\varepsilon_0^2)} \|V_i-V_j\|_q\gtrsim \sqrt{\varepsilon_0^2(1-\varepsilon_0^2)}\|V_iV_i^{\top}-V_jV_j^{\top}\|_q\gtrsim \sqrt{\varepsilon_0^2(1-\varepsilon_0^2)} r^{1/q}
	\end{align*}
	and,  meanwhile, 
	$$
	\|U_iU_i^{\top}-U_jU_j^{\top}\|_{\rm F}\lesssim \|U_i-U_j\|_{\rm F}\leq \varepsilon_0 \|V_i-V_j\|_{\rm F}\leq \sqrt{2r}\varepsilon_0.
	$$
	We then consider a family of distributions as 
	\begin{align*}
		\mcP\big(\mcS_q^{(p)},\lambda,\sigma^2\big) = \{\mcN(0, \Sigma)^{\ot n}: \Sigma = \lambda U U^{\top}+\sigma^2 I_p,\ U\in \mcS_q^{(p)}\} \subset \mcP(\lambda,\sigma^2),
	\end{align*}
	whose cardinality $N:={\rm Card}\Big(\mcP\big(\mcS_q^{(p)},\lambda,\sigma^2\big)\Big)\geq 2^{r(p-r)}$. 
	
	For $i\neq i^{\prime}\in [N]$, the probability measures $P_i = \mcN(0, \Sigma_i)^{\ot n}$ and $P_{i^{\prime}} = \mcN(0, \Sigma_{i^{\prime}})^{\otimes n}$ in $\mcP\big(\mcS_q^{(p)}, \lambda, \sigma^2\big)$ satisfy 
	$$ 
	\sum_{k\in[n]}\operatorname{\operatorname{TV}}\left( \mcN(0, \Sigma_i),  \mcN(0, \Sigma_{i^{\prime}}) \right) \lesssim \frac{n}{2} \sqrt{r}\varepsilon_0 \sqrt{\frac{\lambda^2}{\sigma^2(\sigma^2+\lambda)}}, 
	$$
	and 
	\begin{equation*}
		\max_{i\ne j \in [N]}  \operatorname{KL} \left(\mcN(0, \Sigma_i)^{\ot n} \| \mcN(0, \Sigma_{i^{\prime}})^{\ot n} \right) \lesssim \frac{n}{2} \frac{\lambda^2 }{\sigma^2(\sigma^2+\lambda)} \varepsilon_0^2r.
	\end{equation*}
	
	To invoke Lemma~\ref{lem:dp-fano}, we define the metric $\rho: \OO_{p, r}\times \OO_{p,r}\mapsto \RR^+$ as $\rho(U_i, U_j):=\|U_iU_i^{\top}-U_jU_j^{\top}\|_q$ for any $q\in[1,\infty]$ and take  $\rho_0\asymp \tau\varepsilon_0r^{1/q}$,  
	$$
	l_0 =c_0\frac{n}{2}\frac{\lambda^2 }{\sigma^2(\sigma^2+\lambda)} \varepsilon_0^2r\quad {\rm and}\quad t_0 =  c_0\frac{n}{2}\sqrt{r} \varepsilon_0 \sqrt{\frac{\lambda^2}{\sigma^2(\sigma^2+\lambda)}}
	$$ 
	for some small absolute constant $c_0,\tau>0$.  Then, by Lemma~\ref{lem:dp-fano}, for any $(\varepsilon,\delta)$-DP estimator $\wtU$, 
	\begin{align*}
		&\sup_{P \in \mcP\big(\mcS_q^{(p)},\lambda,\sigma^2\big) } \EE \big\|\tdU\tdU^{\top} - UU^{\top} \big\|_q  \\ 
		& \geqslant \max \left\{\frac{\tau\varepsilon_0 r^{1/q}}{2} \left(1-\frac{  c_0\frac{n}{2}\frac{\lambda^2 }{\sigma^2(\sigma^2+\lambda)} \varepsilon_0^2r +\log 2}{\log N}\right), \frac{\tau \varepsilon_0 r^{1/q}}{4}\left(1 \wedge \frac{N-1}{\exp \left( 4 \varepsilon  c_0\frac{n}{2} \sqrt{r}\varepsilon_0 \sqrt{\frac{\lambda^2}{\sigma^2(\sigma^2+\lambda)}} \right)}\right)\right\},
	\end{align*}
where we have used the condition that 
$$
\delta\leq c_0'\exp\Big\{2\eps-c_0\big(\eps\sqrt{npr}+pr\big)\Big\},
$$
for some small constants $c_0, c_0'>0$.

Recall that $N\geq 2^{rp/2}$ if $p\geq 2r$. We can take 
	$$
	\varepsilon_0 \asymp \sqrt{\frac{\sigma^2(\lambda+\sigma^2)}{\lambda^2}} \sqrt{\frac{p}{n}} + \sqrt{\frac{\sigma^2(\lambda+\sigma^2)}{\lambda^2}}\frac{p\sqrt{r}}{n\varepsilon},
	$$ 
	and get
	\begin{align*}
		\sup_{P \in \mcP\big(\mcS_q^{(p)}, \lambda,\sigma^2\big) } \EE \Big\|\tdU\tdU^{\top} - UU^{\top} \Big\|_q \gtrsim \sqrt{\frac{\sigma^2(\lambda+\sigma^2)}{\lambda^2}}\cdot r^{1/q} \sqrt{\frac{p}{n}} + \sqrt{\frac{\sigma^2(\lambda+\sigma^2)}{\lambda^2}}\cdot r^{\frac{1}{2}+\frac{1}{q}}\frac{p}{n\varepsilon}. 
	\end{align*} 
	Since a trivial upper of $\|\wtU\wtU^{\top}-UU^{\top}\|_q\leq (2r)^{1/q}$ and $\mcP\big(\mcS_q^{(p)}, \lambda, \sigma^2\big)\subset \mcP(\lambda,\sigma^2)$, we conclude that, for any $q\in[1,\infty]$,  
	\begin{align*}
		& \inf_{\tdU } \sup_{P \in \mcP(\lambda,\sigma^2) } \EE\Big\|\tdU\tdU^{\top} - UU^{\top} \Big\|_q \gtrsim \bigg(\frac{\sigma^2}{\lambda}+\sqrt{\frac{\sigma^2}{\lambda}}\bigg)\left (r^{1/q}\sqrt{\frac{p}{n}}+r^{\frac{1}{2}+\frac{1}{q}}\frac{p}{n\varepsilon}\right)\bigwedge r^{1/q},
	\end{align*} 
	where the infimum is taken over all possible $(\varepsilon,\delta)$-DP algorithms. Now it suffices to choose $q=1, 2, \infty$ to obtain the bounds in nuclear norm, Frobenius norm, and spectral norm, respectively.

	\subsection{Proof of Theorem~\ref{thm:dp-Sigma-lwb}}
	
	Note that the two terms in the minimax lower bounds are contributed by estimating the eigenvalues and eigenvectors, separately. We begin with the term related to estimating eigenvectors. Consider a subset $\mcP_1(\lambda,\sigma^2)\subset \mcP(\lambda, \sigma^2)$ defined by
	$$
	\mcP_1(\lambda,\sigma^2):=\bigg\{\mcN(0, \Sigma): \Sigma=\lambda UU^{\top}+\sigma^2 I_p\quad {\rm and}\quad U\in\OO_{p, r}\bigg\},
	$$
	where we assume $\lambda$ and $\sigma^2$ are both known. If $\lambda$ is already known,  it suffices to estimate $UU^{\top}$ differentially privately by an estimator $\wtU \wtU^{\top}$ so that we can construct the covariance matrix estimator $\wtSig=\lambda \wtU \wtU^{\top}+\sigma^2 I_p$.  Therefore, 
	\begin{equation}\label{eq:proof-lwb-Sigma-bd1}
		\begin{aligned}
			\inf_{\wtSig} \underset{P\in\mcP_1(\lambda,\sigma^2)}{\sup} \EE\big\|\tdSig -& \Sigma \big\| 
			\geq \inf_{\wtU} \underset{P\in\mcP_1(\lambda,\sigma^2)}{\sup}\lambda\cdot \EE \big\|\wtU\wtU^{\top} - UU^{\top} \big\| \\ 
			& \gtrsim  \sqrt{\sigma^2(\lambda+\sigma^2)} \Big( \sqrt{\frac{p}{n}}+  \frac{\sqrt{r}p}{n \varepsilon}\Big) \bigwedge \lambda, 
		\end{aligned} 
	\end{equation}  
	and 
	\begin{equation}\label{eq:proof-lwb-Sigma-bd2}
		\begin{aligned}
			\inf_{\wtSig} \underset{P\in\mcP_1(\lambda,\sigma^2)}{\sup} \EE\big\|\tdSig -& \Sigma \big\|_{\rm F}
			\geq \inf_{\wtU} \underset{P\in\mcP_1(\lambda,\sigma^2)}{\sup}\lambda\cdot \EE \big\|\wtU\wtU^{\top} - UU^{\top} \big\|_{\rm F} \\ 
			& \gtrsim  \sqrt{\sigma^2(\lambda+\sigma^2)} \Big( \sqrt{\frac{pr}{n}}+  \frac{rp}{n \varepsilon}\Big) \bigwedge \lambda\sqrt{r}, 
		\end{aligned} 
	\end{equation} 
	where the last inequalities in (\ref{eq:proof-lwb-Sigma-bd1}) and (\ref{eq:proof-lwb-Sigma-bd2}) are both due to Theorem~\ref{thm:dp-pca-lwb}. These establish the second term in the minimax lower bounds of Theorem~\ref{thm:dp-Sigma-lwb}. More generally, we can establish the minimax lower bound in Schatten-$q$ norms:
	\begin{align}
		\inf_{\wtSig} \underset{P\in\mcP_1(\lambda,\sigma^2)}{\sup} \EE\big\|\tdSig -& \Sigma \big\|_{\rm q}
		\geq \inf_{\wtU} \underset{P\in\mcP_1(\lambda,\sigma^2)}{\sup}\lambda\cdot \EE \big\|\wtU\wtU^{\top} - UU^{\top} \big\|_{q} \notag\\ 
		& \gtrsim  \sqrt{\sigma^2(\lambda+\sigma^2)} \Big( r^{1/q}\sqrt{\frac{p}{n}}+  \frac{pr^{\frac{1}{2}+\frac{1}{q}}}{n \varepsilon}\Big) \bigwedge \lambda r^{1/q},  \label{eq:proof-dp-Sigma-lwb-bd3}
	\end{align}
	for any $q\in[1,\infty]$. 
	
	We now establish the first term in the minimax lower bounds $\lambda\big(\sqrt{r/n}+r/(n\varepsilon)\big)$, which is contributed by estimating the singular values. It is unrelated to the nuisance variance $\sigma^2$ and the eigenvectors $U$. Without loss of generality, we can assume $\sigma^2=0$ and $U_{V}$ in the format 
	$$
	U_{V}=\left(
	\begin{array}{c}
		V\\
		{\bf 0}_{(p-r)\times r}
	\end{array} 
	\right),
	$$
	where $V=[V_0, V_{0\perp}]\in\OO_{r,r}$ with some $V_0\in\OO_{r,r/4}$. It is easy to check that $U_V^{\top}U_V=V^{\top}V=I_r$.  Define 
	$$
	\Lambda_0={\rm diag}(\underbrace{2\lambda, \cdots, 2\lambda}_{r/4}, \lambda,\cdots,\lambda),
	$$ 
	which is an $r\times r$ diagonal matrix. For any $V=[V_0, V_{0\perp}]\in\OO_{r,r}$ with $V_0\in\OO_{r, r/4}$, we consider the following covariance matrix 
	$$
	\Sigma_{V_0}:=U_V\Lambda_0 U_V^{\top}=\left(
	\begin{array}{cc}
		\lambda I_r+\lambda V_0V_0^{\top}&{\bf 0}_{r\times (p-r)}\\
		{\bf 0}_{(p-r)\times r}&{\bf 0}_{(p-r)\times (p-r)}, 
	\end{array}
	\right)
	$$
	To this end, we define a subset $\mcP_2(\lambda)\subset \mcP(\lambda,\sigma^2)$ as 
	$$
	\mcP_2(\lambda):=\left\{\mcN(0,  \Sigma_{V_0}): \Sigma_{V_0}= \left(
	\begin{array}{cc}
		\lambda I_r+\lambda V_0V_0^{\top}&{\bf 0}_{r\times (p-r)}\\
		{\bf 0}_{(p-r)\times r}&{\bf 0}_{(p-r)\times (p-r)}
	\end{array}
	\right),\ V_0\in\OO_{r, r/4}\right\}.
	$$
	We are interested in the minimax lower bound:
	$$
	\inf_{\wtSig} \sup_{P\in \mcP_2(\lambda)} \EE \|\wtSig- \Sigma\|\quad {\rm and}\quad \inf_{\wtSig} \sup_{P\in \mcP_2(\lambda)} \EE \|\wtSig- \Sigma\|_{\rm F},
	$$
	where the infimum is taken over all the possible $(\varepsilon,\delta)$-DP algorithms. 
	
	Observe that if a random vector $X=(X_1,\cdots,X_p)^{\top}\sim \mcN(0, \Sigma_{V_0})$, it means that only the first $r$ entries of $X$ are random variables since other variables are simply zeros with probability one. This suggests that it suffices to consider the reduced problem of estimating an $r\times r$ spiked covariance matrix. Towards that end, define a family of distributions of $r$-dimensional random vector
	$$
	\mcP_2^{(0)}(\lambda):=\Big\{\mcN(0, \Sigma_{0V_0}):\ \Sigma_{0V_0}=\lambda I_r+\lambda V_0V_0^{\top},\ V_0\in\OO_{r, r/4} \Big\}. 
	$$
	Given i.i.d. observations $X_{1}^{(0)},\cdots, X_n^{(0)}\sim P\in \mcP_2^{(r)}(\lambda)$, we aim to estimate the covariance matrix $\Sigma_{0}$ with an $(\varepsilon,\delta)$-differentially private algorithm.  Clearly, by definition,
	\begin{equation}\label{eq:proof-dp-Sigma-lwb-bd4}
		\inf_{\wtSig} \sup_{P\in \mcP_2(\lambda)} \EE \|\wtSig- \Sigma\|_q=\inf_{\wtSig_0} \sup_{P\in \mcP_2^{(0)}(\lambda)} \EE \|\wtSig_0- \Sigma_0\|_q, 
	\end{equation}
	for any Schatten-$q$ norms. It therefore suffices to study the RHS of (\ref{eq:proof-dp-Sigma-lwb-bd4}), which is the differentially private minimax lower bound for estimating an $r\times r$ spiked covariance matrix with rank $r/4$. Without loss of generality, we can still assume $\lambda$ is known and we can immediately invoke the bounds (\ref{eq:proof-lwb-Sigma-bd1}), (\ref{eq:proof-lwb-Sigma-bd2}), (\ref{eq:proof-dp-Sigma-lwb-bd3}) by replacing $\sigma^2=\lambda$ and $r\leftarrow r/4$, $p\leftarrow r$ there, and conclude that 
	
	\begin{equation}\label{eq:proof-lwb-Sigma-bd5}
		\begin{aligned}
			\inf_{\wtSig_0} \underset{P\in\mcP_2^{(0)}(\lambda)}{\sup} \EE\big\|\tdSig -& \Sigma \big\| 
			\gtrsim  \lambda \Big( \sqrt{\frac{r}{n}}+  \frac{r^{3/2}}{n \varepsilon}\Big) \bigwedge \lambda, 
		\end{aligned} 
	\end{equation}  
	and 
	\begin{equation}\label{eq:proof-lwb-Sigma-bd6}
		\begin{aligned}
			\inf_{\wtSig_0} \underset{P\in\mcP_2^{(0)}(\lambda)}{\sup} \EE\big\|\tdSig -& \Sigma \big\|_{\rm F}
			\gtrsim  \lambda \Big( \frac{r}{\sqrt{n}}+  \frac{r^{2}}{n \varepsilon}\Big) \bigwedge \lambda\sqrt{r}, 
		\end{aligned} 
	\end{equation}
	or more generally, 
	\begin{equation}\label{eq:proof-lwb-Sigma-bd7}
		\begin{aligned}
			\inf_{\wtSig_0} \underset{P\in\mcP_2^{(0)}(\lambda)}{\sup} \EE\big\|\tdSig -& \Sigma \big\|_{q}
			\gtrsim  \lambda \Big( \frac{r^{\frac{1}{2}+\frac{1}{q}}}{\sqrt{n}}+  \frac{r^{\frac{3}{2}+\frac{1}{q}}}{n \varepsilon}\Big) \bigwedge \lambda r^{1/q}, 
		\end{aligned} 
	\end{equation}
	for any $q\in[1,+\infty]$.

	Finally, putting together (\ref{eq:proof-lwb-Sigma-bd1})-(\ref{eq:proof-dp-Sigma-lwb-bd3}) and (\ref{eq:proof-lwb-Sigma-bd5})-(\ref{eq:proof-lwb-Sigma-bd7}), we conclude that 
	\begin{align*}
		\inf_{\wtSig}\sup_{P\in\mcP(\lambda,\sigma^2)}\ \EE \|\wtSig-\Sigma\|_q\gtrsim \Bigg(\lambda \Big( \frac{r^{\frac{1}{2}+\frac{1}{q}}}{\sqrt{n}}+  \frac{r^{\frac{3}{2}+\frac{1}{q}}}{n \varepsilon}\Big) +\sqrt{\sigma^2(\lambda+\sigma^2)} \Big( r^{1/q}\sqrt{\frac{p}{n}}+  \frac{pr^{\frac{1}{2}+\frac{1}{q}}}{n \varepsilon}\Big) \Bigg)\bigwedge \lambda r^{1/q}. 
	\end{align*}
	Now by setting $q=1,2,\infty$, we complete the proof.

\subsection{Proof of Theorem~\ref{thm:unknown-sigma}}
Recall that $\widehat \Sigma^{(i)}$ is the version of $\widehat\Sigma$ after replacing $X_i$ by the independent copy $X_i'$, and let $\widehat\sigma^{(i)2}$ be defined similarly using $\widehat\Sigma^{(i)}$ instead of $\widehat \Sigma$. Then,
\begin{align*}
\big|\widehat\sigma^2-\widehat \sigma^{(i)2}\big|=&\left|\frac{\sum_{(p\wedge n)/4\leq k\leq 3(p\wedge n)/4} q_k \big(\lambda_k(\widehat\Sigma)-\lambda_k(\widehat\Sigma^{(i)}) \big)}{\sum_{(p\wedge n)/4\leq k\leq 3(p\wedge n)/4} q_k^2} \right|\\
\leq& \left(\frac{\sum_{(p\wedge n)/4\leq k\leq 3(p\wedge n)/4} \big(\lambda_k(\widehat\Sigma)-\lambda_k(\widehat\Sigma^{(i)})\big)^2}{\sum_{(p\wedge n)/4\leq k\leq 3(p\wedge n)/4} q_k^2}\right)^{1/2}\\
\leq& \frac{C_1}{\sqrt{p\wedge n}}\cdot \|\widehat \Sigma-\widehat\Sigma^{(i)}\|_{\rm F}\\
\leq& \frac{C_1}{\sqrt{p\wedge n}}\cdot \frac{\lambda(r+\log n)+\sigma^2(p+\log n)}{n}, 
\end{align*}
where the second inequality is due to Hoffman-Weilandt's inequality and the fact that $q_k\asymp 1$ for all $(p\wedge n)/4\leq k\leq 3(p\wedge n)/4$, and the last inequality holds with probability at least $1-n^{-100}$ in the event $\mcE_0$. Therefore, by Gaussian mechanism, the estimator $\widetilde\sigma^2$ is $(\eps/3, \delta/3)$-DP in the event $\mcE_0$.

\subsection{Proof of Theorem~\ref{thm:unknown-sigma-2}}
By definition of $\widetilde\sigma^2$,  we have 
$$
|\widetilde\sigma^2-\sigma^2|\leq |\widehat \sigma^2-\sigma^2|+|\widetilde\sigma^2-\widehat\sigma^2|. 
$$
By the concentration of normal random variables, we have 
$$
\PP\bigg(|\widetilde\sigma^2-\widehat\sigma^2|\geq C_4'\sqrt{\frac{\log n}{p\wedge n}}\cdot \frac{\lambda(r+\log n)+\sigma^2(p+\log n)}{n}\bigg)\leq n^{-100}.  
$$
Moreover, by Theorem 1 of \cite{ke2023estimation}, we get 
$$
\PP\bigg(\big|\widehat\sigma^2-\sigma^2\big|\leq \frac{C_4\sigma^2}{n^{3/4}}\bigg)\geq 1-n^{-100}.
$$
Therefore, under the condition (\ref{eq:unknown-sigma-cond}), we have
$$
\PP\Big(\frac{9\sigma^2}{10}\leq\widetilde\sigma^2\leq \frac{11\sigma^2}{10}\Big)\geq 1-2n^{-100},
$$
which ensures that the DP-estimate $\widetilde\sigma^2$ and the true $\sigma^2$ are of the same order. Therefore, the outputs of Algorithm ~\ref{alg:DP-PCA} with sensitivity $\Delta_1$ and $\Delta_2$ using $\widetilde\sigma^2$ can still provide the privacy guarantee.  Moreover,  the same upper bound as in Theorem~\ref{thm:dp-pca-upb} still holds for $\widetilde U\widetilde U^{\top}$ except that the constants in $\log(C/\delta)$ need to be updated.

Note that 
$$
\big\|\widetilde\Sigma-\Sigma \big\|\leq \big\|\widetilde U \widetilde \Lambda\widetilde U^{\top} -U\Lambda U^{\top}\big\|+|\widetilde\sigma^2-\sigma^2|,
$$
where the above analysis shows that 
$$
\PP\bigg(|\widetilde\sigma^2-\sigma^2|\leq C_4'\sqrt{\frac{\log n}{p\wedge n}}\cdot \frac{\lambda(r+\log n)+\sigma^2(p+\log n)}{n\eps}\sqrt{\log\frac{4}{\delta}}+\frac{C_4\sigma^2}{n^{3/4}}\bigg)\geq 1-2n^{-100}. 
$$
Moreover, 
\begin{align*}
\|\wtU \wtLambda& \wtU^{\top}-U\Lambda U^{\top}\|=\big\|\wtU\big(\wtU^{\top}(\whSig-\sigma^2I_p+E) \wtU\big) \wtU^{\top}-UU^{\top}(\Sigma-\sigma^2 I_p) UU^{\top} \big\|\\
		\leq& \big\|\wtU \wtU^{\top}(\whSig-\sigma^2 I_p)\wtU\wtU^{\top}-UU^{\top}(\Sigma-\sigma^2I_p)UU^{\top}\big\|+\|E\|+|\widetilde\sigma^2-\sigma^2|\\
		\lesssim& \sqrt{\sigma^2(\lambda+\sigma^2)} \bigg(\sqrt{\frac{p}{n}}+\frac{p\sqrt{(r+\log n)}}{n\varepsilon}\sqrt{\log\frac{4}{\delta}}\bigg)+\lambda\sqrt{\frac{r}{n}}+|\widetilde\sigma^2-\sigma^2|\\
		\lesssim&C_{5,\gamma} \Bigg(\lambda\bigg(\sqrt{\frac{r}{n}}+\frac{(r+\log n)^{3/2}}{n\varepsilon}\cdot \sqrt{\log\frac{4}{\delta}}\bigg)+\sqrt{\sigma^2(\lambda+\sigma^2)} \bigg(\sqrt{\frac{p}{n}}+\frac{p\sqrt{(r+\log n)}}{n\varepsilon}\sqrt{\log\frac{4}{\delta}}\bigg)\Bigg),
\end{align*}
where the second inequality holds in the event $\mcE_{\Delta}$ as in the proof of Theorem~\ref{thm:dp-Sigma-upb}.

The upper bound for Frobenius-norm error can be established similarly by noticing that 
$$
\big\|\widetilde\Sigma-\Sigma \big\|_{\rm F}\leq \|\widetilde U\widetilde\Lambda\widetilde U^{\top}-U\Lambda U^{\top}\|_{\rm F}+\sqrt{p}|\widetilde\sigma^2-\sigma^2|.
$$

\section{Proofs of Lemmas}\label{sec:proofs} 
\subsection{Proof of Lemma~\ref{lem:DP-alg}}
	
	A non-private estimator of the spectral projector $UU^{\top}$ is the empirical spectral projector $\whU\whU^{\top}$. To privatize the spectral projector $\whU\whU^{\top}$, we introduce randomness according to Gaussian mechanism by Lemma~\ref{lem:gaussian_mechanism}. Let $Z$ be a $p\times p$ symmetric matrix with i.i.d. entries 
	$$
	\left(Z\right)_{ij}\sim \mcN\bigg(0, \frac{8\Delta_1^2}{\varepsilon^2}\log\frac{2.5}{\delta}\bigg), \quad \forall 1\leq i \leq j \leq p. 
	$$ 
	We define a randomized algorithm $\whP$ in the way that  
	$$
	\whP:= \whU\whU + Z. 
	$$ 
	In order to ensure that $\whP$ is $(\varepsilon/2, \delta/2)$-DP, by Lemma~\ref{lem:gaussian_mechanism}, it suffices to choose $\Delta_1$ as follows 
	$$
	\Delta_1 = C_1\bigg(\frac{\sigma^2}{\lambda}+\sqrt{\frac{\sigma^2}{\lambda}}\bigg) \frac{\sqrt{p(r+\log n)}}{n},   
	$$
	for a large enough but absolute constant $C_1>0$. By Lemma~\ref{lem:gaussian_mechanism} and Lemma~\ref{lem:sense-U},  $\whP$ is $(\varepsilon/2, \delta/2)$-DP with probability at least $1-e^{-c_1p}-3n^{-9}-10^{-20\tilde{r}}$.  Finally, by the post-processing property of differentially private algorithm, we conclude that $\wtU\wtU^{\top}$ is an $(\varepsilon/2, \delta/2)$-DP estimator with the same probability. 
	
	Similarly, according to Lemma~\ref{lem:gaussian_mechanism} and Lemma~\ref{lem:sense-lambda}, $\wtLambda$ is $(\varepsilon/2, \delta/2)$-DP with probability at least $1-n^{-10}$ if we choose a large $\Delta_2$ as 
	$$
	\Delta_2=C_2\frac{\lambda(r+\log n)+\sigma^2(p+\log n)}{n}
	$$ 
	for some large absolute constant $C_2>0$.  This is due to the sensitivity $ \wtU^{\top}(\whSig-\sigma^2 I_p)\wtU$ is bounded by (note that $\wtU$ is already differentially private)
	\begin{align*}
		\big\| \wtU^{\top}(\whSig-\whSig^{(i)})\wtU\big\|_{\rm F}\leq \|\whSig-\whSig^{(i)}\|_{\rm F}\leq \Delta_2,
	\end{align*}
	where the last inequality holds in event $\mcE^{0}$ by Lemma~\ref{lem:con-X-norm}.   Then, the joint mechanism $\widetilde{J}=(\wtU, \wtLambda)$ is $(\varepsilon, \delta)$-DP with probability at least $1-4n^{-100}-e^{-c_1(n\wedge p)}$, by the post-processing property of differentially private algorithms. Therefore, we conclude that $\wtSig=\wtU\wtLambda \wtU^{\top}+\sigma^2I_p$ is $(\varepsilon, \delta)$-DP with the same probability.

	
\subsection{Proof of Lemma~\ref{lem:sense-U}}
	\label{appendix_lem:sense-U}
	
	Let the events $\mcE_0, \mcE_{\Delta}, \mcE_1$ be defined as in Section~\ref{sec:tech-lem}.  The following analysis proceeds mainly on the event $\mcE^{\ast}:=\mcE_0\cap \mcE_{\Delta}\cap \mcE_1$,  which occurs with probability at least $1-e^{-c_1(n\wedge p)}-3n^{-99}$. 
	
	On the event $\mcE^*$ and under the conditions that $\Sigma \in \Theta\left(p, r, \lambda, \sigma^2\right)$, $n\geq C_1(r\log n+\log^2n)$, $2r\leq p$ and $\lambda/\sigma^2\geq C_1\big(p/n+\sqrt{p/n}\big)$ for a large absolute constant $C_1>0$, we have  
	$
	\lambda_r \geq 5 \left(\|\Delta\|\vee \|\Delta^{(i)}\| \right). 
	$ 
	Therefore, we are able to apply Lemma \ref{lem:spectral-formula} to obtain 
	$$
	\whU\whU^{\top} - UU^{\top} = \sum_{k \geq 1} \mathcal{S}_{\Sigma, k}(\Delta)\quad {\rm and}\quad \whU^{(i)}\whU^{(i)\top} - UU^{\top} = \sum_{k \geq 1} \mathcal{S}_{\Sigma, k}(\Delta^{(i)}).  
	$$ 
	The explicit formula of spectral projectors implies that
	\begin{equation}\label{first_oder_term_plus_higher_order_term}
		\begin{aligned}
			\|\whU\whU^{\top} -\whU^{(i)}\whU^{(i)\top} \|_{\rm F} &= \bigg\| \sum_{k \geq 1} \mathcal{S}_{\Sigma, k}(\Delta) - \sum_{k \geq 1} \mathcal{S}_{\Sigma, k}(\Delta^{(i)}) \bigg\|_{\rm F} \\
			& \leq \| \mathcal{S}_{\Sigma, 1}(\Delta) - \mathcal{S}_{\Sigma, 1}(\Delta^{(i)}) \|_{\rm F} + \bigg\| \sum_{k \geq 2} \mathcal{S}_{\Sigma, k}(\Delta) - \sum_{k \geq 2} \mathcal{S}_{\Sigma, k}(\Delta^{(i)}) \bigg\|_{\rm F}. 
		\end{aligned}
	\end{equation}
	We now bound the first-order term $\| \mathcal{S}_{\Sigma, 1}(\Delta) - \mathcal{S}_{\Sigma, 1}(\Delta^{(i)}) \|_{\rm F} $ and the higher order term  $ \| \sum_{k \geq 2} \mathcal{S}_{\Sigma, k}(\Delta) - \sum_{k \geq 2} \mathcal{S}_{\Sigma, k}(\Delta^{(i)}) \|_{\rm F}$, separately. 
	
	\vspace{0.2cm}
	
	\noindent\emph{Step 1: bounding the first order term.} By the definitions of $ \mathcal{S}_{\Sigma, 1}(\Delta)$ and $ \mathcal{S}_{\Sigma, 1}(\Delta^{(i)}) $, 
	\begin{align}
		\max_{i\in[n]}\| \mathcal{S}_{\Sigma, 1}(\Delta) - \mathcal{S}_{\Sigma, 1}(\Delta^{(i)}) \|_{\rm F} \leq& \max_{i\in[n]}\| Q^{-1} (\Delta - \Delta^{(i)})  Q^{\perp}\|_{\rm F} + \|Q^{\perp} (\Delta - \Delta^{(i)}) Q^{-1} \|_{\rm F}\notag\\
		\leq&\frac{2}{n}\max_{i\in[n]}\Big(\|\Lambda^{-1}U^{\top}X_iX_i^{\top}U_{\perp}\|+\|\Lambda^{-1}U^{\top}X_i'X_i'^{\top}U_{\perp}\|\Big)\notag\\
		\leq& C_3\sqrt{\frac{\sigma^2(\lambda+\sigma^2)}{\lambda^2}}\cdot \frac{\sqrt{p(r+\log n)}}{n},\label{eq:proof-S-first}
	\end{align}
	where the last inequality holds on $\mcE^{\ast}$ for all $i\in[n]$ based on Lemma~\ref{lem:con-bounds}.

	\vspace{0.2cm}
	
	\noindent\emph{Step 2: bounding the higher-order terms}.  	Let $I_k$ be the index set for terms in $\mcS_{\Sigma, k}$
	$$
	I_k = \bigg\{\bfs: \mathbf{s}=\left(s_1, \ldots, s_{k+1}\right), \sum_{m=1}^{k+1} s_m = k, s_m\geq 0, \quad \forall m\in[k+1] \bigg\}, 
	$$
	with the cardinality  $|I_k| = \binom{2k}{k}$. We define
	\begin{align*}
		\mathcal{T}_{\Sigma, k, \bfs, l}(\Delta - \Delta^{(i)}) := Q^{-s_1}\Delta^{(i)} Q^{-s_2} \cdots Q^{-s_{l}} (\Delta - \Delta^{(i)})Q^{s_{l+1}} \cdots Q^{-s_{k}}\Delta Q^{s_{k+1}}, 
	\end{align*}
	for $k\geq 2$, $\bfs = (s_1, \cdots, s_{k+1})\in I_k$ and $l\in[k]$. Since  $|I_k| = \binom{2k}{k}$, the higher order terms can be bounded as follows 
	\begin{equation}\label{higher_order_in_terms_of_t}
		\begin{aligned}
			\Big\|\sum_{k \geq 2} \mathcal{S}_{\Sigma, k}(\Delta) - \sum_{k \geq 2} \mathcal{S}_{\Sigma, k}(\Delta^{(i)})\Big\|_{\rm F} 
			& =\Big\| \sum_{k \geq 2} \sum_{\bfs \in I_k} \sum_{l\in[k]}  \mcT_{\Sigma, k, \bfs, l}(\Delta - \Delta^{(i)})\Big\|_{\rm F} \\ 
			& \leq \sum_{k \geq 2}  \binom{2k}{k}  \max_{\bfs\in I_k}\sum_{l\in[k]} \| \mcT_{\Sigma, k, \bfs, l}(\Delta - \Delta^{(i)})\|_{\rm F}. 
		\end{aligned} 
	\end{equation}
	It suffices to upper bound $\|\mcT_{\Sigma,k,\bfs,l}(\Delta-\Delta^{(i)})\|_{\rm F}$ for any $\bfs\in I_k$. Denote
	$$
	D_{\max}:=C_2\sqrt{\frac{(\lambda+\sigma^2)(r\lambda+p\sigma^2)}{n}}+\frac{\lambda_r}{10}, 
	$$
	the upper bound appeared in the event $\mcE_{\Delta}$ so that $\|\Delta\|+\max_{i\in[n]}\|\Delta^{(i)}\|\leq D_{\max}$ in the event $\mcE_{\Delta}$.   
	
	\begin{lemma}
		\label{lem:T-highterm}
		Under the conditions of Lemma~\ref{lem:sense-U}, the following bound holds in the event $\mcE^{\ast}$ 	for all $k\geq 2$ and $\bfs \in I_k$.
		\begin{align*}
			&  \sum_{l\in[k]}\| \mathcal{T}_{\Sigma, k, \bfs, l}(\Delta - \Delta^{(i)}) \|_{\rm F} \leq C_6 \frac{(3+k)k}{2} \left(\frac{D_{\max}}{\lambda_r} \right)^{k-2}\cdot \frac{\sigma^2}{\lambda}\Big(\sqrt{\frac{p}{n}}+\frac{p}{n}\Big) \cdot\sqrt{\frac{\sigma^2 (\lambda+\sigma^2) }{\lambda^2}} \frac{\sqrt{p(r+\log n)}}{n}, 
		\end{align*}
		where $C_6>0$ is an absolute constant. 
	\end{lemma}

	We now continue from (\ref{higher_order_in_terms_of_t}) and get for all $i\in[n]$
	\begin{align}
		\Big\|\sum_{k \geq 2} \mathcal{S}_{\Sigma, k}&(\Delta) - \sum_{k \geq 2} \mathcal{S}_{\Sigma, k}(\Delta^{(i)})\Big\|_{\rm F} 
		\leq \sum_{k\geq 2} {2k\choose k}  \max_{\bfs \in I_k}\sum_{l\in[k]}\| \mathcal{T}_{\Sigma, k, \bfs, l}(\Delta - \Delta^{(i)}) \|_{\rm F} \notag\\
		&\lesssim \sum_{k\geq 2} {2k \choose k} \frac{(3+k)k}{2} \left(\frac{D_{\max}}{\lambda_r} \right)^{k-2}\cdot \frac{\sigma^2}{\lambda}\Big(\sqrt{\frac{p}{n}}+\frac{p}{n}\Big) \cdot\sqrt{\frac{\sigma^2 (\lambda+\sigma^2) }{\lambda^2}} \frac{\sqrt{p(r+\log n)}}{n}\notag\\
		&\lesssim \frac{\sigma^2}{\lambda}\Big(\sqrt{\frac{p}{n}}+\frac{p}{n}\Big) \cdot\sqrt{\frac{\sigma^2 (\lambda+\sigma^2) }{\lambda^2}} \frac{\sqrt{p(r+\log n)}}{n},\label{eq:sum-Sk-higher}
	\end{align}
	where the last inequality holds if $\lambda\geq C_7D_{\max}$ for a large enough constant $C_7>0$. 
	
	Combining (\ref{eq:proof-S-first}) and (\ref{eq:sum-Sk-higher}) together with the condition $\lambda/\sigma^2\geq C_3(\sqrt{p/n}+p/n)$,  we get in event $\mcE^{\ast}$ that
	\begin{align*}
		\max_{i\in[n]}\|\whU\whU^{\top} -\whU^{(i)}\whU^{(i)\top} \|_{\rm F} \lesssim \sqrt{\frac{\sigma^2 (\lambda+ \sigma^2)}{\lambda^2}}\frac{\sqrt{p(r+\log n)}}{n}. 
	\end{align*}

	\subsection{Proof of Lemma~\ref{lem:sense-lambda}}
	\label{sec: proofs_DP-SPIKED-COV_is_DP}
	
	Recall that,  by definition of $\Delta$ and $\Delta^{(i)}$,  we can write
	$$
	\whSig^{(i)}=\whSig-\frac{1}{n}X_iX_i^{\top}+\frac{1}{n}X_i'X_i'^{\top}. 
	$$
	By Hoffman-Weilandt's inequality,  we have 
	\begin{align*}
		\sum_{k=1}^p \big(\lambda_k(\whSig)-\lambda_k(\whSig^{(i)})\big)^2\leq& \|\whSig-\whSig^{(i)}\|_{\rm F}^2
		\leq 2\|\whSig-\whSig^{(i)}\|^2\leq \frac{4}{n^2}\Big(\|X_i\|^2+\|X_i'\|^2\Big)^2\\
		\leq& C_4\bigg(\frac{\lambda(r+\log n)+(p+\log n)\sigma^2}{n}\bigg)^2,
	\end{align*}
	where the last inequality holds in the event $\mcE_0$ defined in Lemma~\ref{lem:con-X-norm}.   This completes the proof.

	\subsection{Proof of Lemma~\ref{lem:spectral-formula}}
		
\begin{figure}[h!]
			\centering
			\includegraphics[scale=0.5]{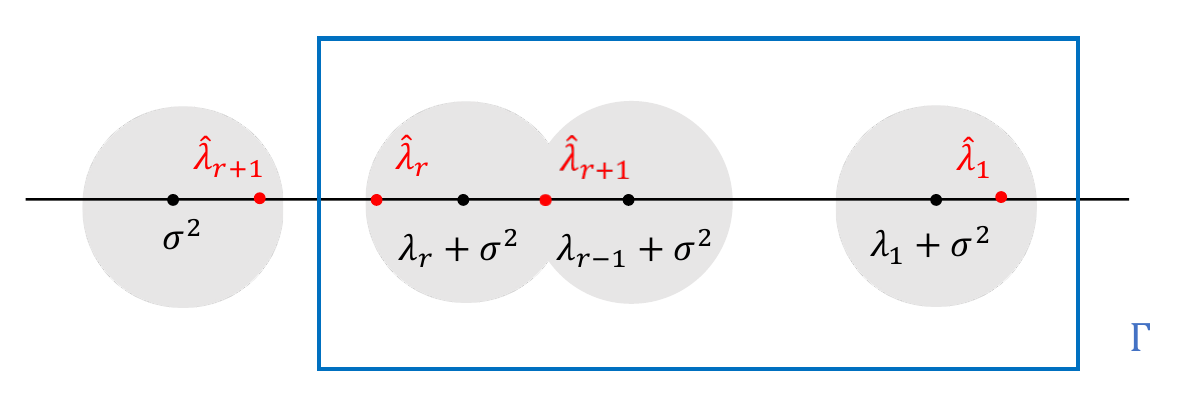}
			\caption{The contour plot of $\Gamma$. All grey balls share the same radius $\|\Delta\|$. }
			\label{contour}
		\end{figure}

		The proof is in spirit similar to the proof for the representation formula in \cite{xia2021normal}. Let $\{\tilde\lambda_i, u_i\}_{i\in[p]}$ be the the singular values and singular vectors of $\Sigma$ where $\tilde\lambda_i = \lambda_i + \sigma^2$ for $i\in[r]$ and $\tilde\lambda_i = \sigma^2$ for $i>r$. Let $\{\widehat{\lambda}_i, \whu_i\}_{i=1}^{p}$ denote the singular values and singular vectors of $\whSig$.  Due to Weyl's Lemma, for all $i\in[p]$,
		$$
		|\hat\lambda_i-\tilde\lambda_i|\leq \|\whSig - \Sigma\| = \|\Delta\|, 
		$$
		and thus $\whl_i$ must lie within $I(\tilde\lambda_i, \|\Delta\|) :=\left[\tilde\lambda_i-\|\Delta\|, \tilde\lambda_i+\|\Delta\|\right]$, a closed interval centered at $\tilde\lambda_i$ with half width $\|\Delta\|$. Under the condition that $\tilde\lambda_r-\tilde\lambda_{r+1}> 2\|\Delta\|$, i.e. $\frac{\lambda_r}{2}\geq \|\Delta\|$, we have $I(\tilde\lambda_r, \|\Delta\|) \cap I(\tilde\lambda_{r+1}, \|\Delta\|) = \emptyset $ and therefore, there exist a contour $\Gamma$ (see Figure \ref{contour}) on the complex plane such that 
		$$
		\{\tilde\lambda_i\}_{i\in[r]}\cup \{\whl_i\}_{i\in[r]} \subset  \Gamma_{D}
		,$$ 
		but 
		$$
		\{\tilde\lambda_i\}_{i\in[p]\setminus[r]}\cup \{\whl_i\}_{i\in[p]\setminus[r]} \subset \Gamma_{D}^{\complement}, 
		$$
		where $\Gamma_{D}$ is the open region enclosed by $\Gamma$, i.e. $\Gamma = \partial \Gamma_{D}$. 
		By Cauchy's integral formula, 
		\begin{align*}
			& \frac{1}{2\pi i}\oint_{\Gamma}(\eta I-\hat \Sigma)^{-1}d\eta\\ 
			&=\sum_{i=1}^r\frac{1}{2\pi i}\oint_{\Gamma}\frac{d\eta}{\eta-\hat \lambda_i}(\whu_i \whu_i^{\top })+\sum_{i=r+1}^p\frac{1}{2\pi i}\oint_{\Gamma}\frac{d\eta}{\eta-\hat \lambda_i}(\whu_i \whu_i^{\top })\\
			&=\sum_{i=1}^r \whu_i\whu_i^{\top}=\whU\whU^{\top}.
		\end{align*}
		As a result, we have
		\begin{equation}\label{hat_U_proj}
			\whU\whU^{\top }=\frac{1}{2\pi i}\oint_{\Gamma}(\eta I-\hat \Sigma)^{-1}d\eta, 
		\end{equation}
		and similarly, 
		\begin{equation}\label{U_proj}
			UU^{\top} = \frac{1}{2\pi i}\oint_{\Gamma}(\eta I- \Sigma)^{-1}d\eta.
		\end{equation}
		We denote $\mcR_{\Sigma}(\eta):=(\eta I- \Sigma)^{-1}$ and thus 
		\begin{align*}
			(\eta I-\hat \Sigma)^{-1}=(\eta I - \Sigma-  \Delta )^{-1}=&\big[(\eta I- \Sigma)\big(I-\mcR_{\Sigma}(\eta) \Delta \big)\big]^{-1}\\
			=&\big(I-\mcR_{\Sigma}(\eta) \Delta \big)^{-1}\mcR_{\Sigma}(\eta). 
		\end{align*}
		Since 
		$$
		\big\|\mcR_{\Sigma}(\eta) \Delta  \big\|\leq \|\mcR_{\Sigma}(\eta)\|\| \Delta \|\leq \frac{2\| \Delta \|}{\tilde\lambda_r}< 1.
		$$
		the Neumann series of $\big(I-\mcR_{\Sigma}(\eta) \Delta \big)^{-1}$ is 
		\begin{equation}\label{eq:neumann}
			\big(I-\mcR_{\Sigma}(\eta) \Delta \big)^{-1}=I+\sum_{k\geq 1}[\mcR_{\Sigma}(\eta) \Delta ]^k.
		\end{equation}
		Plugging (\ref{eq:neumann}) and  $\mcR_{\Sigma}(\eta)=(\eta I- \Sigma)^{-1}$ into (\ref{hat_U_proj}), we have
		\begin{align*}
			\whU\whU^{\top }=&\frac{1}{2\pi i}\oint_{\Gamma}(\eta I-\hat \Sigma)^{-1}d\eta\\
			=&\frac{1}{2\pi i}\oint_{\Gamma}\mcR_{\Sigma}(\eta) d\eta +\sum_{k\geq 1}\frac{1}{2\pi i}\oint_{\Gamma}\big[\mcR_{\Sigma}(\eta) \Delta \big]^k\mcR_{\Sigma}(\eta)d\eta \\
			=& UU^{\top} + \sum_{k\geq 1}\frac{1}{2\pi i}\oint_{\Gamma}\big[\mcR_{\Sigma}(\eta) \Delta \big]^k\mcR_{\Sigma}(\eta)d\eta. 
		\end{align*}
		Then we can write
		$$
		\whU\whU^{\top }-UU^{\top} =\sum_{k\geq 1}\frac{1}{2\pi i}\oint_{\Gamma}\big[\mcR_{\Sigma}(\eta) \Delta \big]^k\mcR_{\Sigma}(\eta)d\eta. 
		$$
		We denote the $k$-th order perturbation as 
		\begin{equation}\label{eq:calS_Mk}
			\mcS_{\Sigma,k}( \Delta ) := \frac{1}{2\pi i}\oint_{\Gamma}\big[\mcR_{\Sigma}(\eta) \Delta \big]^k\mcR_{\Sigma}(\eta)d\eta
		\end{equation}
		for all integer $k\geq 1$ and hence obtain 
		\begin{equation}\label{eq:hatTheta-Theta}
			\whU\whU^{\top }-UU^{\top}=\sum_{k\geq 1}\mcS_{\Sigma,k}( \Delta ).
		\end{equation}
		We derive the explicit formulas for $\mcS_{\Sigma,1}( \Delta )$ and $\mcS_{\Sigma,2}( \Delta )$ as an appetizer served before showing the explicit formulas for $\mcS_{\Sigma,k}( \Delta )$ with general integer $k\geq 1$. Note that 
		$$
		\mcR_{\Sigma}(\eta)=\sum_{j=1}^{p}\frac{1}{\eta-\tilde\lambda_j}u_ju_j^{\top}, 
		$$
		where $\tilde\lambda_i = \lambda_i + \sigma^2$ for $i\in[r]$ and $\tilde\lambda_i = \sigma^2$, for $i\in[p]\setminus[r]$. Let $P_j=u_ju_j^{\top }$ be the spectral projector onto $u_j$, for all $j \in[p]$.

		\paragraph{Derivation of $\mcS_{\Sigma,1}( \Delta )$.}
		By the definition of $\mcS_{\Sigma,1}( \Delta )$,
		\begin{align}
			\mcS_{\Sigma,1}( \Delta )=&\frac{1}{2\pi i}\oint_{\Gamma}\mcR_{\Sigma}(\eta) \Delta \mcR_{\Sigma}(\eta)d\eta\nonumber\\
			=&\sum_{j_1=1}^{p}\sum_{j_2=1}^{p}\frac{1}{2\pi i}\oint_{\Gamma}\frac{d\eta}{(\eta-\tilde\lambda_{j_1})(\eta-\tilde\lambda_{j_2})}P_{j_1} \Delta P_{j_2}.\label{eq:lemSM1_eq1}
		\end{align}
		{\it Case 1:} When both $j_1$ and $j_2$ are greater than $r$, the contour integral in (\ref{eq:lemSM1_eq1}) is zero according to Cauchy's integral formula.\\
		{\it Case 2}: When only one of $j_1$ and $j_2$ is greater than $r$. W.L.O.G, we assume $j_2>r$, then  
		\begin{align*}
			& \sum_{j_1=1}^r\sum_{j_2>r}^p\frac{1}{2\pi i}\oint_{\Gamma} \frac{d\eta}{(\eta-\tilde\lambda_{j_1})(\eta-\tilde\lambda_{j_2})}P_{j_1} \Delta P_{j_2} \\
			& =\sum_{j_1=1}^r\sum_{j_2>r}(\tilde\lambda_{j_1}-\tilde\lambda_{j_2})^{-1}P_{j_1} \Delta P_{j_2} \\
			& = \sum_{j_1=1}^r\sum_{j_2>r}(\lambda_{j_1})^{-1}P_{j_1} \Delta P_{j_2} \\ 
			& =Q^{-1} \Delta Q^{\perp}, 
		\end{align*} 
		where and $Q^{-1} = U\Lambda^{-1} U^{\top}$ and $Q^{\perp}=U_{\perp}U_{\perp}^{\top}$. \\
		{\it Case 3}: When none of $j_1$ or $j_2$ is greater than $r$, the contour integral in (\ref{eq:lemSM1_eq1}) is zero. 
		
		In summary, $\mcS_{\Sigma,1}( \Delta )=Q^{-1} \Delta Q^{\perp}+Q^{\perp} \Delta Q^{-1}$.

		\paragraph{Derivation of $\mcS_{\Sigma,2}( \Delta )$.}
		By the definition of $\mcS_{\Sigma,2}( \Delta )$, 
		\begin{align}
			\mcS_{\Sigma,2}( \Delta )=&\frac{1}{2\pi i}\oint_{\Gamma}\mcR_{\Sigma}(\eta) \Delta \mcR_{\Sigma}(\eta) \Delta \mcR_{\Sigma}(\eta)d\eta\nonumber\\
			=&\sum_{j_1=1}^{p}\sum_{j_2=1}^{p}\sum_{j_3=1}^{p}\frac{1}{2\pi i}\oint_{\Gamma}\frac{d\eta}{(\eta-\tilde\lambda_{j_1})(\eta-\tilde\lambda_{j_2})(\eta-\tilde\lambda_{j_3})}P_{j_1} \Delta P_{j_2} \Delta P_{j_3}.\label{eq:lemSM2_eq1}
		\end{align}
		{\it Case 1}: When all $j_1, j_2, j_3$ are greater than $r$, the contour integral in (\ref{eq:lemSM2_eq1}) is zero by Cauchy's integral formula. \\
		{\it Case 2}: When two of $j_1,j_2,j_3$ are greater than $r$. W.L.O.G., we assume $j_1\leq r$ and $j_2,j_3>r$, then
		\begin{align*}
			& \sum_{j_1=1}^{r}\sum_{j_2,j_3>r}^{p}\frac{1}{2\pi i }\oint_{\Gamma}\frac{d\eta}{(\eta-\tilde\lambda_{j_1})(\eta-\tilde\lambda_{j_2})(\eta-\tilde\lambda_{j_3})}P_{j_1} \Delta P_{j_2} \Delta P_{j_3}\\
			=&\sum_{j_1=1}^{r}\sum_{j_2,j_3>r}^{p}\frac{1}{(\tilde\lambda_{j_1}-\tilde\lambda_{j_2})(\tilde\lambda_{j_1}-\tilde\lambda_{j_3})}P_{j_1} \Delta P_{j_2} \Delta P_{j_3}=Q^{-2} \Delta Q^{\perp} \Delta Q^{\perp} \\
			=&\sum_{j_1=1}^{r}\sum_{j_2,j_3>r}^{p}\frac{1}{\lambda_{j_1}^2}P_{j_1} \Delta P_{j_2} \Delta P_{j_3}=Q^{-2} \Delta Q^{\perp} \Delta Q^{\perp}. 
		\end{align*}
		{\it Case 3}: one of $j_1,j_2,j_3$ is greater than $r$. W.L.O.G., let $j_1,j_2\leq r$ and $j_3>r$, we get
		\begin{align*}
			&\sum_{j_1,j_2=1}^{r}\sum_{j_3>r}^{p}\frac{1}{2\pi i}\oint_{\Gamma}\frac{d\eta}{(\eta-\tilde\lambda_{j_1})(\eta-\tilde\lambda_{j_2})(\eta-\tilde\lambda_{j_3})}P_{j_1} \Delta P_{j_2} \Delta P_{j_3}\\
			=&\sum_{j_1=j_2=1}^{r}\sum_{j_3>r}^{p}\frac{1}{2\pi i}\oint_{\Gamma}\frac{d\eta}{(\eta-\tilde\lambda_{j_1})^2 (\eta-\tilde\lambda_{j_3})}P_{j_1} \Delta P_{j_1} \Delta P_{j_3}\\
			&\quad +\sum_{j_1\neq j_2\geq 1}^{r}\sum_{j_3>r}^{p}\frac{1}{2\pi i}\oint_{\Gamma}\frac{ d\eta}{(\eta-\tilde\lambda_{j_1})(\eta-\tilde\lambda_{j_2})(\eta-\tilde\lambda_{j_3})}P_{j_1} \Delta P_{j_2} \Delta P_{j_3}\\
			=&-\sum_{j_1=1}^r(\tilde\lambda_{j_1}-\tilde\lambda_{j_3})^{-2}P_{j_1} \Delta P_{j_1} \Delta Q^{\perp}-\sum_{j_1\neq j_2\geq 1}^{r}(\tilde\lambda_{j_1}-\tilde\lambda_{j_3})^{-1}(\tilde\lambda_{j_2}-\tilde\lambda_{j_3})^{-1}P_{j_1} \Delta P_{j_2} \Delta Q^{\perp}\\
			=&-\sum_{j_1=1}^r\lambda_{j_1}^{-2}P_{j_1} \Delta P_{j_1} \Delta Q^{\perp}-\sum_{j_1\neq j_2\geq 1}^{r}\lambda_{j_1}^{-1}\lambda_{j_2}^{-1}P_{j_1} \Delta P_{j_2} \Delta Q^{\perp}\\
			=&-Q^{-1} \Delta Q^{-1} \Delta Q^{\perp}.
		\end{align*}
		{\it Case 4:} When none of $j_1,j_2,j_3$ is greater than $r$, the contour integral in (\ref{eq:lemSM2_eq1}) is zero. 
		
		In summary, 
		\begin{align*}
			\mcS_{\Sigma,2}( \Delta )=&\big(Q^{-2} \Delta Q^{\perp} \Delta Q^{\perp}+Q^{\perp} \Delta Q^{-2} \Delta Q^{\perp}+Q^{\perp} \Delta Q^{\perp} \Delta Q^{-2}\big)\\
			-&\big(Q^{\perp} \Delta Q^{-1} \Delta Q^{-1}+Q^{-1} \Delta Q^{\perp} \Delta Q^{-1}+Q^{-1} \Delta Q^{-1} \Delta Q^{\perp}\big).
		\end{align*}
		
		\paragraph{Derivation of $\mcS_{\Sigma,k}( \Delta )$ for general k.}
		By the definition of $\mcS_{\Sigma,k}( \Delta )$, 
		\begin{align}
			\mcS_{\Sigma,k}( \Delta )=\sum_{j_1,\cdots,j_{k+1}\geq 1}^p\frac{1}{2\pi i}\oint_{\Gamma}\Big(\prod_{i=1}^{k+1}\frac{1}{\eta-\tilde\lambda_{j_i}}\Big)d\eta P_{j_1} \Delta P_{j_2} \Delta \cdots P_{j_k} \Delta P_{j_{k+1}}.\label{eq:lemSMk_eq1}
		\end{align}
		In order to deal with each component in the summations (\ref{eq:lemSMk_eq1}), we first consider some special cases for ease of the notation. W.L.O.G., we consider the cases where some $\bar{k}$ indices from $\{j_1,\cdots,j_{k+1}\}$ are not larger than $r$. More specifically, we restrict our discussion to the case where $j_1,\cdots,j_{\bar{k}}\leq r$ and $j_{\bar{k}+1},\cdots,j_{k+1}>r$. By Cauchy integral formula, the integral in (\ref{eq:lemSMk_eq1}) is zero once $\bar{k}=0$ or $\bar{k}=k+1$ and thus we focus on the non-trivial case $\bar{k} \in [k-1]$.
		Note that 
		\begin{align*}
			& \sum_{j_1,\cdots,j_{\bar{k}}\geq 1}^r \sum_{j_{\bar{k}+1},\cdots,j_{k+1}>r}^p\frac{1}{2\pi i}\oint_{\Gamma}\Big(\prod_{i=1}^{p}\frac{1}{\eta-\tilde\lambda_{j_i}}\Big) d\eta P_{j_1} \Delta P_{j_2} \Delta \cdots P_{j_k} \Delta P_{j_{k+1}}\\
			& = \sum_{j_1,\cdots,j_{\bar{k}}\geq 1}^r \sum_{j_{\bar{k}+1},\cdots,j_{k+1}>r}^p\frac{1}{2\pi i}\oint_{\Gamma}\Big(\prod_{i=1}^{\bar{k}}\frac{1}{\eta-(\lambda_{j_i}+\sigma^2)}\Big)(\eta-\sigma^2)^{\bar{k}-k-1}d\eta P_{j_1} \Delta P_{j_2} \Delta \cdots P_{j_k} \Delta P_{j_{k+1}}\\
			& = \sum_{j_1,\cdots,j_{\bar{k}}\geq 1}^r\frac{1}{2\pi i}\oint_{\Gamma}\Big(\prod_{i=1}^{\bar{k}}\frac{1}{\eta-(\lambda_{j_i}+\sigma^2)}\Big)(\eta-\sigma^2)^{\bar{k}-k-1}d\eta P_{j_1} \Delta P_{j_2} \Delta \cdots P_{j_{\bar{k}}} \Delta Q^{\perp} \Delta \cdots  \Delta Q^{\perp}.
		\end{align*}
		Recall that our goal is to prove
		$$
		\mcS_{\Sigma,k}( \Delta )=\sum_{\bfs:s_1+\cdots+s_{k+1}=k}(-1)^{1+\tau(\bfs)}\cdot Q^{-s_1} \Delta Q^{-s_2} \Delta \cdots  \Delta Q^{-s_{k+1}}.
		$$
		Accordingly, in the above summations, we consider the components, where $s_1,\cdots,s_{\bar{k}}\geq 1$ and $s_{\bar{k}+1}=\cdots=s_{k+1}=0$, namely, 
		$$
		\sum_{\substack{s_1+\cdots+s_{\bar{k}}=k\\ s_j\geq 1}}(-1)^{\bar{k}+1}Q^{-s_1} \Delta \cdots  \Delta Q^{-s_{\bar{k}}} \Delta Q^{\perp}\cdots  \Delta Q^{\perp}.
		$$
		It turns out that we need to prove
		\begin{align*}
			& \sum_{j_1,\cdots,j_{\bar{k}}\geq 1}^r\frac{1}{2\pi i}\oint_{\Gamma}\Big(\prod_{i=1}^{\bar{k}}\frac{1}{\eta-(\lambda_{j_i}+\sigma^2)}\Big)(\eta-\sigma^2)^{\bar{k}-k-1}d\eta P_{j_1} \Delta P_{j_2} \Delta \cdots P_{j_{\bar{k}}}\nonumber\\
			& = \sum_{j_1,\cdots,j_{\bar{k}}\geq 1}^r\sum_{\substack{s_1+\cdots+s_{\bar{k}}=k\\ s_j\geq 1}}(-1)^{\bar{k}+1}\frac{1}{(\lambda_{j_1})^{s_1}\cdots (\lambda_{j_{\bar{k}}})^{s_{\bar{k}}}} P_{j_1} \Delta P_{j_2} \Delta \cdots  \Delta P_{j_{\bar{k}}}.
		\end{align*}
		It suffices to prove that for all $\bfj=(j_1,\dots,j_{\bar{k}})\in\{1,\cdots,r\}^{\bar{k}}$, 
		\begin{equation}
			\begin{aligned}
				& \frac{1}{2\pi i}\oint_{\Gamma}\Big(\prod_{i=1}^{\bar{k}}\frac{1}{\eta-(\lambda_{j_i}+\sigma^2)}\Big)(\eta-\sigma^2)^{\bar{k}-k-1}d\eta \\
				&=\sum_{\substack{s_1+\cdots+s_{\bar{k}}=k\\ s_j\geq 1}}(-1)^{\bar{k}+1}\frac{1}{(\lambda_{j_1})^{s_1}\cdots (\lambda_{j_{\bar{k}}})^{s_{\bar{k}}}} . \label{eq:lemSMk_eq2}
			\end{aligned}
		\end{equation}

		To prove (\ref{eq:lemSMk_eq2}), we rewrite its right hand side. Given any $\bfj=(j_1,\cdots,j_{\bar{k}})\in\{1,\cdots,r\}^{\bar{k}}$, we define
		$$
		V_i(\bfj):=\big\{1\leq t\leq \bar{k}: j_t=i\big\} \quad \forall i\in[r], 
		$$
		as a set that contains all location $l\in[k+1]$ such that $\lambda_{j_l}=\lambda_i$. For simplicity, we also denote $v_i(\bfj)=\left|V_i(\bfj)\right|$. Then, the right hand side of (\ref{eq:lemSMk_eq2}) is written as
		\begin{align*}
			& \sum_{\substack{s_1+\cdots+s_{\bar{k}}=k\\ s_j\geq 1}}(-1)^{\bar{k}+1}\frac{1}{(\lambda_{j_1})^{s_1}\cdots (\lambda_{j_{\bar{k}}})^{s_{\bar{k}}}} \\ &=(-1)^{\bar{k}+1}\sum_{\substack{s_1+\cdots+s_{\bar{k}}=k\\s_j\geq 1}}\lambda_1^{-\sum_{p\in V_1(\bfj)}s_p}\cdots \lambda_r^{-\sum_{p\in V_r(\bfj)}s_p}.
		\end{align*}
		Now, we denote $t_i(\bfj)=\sum_{p\in V_i(\bfj)}s_p$ for all $i\in[r]$ and rewrite the above equation as 
		\begin{align*}
			& \sum_{\substack{s_1+\cdots+s_{\bar{k}}=k\\ s_j\geq 1}}(-1)^{\bar{k}+1}\frac{1}{\lambda_{j_1}^{s_1}\cdots \lambda_{j_{\bar{k}}}^{s_{\bar{k}}}}  \\ & =(-1)^{\bar{k}+1}\sum_{\substack{t_1(\bfj)+\cdots+t_r(\bfj)=k\\ t_i(\bfj)\geq v_i(\bfj)\\t_i(\bfj)=0\ {\rm if}\ v_i(\bfj)=0} }\prod_{i: v_i(\bfj)\geq 1}{ \binom{t_i(\bfj)-1}{v_i(\bfj)-1} }\lambda_i^{-t_i(\bfj)} \\
			& =(-1)^{\bar{k}+1}\sum_{\substack{t_1(\bfj)+\cdots+t_r(\bfj) \\ =k-\bar{k}\\ t_i(\bfj)=0 \;\text{if}\; v_i(\bfj)=0} }\prod_{i: v_i(\bfj)\geq 1}{ \binom{t_i(\bfj)+v_i(\bfj)-1}{v_i(\bfj)-1}\lambda_i^{-t_i(\bfj)-v_i(\bfj)}}
		\end{align*}
		where the last equality is due to the fact $v_1(\bfj)+\cdots+v_r(\bfj)=\bar{k}$. Similarly, the left hand side of (\ref{eq:lemSMk_eq2}) can be written as 
		\begin{align*}
			& \frac{1}{2\pi i}\oint_{\Gamma}\frac{d\eta}{(\eta-(\lambda_{j_1}+\sigma^2))\cdots(\eta-(\lambda_{j_{\bar{k}}}+\sigma^2))(\eta-\sigma^2)^{k+1-\bar{k}}}\\ &=\frac{1}{2\pi i}\oint_{\Gamma}\frac{d\eta}{(\eta-(\lambda_{1}+\sigma^2))^{v_1(\bfj)}\cdots(\eta-(\lambda_r+\sigma^2))^{v_r(\bfj)}(\eta-\sigma^2)^{k+1-\bar{k}}}.
		\end{align*}
		Therefore, in order to prove (\ref{eq:lemSMk_eq2}), it suffices to prove that for any $\bfj=(j_1,\cdots,j_{\bar{k}})$ the following equality holds
		\begin{equation}
			\begin{aligned}
				\label{aim}
				& \frac{1}{2\pi i}\oint_{\Gamma}\frac{d\eta}{(\eta-(\lambda_{1}+\sigma^2))^{v_1}\cdots(\eta-(\lambda_r+\sigma^2))^{v_r}(\eta-\sigma^2)^{k+1-\bar{k}}}\\ &=(-1)^{\bar{k}+1}\sum_{\substack{t_1+\cdots+t_r=k-\bar{k}\\t_i=0\ {\rm if}\ v_i=0} }\prod_{i: v_i\geq 1}^r{ \binom{t_i+v_i-1}{v_i-1} }\lambda_i^{-t_i-v_i}, 
			\end{aligned}
		\end{equation}

		where we omitted the index $\bfj$ in definitions of $v_i(\bfj)$ and $t_i(\bfj)$ without causing any confusions. The non-negative numbers $v_1+\cdots+v_{r}=\bar{k}$. Let 
		$$
		\varphi(\eta)=\frac{1}{(\eta-(\lambda_{1}+\sigma^2))^{v_1}\cdots(\eta-(\lambda_r+\sigma^2))^{v_r}(\eta-\sigma^2)^{k+1-\bar{k}}}, 
		$$
		be a function of $\eta$. According to the Residue theorem,
		$$
		\frac{1}{2\pi i}\oint_{\Gamma}\varphi(\eta)d\eta=-{\rm Res}(\varphi,\eta=\infty)-{\rm Res}(\varphi,\eta=\sigma^2).
		$$
		Note that ${\rm Res}(\varphi, \eta=\infty)=0$ and it suffices to calculate ${\rm Res}(\varphi,\eta=\sigma^2)$. Let $\gamma_{\sigma^2} $ be a contour plot around $\eta = \sigma^2$ such that none of $\{\lambda_i+\sigma^2\}_{i\in[r]}$ lies in $\partial \gamma_{\sigma^2} $. Then,
		\begin{align*}
			{\rm Res}(\varphi,\eta=\sigma^2)=\frac{1}{2\pi i}\oint_{\gamma_{\sigma^2} }\varphi(\eta)d\eta.
		\end{align*}
		By Cauchy integral formula for derivatives, we have
		\begin{align*}
			& {\rm Res}(\varphi,\eta=\sigma^2) \\
			& = \varphi^{(k-\bar{k})}(\eta)\Big|_{\eta = \sigma^2} \\
			&=\frac{1}{(k-\bar{k})!}\Big[\prod_{i: v_i\geq 1}^r(\eta-(\lambda_{i}+\sigma^2))^{-v_i}\Big]^{(k-\bar{k})}\Big|_{\eta=\sigma^2}
		\end{align*}
		where $\varphi(\eta)^{(k-\bar{k})}$ is the $k-\bar{k}$-th order differentiation of $\varphi(\eta)$. Further, by general Leibniz rule, 
		\begin{align*}
			& {\rm Res}(\varphi,\eta=0)\\
			& =\frac{1}{(k-\bar{k})!}\sum_{\substack{t_1+\cdots+t_r=k-\bar{k}\\ t_i=0\ {\rm if}\ v_i=0}} \frac{(k-\bar{k})!}{t_1!t_2!\cdots t_r!}\prod_{i: v_i\geq 1}^r\Big[(\eta-(\lambda_{i}+\sigma^2)^{-v_i}\Big]^{(t_i)}\Big|_{\eta=\sigma^2}\\
			&= (-1)^{k-\bar{k}}\sum_{\substack{t_1+\cdots+t_r=k-\bar{k}\\ t_i=0\ {\rm if}\ v_i=0}} \prod_{i:v_i\geq 1}^r\frac{v_i(v_i+1)\cdots(v_i+t_i-1)}{t_i!}(-\lambda_i)^{-v_i-t_i}\\
			&= (-1)^{k-\bar{k}}\sum_{\substack{t_1+\cdots+t_r=k-\bar{k}\\ t_i=0\ {\rm if}\ v_i=0}} \prod_{i:v_i\geq 1}^r{ \binom{t_i+v_i-1}{v_i-1} }(-\lambda_i)^{-v_i-t_i}\\
			&= (-1)^{2k-\bar{k}}\sum_{\substack{t_1+\cdots+t_r=k-\bar{k}\\ t_i=0\ {\rm if}\ v_i=0}} \prod_{i:v_i\geq 1}^r{\binom{t_i+v_i-1}{v_i-1} }(\lambda_i)^{-v_i-t_i}.
		\end{align*}
		Therefore, 
		$$
		\frac{1}{2\pi i}\oint_{\Gamma}\varphi(\eta)d\eta=(-1)^{\bar{k}+1}\sum_{\substack{t_1+\cdots+t_r=k-\bar{k}\\ t_i=0\ {\rm if}\ v_i=0}} \prod_{i:v_i\geq 1}^r{\binom{t_i+v_i-1}{v_i-1} }\lambda_i^{-v_i-t_i}
		$$
		which shows (\ref{aim}) is true and the proof is completed.

	\subsection{Proof of Lemma~\ref{lem:dp-fano}}
	\label{sec:dp-fano}
	
	Recall that the Kullback-Leibler divergence and total variation distance between two probability measures defined in the same measure space $(\Omega, \mcF) $ are defined by 
	$$
	{\rm TV}(\mu \| \nu) := \int_{\Omega} \log \left(\frac{\mu(d x)}{\nu(d x)}\right) \mu(d x)\quad{\rm and}\quad {\rm TV}(\mu,\nu) := \sup_{E\in \mcF} | \mu(E)-\nu(E)|.
	$$
	By definition,  for any $(\varepsilon, 0)$-DP algorithm $A$,  we have 
	\begin{equation*}
		\begin{aligned}
			\sup_{P\in \mcP} \EE _{A} \; \rho( A,\theta(P)) 
			&\geq \sup_{P \in \mcQ} \EE _{A} \; \rho( A,\theta(P))\\
			&\geq \frac{1}{N} \sum_{i=1}^N  \frac{\rho_0}{2}\EE _{A} \mathbb{I}\left(\rho\left( A , \theta(P_i)\right)>\frac{\rho_0}{2}\right). 
		\end{aligned}
	\end{equation*}
	Following the proof of generalized Fano's Lemma \citep{devroye1987course, yu1997assouad}, we further reduce the estimation problem into hypothesis testing. 		Let $h : \Theta \rightarrow \Theta_{\mcQ}$ be defined by $h (\theta)=\underset{\theta_i \in \Theta_{\mcQ}}{\operatorname{argmin}} \; \rho\left(\theta_i, \theta\right)$. 	For all $i, j \in [N]$, define
	$$
	p_{ji}:= \EE _{A, X\sim P_j} \mathbb{I} (h (A)= \theta(P_i) ) ,
	$$
	measuring the probability that the algorithm $A$ outputs an estimator closer to $\theta(P_i)$ using data actually sampled from distribution $P_j$. The total type-I error of testing $H_0: X\sim P_j$ versus $H_1: X\sim P\in \mcP\setminus P_j$ is 
	$$
	\beta_j := \EE_{A, X\sim P_j} \mathbb{I} (h ( A ) \ne \theta(P_j)) = \sum_{i\ne j} p_{ji}. 
	$$
	Observe that 
	\begin{equation*}
		\begin{aligned}
			& \frac{1}{N} \sum_{i=1}^N  \frac{\rho_0}{2}\EE _{A} \mathbb{I}\left(\rho\left( A , \theta(P_i)\right)>\frac{\rho_0}{2}\right) \geqslant \frac{\rho_0}{2 N} \sum_{i=1}^N  \EE_{A, X\sim P_i} \mathbb{I}(h \left( A \right) \neq \theta(P_i)) = \frac{\rho_0}{2 N} \sum_{i=1}^N \beta_i.
		\end{aligned}
	\end{equation*}	
	Therefore,
	\begin{align*}
		& \sup_{P\in \mcP} \EE _{A} \; \rho( A,\theta(P)) \geq  \frac{\rho_0}{2 N} \sum_{i=1}^N \beta_i. 
	\end{align*}
	It suffices to find a lower bound for $\sum_{i=1}^N \beta_i $. Towards that end, the following lemma is needed. Its proof is relegated to Section~\ref{sec:lem-bridge}
	
	\begin{lemma}\label{bridge}
		Let $\mu = \mu_1 \times \cdots \times \mu_n$ and $\nu = \nu_1 \times \cdots \times \nu_n $ be two probability measures defined on the same measurable space $(\Omega,  \mcF)$. The total variation distance is denoted by $t_k = \operatorname{\operatorname{TV}}(\mu_k, \nu_k)$ for all $k \in [n]$.
		Let $A: \Omega\mapsto \Theta$ be an $(\varepsilon,\delta)$-DP random algorithm. 	Then, for any subset $E\subset\Theta$, 
		\begin{align*}
			\EE_{A, X\sim \mu} \mathbb{I}\big(  A(X)  \in E\big) & \leq \exp\bigg(4 \varepsilon \sum_{k=1}^n t_k\bigg)	\left(\EE_{A, X\sim \nu}\mathbb{I}\big(  A(X)  \in E\big) + \frac{2\delta}{e^{\varepsilon}-1}\right), 
		\end{align*}
	\end{lemma}

	Now consider any pair $P_i, P_j \in \mcQ$ where $P_i = \mu_i^{(1)}\times \cdots, \times \mu_i^{(n)}$ and $P_j = \mu_j^{(1)}\times \cdots, \times  \mu_j^{(n)}$.  By Lemma \ref{bridge} and the fact $p_{ii}  + \sum_{j \ne i} p_{ij}   = 1$ , we have 
	\begin{align*}
		p_{ji} & \geqslant \exp \left(- 4 \varepsilon \sum_{k\in[n]} \operatorname{TV}\left(\mu_i^{(k)} , \mu_j^{(k)}\right) \right) p_{ii} - \frac{2\delta}{e^{\varepsilon}-1} \\
		& =\exp \left(- 4 \varepsilon \sum_{k\in[n]} \operatorname{TV}\left(\mu_i^{(k)} , \mu_j^{(k)}\right)  \right)\left(1 - \sum_{j \ne i} p_{ij}\right) -\frac{2\delta}{e^{\varepsilon}-1}\\
		& =\exp \left(- 4 \varepsilon \sum_{k\in[n]} \operatorname{TV}\left(\mu_i^{(k)} , \mu_j^{(k)}\right)  \right)\left(1 - \beta_i \right) -\frac{2\delta}{e^{\varepsilon}-1} \\
		& \geq \exp \left(- 4 \varepsilon t_0 \right)\left(1 - \beta_i\right)  -\frac{2 \delta}{e^{\varepsilon}-1}. 
	\end{align*}
	Taking summation over $i\in [N] \setminus \{ j \}$, we have
	\begin{equation*}
		\beta_j   = \sum_{i\ne j}p_{ji}  \geqslant \exp \left(- 4 \varepsilon t_0 \right) \left(N-1-\sum_{i \neq j} \beta_i\right) - \frac{2(N-1)\delta}{e^{\varepsilon}-1}  
	\end{equation*}
	Further summing over $j\in [N] $, we have
	\begin{equation*}
		\sum_{j=1}^N \beta_j  \geqslant \frac{N(N-1)}{(N-1)+\exp \left(  4 \varepsilon t_0 \right)}\cdot\left(1-\frac{2\delta e^{4\varepsilon t_0}}{e^{\varepsilon}-1}\right). 
	\end{equation*}
	Therefore, 
\begin{align*}
			\sup_{P\in \mcP} \EE _{A} \; \rho( A,\theta(P)) &\geq \frac{\rho_0}{2}\frac{(N-1)}{(N-1)+\exp \left( 4 \varepsilon t_0 \right)}\cdot \left(1-\frac{2\delta e^{4\varepsilon t_0}}{e^{\varepsilon}-1}\right)\\
			 &\geq \frac{\rho_0}{4} \left( 1 \wedge \frac{N-1}{\exp \left( 4 \varepsilon t_0 \right) } \right) \left(1-\frac{2\delta e^{4\varepsilon t_0}}{e^{\varepsilon}-1}\right). 
\end{align*}

	Combining the above lower bound with classical generalized Fano's Lemma (e.g.,  \cite{tsybakov2008intro}), we have 
	\begin{align*}
		\inf_{ A\in \mcA_{\varepsilon,\delta}(\mcP)} & \sup_{P\in\mcP } \EE_{A} \; \rho( A, \theta(P)) \geqslant \max \left\{\frac{\rho_0}{2}\left(1-\frac{ l_0 +\log 2}{\log N}\right), \frac{\rho_0}{4}\left(1 \wedge \frac{N-1}{\exp \left( 4 \varepsilon t_0 \right)}\right)\left(1-\frac{2\delta e^{4\varepsilon t_0}}{e^{\varepsilon}-1}\right)\right\}. 
	\end{align*}

	\subsection{Proof of Corollary~\ref{cor:sub-Gaussian}}
	It suffices to establish the concentration bounds as in Lemmas~\ref{lem:con-X-norm} and \ref{lem:con-bounds}.  Note that the while Lemma~\ref{lem:con-hatSigma} is stated in Gaussian distributions,  but the claimed bounds still hold for sub-Gaussian distributions.  See \cite{koltchinskii2017concentration} for more details.  
	
	As a result,  it  is easy to check that the upper bounds of $\max_{i\in[n]}\|X_i\|,  \max_{i\in [n]}\|U^{\top}X_i\|$,  and $\max_{i\in[n]}\|U_{\perp}^{\top}X_i\|^2$ stated in Lemma~\ref{lem:con-X-norm} still hold for sub-Gaussian distribution.  Similarly,  the upper bound of $\|\Delta\|$ in Lemma~\ref{lem:con-bounds} still holds.  However,  we shall pay specific attentions to the bounds in event $\mcE_1$.   
	
	\begin{lemma}\label{lem:subG-bounds}
		Suppose that $X_1,X_1',  \cdots,  X_n,  X_n'$  follow sub-Gaussian distribution with the proxy-covariance matrix $\Sigma\in\Theta(\lambda,\sigma^2)$,  $n\geq C_1\big(r\log (p+n)\log^2r+\log^2n\big)$,  $2r\leq p$,  and $\lambda/\sigma^2\geq C_1p/n$ for some absolute constant $C_1>0$.  There exist absolute constants $c_2,  C_3>0$ such that the event 
		\begin{align*}
			\mcE_1':=\Bigg\{\big\|U^{\top}&\Delta U_{\perp} \big\|+\max_{i\in[n]}\|U^{\top}\Delta^{(i)}U_{\perp}\|\leq C_3\sqrt{\frac{\sigma^2(\lambda+\sigma^2)p\log(p+n)}{n}} \Bigg\}\\
			&\bigcap \Bigg\{\max_{i\in[n]}\big\|U^{\top}(X_iX_i^{\top}/n)U_{\perp}\big\|+\big\|U^{\top}(X_i'X_i'^{\top}/n)U_{\perp}\big\|\leq C_3\frac{\sqrt{\sigma^2(\lambda+\sigma^2)p(r+\log n)}}{n} \Bigg\}\\
			&\bigcap \Bigg\{\max_{i\in[n]}\Big\|U^{\top}\Big(\frac{1}{n}\sum_{j\neq i}X_jX_j^{\top}\Big)U_{\perp}U_{\perp}^{\top}X_i \Big\|\leq C_3\sigma\cdot \sqrt{\frac{\sigma^2(\lambda+\sigma^2)p(r+\log n)\log^2(p+n)}{n}} \Bigg\}, 
		\end{align*}
		holds with probability $\PP(\mcE_1')\geq 1-2n^{-100}$. Meanwhile,
		$$
		\EE \|U^{\top}\Delta U_{\perp}\|\leq C_3\sqrt{\frac{\sigma^2(\lambda+\sigma^2)p\log p}{n}}.
		$$
	\end{lemma}
	We then study the sensitivity of eigenvectors and eigenvalues by bounding $\|\wtU\wtU^{\top}-\wtU^{(i)}\wtU^{(i)\top}\|_{\rm F}$ and $\sum_{k=1}^p \big(\lambda_k(\whSig)-\lambda_k(\whSig^{(i)})\big)^2$,  whose proofs are similar to those of Lemmas~\ref{lem:sense-U} and ~\ref{lem:sense-lambda}.  We can show that  if  $\lambda/\sigma^2\geq C_3\big(\sqrt{p/n}+p/n\big)\log(p+n)$,  then in the event $\mcE^{\ast}:=\mcE_0\cap \mcE_1'\cap \mcE_{\Delta}$, 
	$$
	\max_{i\in[n]}\big\|\wtU\wtU^{\top}-\wtU^{(i)}\wtU^{(i)\top} \big\|_{\rm F}\lesssim \sqrt{\frac{\sigma^2(\lambda+\sigma^2)}{\lambda^2}}\cdot \frac{\sqrt{p(r+\log n)}}{n}.
	$$
	Basically,  the sensitivity upper bound is the same as in the Gaussian case since it is determined by the first order term $\max_{i\in[n]}\big\|U^{\top}(X_iX_i^{\top}/n)U_{\perp} \big\|=\max_{i\in[n]}\|U^{\top}X_i\|\|U_{\perp}^{\top}X_i\|/n$.  Similarly,  
	$$
	\sum_{k=1}^p \big(\lambda_k(\whSig)-\lambda_k(\whSig^{(i)})\big)^2\leq 2\|\whSig-\whSig^{(i)}\|^2\leq C_4\bigg(\frac{\lambda(r+\log n)+\sigma^2(p+\log n)}{n}\bigg)^2,
	$$
	where we used the upper bounds of $\|X_i\|^2$ and $\|X_i'\|^2$ stated in Lemma~\ref{lem:con-X-norm}.  The rest of the proof is skipped.

\subsection{Proof of Lemma \ref{lem:T-highterm}}
We begin with discussing the case  $\bfs = (k, 0, \cdots, 0)$ and $l=1$. By definition, we have
		\begin{align*}
			\| \mathcal{T}_{\Sigma, k, \mathbf{s}, 1}\left(\Delta-\Delta^{(i)}\right) \|_{\rm F} 
			& = \| Q^{-k} \left(\Delta-\Delta^{(i)}\right) Q^{\perp} \Delta Q^{\perp} \cdots Q^{\perp} \Delta Q^{\perp}\|_{\rm F}  \\
			& \leq \frac{1}{\lambda_r} \|U \left(\Delta-\Delta^{(i)}\right) U_{\perp}^{\top} \|_{\rm F} \left(\frac{\|\Delta\|}{\lambda_r} \right)^{k-1}. 
		\end{align*} 
		By Lemma~\ref{lem:con-bounds}, in the event $\mcE^{\ast}$, we get
		\begin{align*}
			\| \mathcal{T}_{\Sigma, k, \mathbf{s}, 1}\left(\Delta-\Delta^{(i)}\right) \|_{\rm F} \lesssim  \sqrt{\frac{\sigma^2 (\lambda+\sigma^2) }{\lambda^2}} \frac{\sqrt{p(r+\log n)}}{n} \left(\frac{D_{\max}}{\lambda_r} \right)^{k-3}. 
		\end{align*} 
		For the cases $\bfs = (k, 0, \cdots, 0)$ and $2 \leq l \leq k$, 
		\begin{align*}
			\| \mathcal{T}_{\Sigma, k, \mathbf{s}, l}\left(\Delta-\Delta^{(i)}\right) \|_{\rm F} &= \| Q^{-k} \Delta^{(i)} Q^{\perp} \cdots Q^{\perp} \left(\Delta-\Delta^{(i)}\right) Q^{\perp} \cdots Q^{\perp} \Delta Q^{\perp}\|_{\rm F} \\
			& \leq \Big\| Q^{-k} \Delta^{(i)} Q^{\perp} \cdots Q^{\perp} \frac{1}{n} X_iX_i^{\top}Q^{\perp} \cdots Q^{\perp} \Delta Q^{\perp}\Big\|_{\rm F} \\
			& \quad +  \Big\| Q^{-k} \Delta^{(i)} Q^{\perp} \cdots Q^{\perp} \frac{1}{n} X_i^{\prime}X_i^{\prime \top}Q^{\perp} \cdots Q^{\perp} \Delta Q^{\perp}\Big\|_{\rm F}.
		\end{align*}
		We need to control each term in the RHS of above inequality. The proof of the following lemma can be found in Section~\ref{sec:proof-tech-lem} of the Appendix. 
		
		\begin{lemma}\label{long_sequence}
			For all $k\geq 2$ and $2\leq l\leq k$, the following bounds hold in the event $\mcE^{\ast}$. 
			\begin{align*}
				\max_{i\in[n]}\Big\| Q^{-k} \Delta^{(i)} Q^{\perp} \cdots Q^{\perp} &\frac{1}{n} X_iX_i^{\top}Q^{\perp} \cdots Q^{\perp} \Delta Q^{\perp}\Big\|_{\rm F}\\
				& \lesssim  \left(\frac{D_{\max}}{\lambda_r} \right)^{k-2}\cdot \frac{\sigma^2}{\lambda}\sqrt{\frac{p}{n}} \cdot\sqrt{\frac{\sigma^2 (\lambda+\sigma^2) }{\lambda^2}} \frac{\sqrt{p(r+\log n)}}{n}, 
			\end{align*}
			and moreover
			\begin{align*}
				\max_{i\in[n]}\Big\| Q^{-k} \Delta^{(i)} Q^{\perp} \cdots &Q^{\perp} \frac{1}{n} X_i^{\prime}X_i^{ \prime\top}Q^{\perp} \cdots Q^{\perp} \Delta Q^{\perp} \Big\|_{\rm F} \\
				\lesssim& l \left(\frac{D_{\max}}{\lambda_r} \right)^{k-2}\cdot \frac{\sigma^2}{\lambda}\Big(\sqrt{\frac{p}{n}}+\frac{p}{n}\Big) \cdot\sqrt{\frac{\sigma^2 (\lambda+\sigma^2) }{\lambda^2}} \frac{\sqrt{p(r+\log n)}}{n}. 
			\end{align*}
			Note that the $Q^{\perp}$ term appears $l-1$ times before $X_iX_i^{\top}$ or $X_iX_i'^{\top}$ in the above product sequences of matrices. 
		\end{lemma}
		
		According to Lemma~\ref{long_sequence}, for $2\leq l\leq k$
		\begin{align*}
			& \| \mathcal{T}_{\Sigma, k, \mathbf{s}, l}\left(\Delta-\Delta^{(i)}\right) \|_{\rm F} \lesssim  (l+1) \left(\frac{D_{\max}}{\lambda_r} \right)^{k-2}\cdot \frac{\sigma^2}{\lambda}\Big(\sqrt{\frac{p}{n}}+\frac{p}{n}\Big) \cdot\sqrt{\frac{\sigma^2 (\lambda+\sigma^2) }{\lambda^2}} \frac{\sqrt{p(r+\log n)}}{n}. 
		\end{align*}
		Taking the summation over $l\in[k]$, we have 
		\begin{align*}
			\sum_{l\in[k]} \| \mathcal{T}_{\Sigma, k, \mathbf{s}, l}\left(\Delta-\Delta^{(i)}\right) \|_{\rm F} 
			& \lesssim \frac{(3+k)k}{2} \left(\frac{D_{\max}}{\lambda_r} \right)^{k-2}\cdot \frac{\sigma^2}{\lambda}\Big(\sqrt{\frac{p}{n}}+\frac{p}{n}\Big) \cdot\sqrt{\frac{\sigma^2 (\lambda+\sigma^2) }{\lambda^2}} \frac{\sqrt{p(r+\log n)}}{n},
		\end{align*} 
		which holds for $\bfs=\{k,0,\cdots,0\}$ and all $i\in[n]$.
		
		We now argue that the above bound holds for all $\bfs:=\{s_1,\cdots,s_{k+1}\}$ such that $\sum_j s_j=k$ and $s_j$'s are non-negative integers. 
		For any $\bfs\in I_k\setminus \{k, 0, \cdots, 0\}$, there exists at least two subsequence within $\mathcal{T}_{\Sigma, k, \mathbf{s}, l}\left(\Delta-\Delta^{(i)}\right) $ such that the sequence is starting from $UU^{\top}$ and ends with all other projectors being $U_{\perp}U_{\perp}^{\top}$ (Note that the other cases involving $U^{\top}(\Delta-\Delta^{(i)})U$ are smaller terms). Without loss of the generality, we discuss one subsequence where the first projector is $Q^{-\bar{k}}$ with $\bar{k}>0$ and the remaining projectors are all $Q_{\perp}$. Suppose that the target subsequence contains $t+1$ projectors $(t\leq k)$ with indices $(s_{m+1}, \cdots, s_{m+t}, s_{m+t+1})$ for some integer $m\geq 0$ satisfying $\sum_{j = m+1}^{m+t+1} s_j = \bar{k}$ and $\sum_{j=1}^{m} s_j + \sum_{j=m+t+2}^{j=k+1} s_j = k- \bar{k}$.  
		
		When $m+1\leq l \leq m+t+1$, the target subsequence times is upper bounded with the argument for the case $\bfs = (t, 0, \cdots, 0)$.  The remaining part in  $\mathcal{T}_{\Sigma, k, \mathbf{s}, l}\left(\Delta-\Delta^{(i)}\right) $  contains $k-t$ projectors and is upper bounded by $D_{\max}^{k-t}/\lambda^{k-\bar k}$. Therefore, 
		\begin{align*}
			\| \mathcal{T}_{\Sigma, k, \mathbf{s}, l}&\left(\Delta-\Delta^{(i)}\right) \|_{\rm F} \\ 
			& \lesssim \frac{D_{\max}^{k-t}}{\lambda^{k-\bar{k}}}\cdot \lambda^{t-\bar{k}}\cdot  (l+1) \left(\frac{D_{\max}}{\lambda_r} \right)^{k-2}\cdot \frac{\sigma^2}{\lambda}\Big(\sqrt{\frac{p}{n}}+\frac{p}{n}\Big) \cdot\sqrt{\frac{\sigma^2 (\lambda+\sigma^2) }{\lambda^2}} \frac{\sqrt{p(r+\log n)}}{n}\\
			& \lesssim  (l+1) \left(\frac{D_{\max}}{\lambda_r} \right)^{k-2}\cdot \frac{\sigma^2}{\lambda}\Big(\sqrt{\frac{p}{n}}+\frac{p}{n}\Big) \cdot\sqrt{\frac{\sigma^2 (\lambda+\sigma^2) }{\lambda^2}} \frac{\sqrt{p(r+\log n)}}{n}.   
		\end{align*} 
		The same argument can be applied to the case $l\in[k+1]\setminus [m+1: m+t+1]$ and is skipped here.  Finally we conclude that in event $\mcE^{\ast}$,
		\begin{align*}
			&  \sum_{l\in[k]}\| \mathcal{T}_{\Sigma, k, \bfs, l}(\Delta - \Delta^{(i)}) \|_{\rm F} \lesssim  \frac{(3+k)k}{2} \left(\frac{D_{\max}}{\lambda_r} \right)^{k-2}\cdot \frac{\sigma^2}{\lambda}\Big(\sqrt{\frac{p}{n}}+\frac{p}{n}\Big) \cdot\sqrt{\frac{\sigma^2 (\lambda+\sigma^2) }{\lambda^2}} \frac{\sqrt{p(r+\log n)}}{n}, 
		\end{align*}
		for $k\geq 2$ and any $\bfs \in I_k$.  This concludes the proof of Lemma~\ref{lem:T-highterm}.

		\section{Proofs of Technical Lemmas} 
		\label{sec:proof-tech-lem}
		
		
		\subsection{Proof of Lemma~\ref{lem:con-X-norm}}
		
		By the orthogonal invariance of normal distribution,  we can simply assume $\Sigma={\rm diag}(\lambda_1,\cdots,\lambda_p)$ and as a result,
		$$
		\|X\|^2\stackrel{{\rm d.}}{=}\lambda_1 Z_1^2+\cdots+\lambda_p Z_p^2,
		$$
		where $Z_1,\cdots, Z_p$ are i.i.d. standard normal random variables.  Note that $\lambda_i Z_i^2$ is sub-exponential ${\rm SE}(C_1\lambda_i^2,  C_2\lambda_i)$ for some absolute constants $C_1,C_2>0$ meaning that 
		$$
		\PP\big(|\lambda_i Z_i^2-\lambda_i|\geq t\big)\leq \exp\bigg(-c_1\min\bigg\{\frac{t^2}{\lambda_i^2}, \ \frac{t}{\lambda_i} \bigg\}\bigg), 
		$$
		for any $t>0$.   By the composition property of sub-exponential random variables,  we have 
		$$
		\lambda_1 Z_1^2+\cdots+\lambda_p Z_p^2\in {\rm SE}\Big(C_1\sum_{i=1}^p\lambda_i^2,  C_2\lambda_1\Big).
		$$
		Therefore,  we conclude that 
		$$
		\PP\bigg(\Big|\|X\|^2-{\rm tr}(\Sigma) \Big|\leq C_1\Big(u\sum_{i=1}^p \lambda_i^2\Big)^{1/2}+C_2\lambda_1 u\bigg)\geq 1-e^{-c_1u}, 
		$$
		for any $u>0$.  The proof of the upper bound $\|X_i\|^2, \|U^{\top}X_i\|^2$ and $\|U_{\perp}^{\top}X_i\|^2$ is straightforward and skipped here.

		\subsection{Proof of Lemma~\ref{lem:con-bounds}}
		
		Recall that $\Delta=\whSig-\Sigma$ and $\Delta^{(i)}=\whSig^{(i)}-\Sigma$.  By Lemma~\ref{lem:con-hatSigma},  we have 
		$$
		\EE\|\Delta\|+\max_{i\in[n]}\|\Delta^{(i)}\|\leq C_1\Bigg(\sqrt{\frac{(\lambda+\sigma^2)(r\lambda+p\sigma^2)}{n}}+\frac{r\lambda+p\sigma^2}{n}\Bigg)\leq C_2\sqrt{\frac{(\lambda+\sigma^2)(r\lambda+p\sigma^2)}{n}},
		$$
		where the last inequality is due to $C_3r\leq n$ and $\lambda/\sigma^2\geq C_3p/n$ for some large constant $C_3>0$.  These conditions also imply that ${\rm tr}(\Sigma)\leq n\|\Sigma\|$. Therefore,   by Lemma~\ref{lem:con-hatSigma},  we get 
		\begin{align*}
			\PP\Bigg(\Big|\|\Delta\|-\EE \|\Delta\|\Big|\leq C_4(\lambda+\sigma^2)\sqrt{\frac{t}{n}}\Bigg)\geq 1-e^{-t},
		\end{align*}
which holds as long as $0<t<c_0n$	for a small $c_0>0$.   By choosing $t=c_0(n\wedge p)$,  we conclude that 
		$$
		\PP\bigg(\|\Delta\|\leq C_4\sqrt{\frac{(\lambda+\sigma^2)(r\lambda+p\sigma^2)}{n}}+\frac{\lambda}{10}\bigg)\geq 1-e^{-c_0(n\wedge p)}.
		$$
Similar bounds can be also derived for each $\|\Delta^{(i)}\|$. By taking a union bound, we conclude that there exists an event $\mcE_{\Delta}$ with probability $\PP(\mcE_{\Delta})\geq 1-e^{-c_0'(n\wedge p)}$ under which the following bound holds
		$$
		\max_{i\in[n]} \|\Delta\|+\|\Delta^{(i)}\|\leq C_6\sqrt{\frac{(\lambda+\sigma^2)(r\lambda+p\sigma^2)}{n}}+\frac{\lambda}{10}.
		$$

		Now we study the bounds in event $\mcE_1$ and begin with $U^{\top}\Delta U_{\perp}$. By definition, we write
		$$
		U^{\top}\Delta U_{\perp}=\frac{1}{n}\sum_{i=1}^n U^{\top}X_iX_i^{\top}U_{\perp}. 
		$$
		Since $X_i\sim \mcN(0,\Sigma)$ and by the spiked structure of $\Sigma$,  we have 
		$$
		U^{\top}X_i\sim \mcN\big(0, \Lambda_r\big)\quad {\rm and}\quad U_{\perp}^{\top}X_i\sim \mcN\big(0, \sigma^2I_{p-r}\big),
		$$
		where $\Lambda_r:={\rm diag}\big(\lambda_1+\sigma^2,\cdots, \lambda_r+\sigma^2\big)$. Let $Z_i\sim \mcN(0, I_r), Y_i\sim \mcN(0, I_{p-r}), i\in[n]$ be independent Gaussian random vectors. Then
		$$
		U^{\top}X_i\stackrel{{\rm d.}}{=} \Lambda_r^{1/2}Z_i\quad {\rm and }\quad U_{\perp}^{\top}X_i\stackrel{{\rm d.}}{=} \sigma Y_i,\quad i\in[n]. 
		$$
		As a result, 
		$$
		U^{\top}\Delta U_{\perp}\stackrel{{\rm d.}}{=} \sigma \Lambda_r^{1/2} \frac{1}{n}\sum_{i=1}^n Z_iY_i^{\top}\Longrightarrow \big\|U^{\top}\Delta U_{\perp}\big\|\leq \frac{\sqrt{\sigma^2(\lambda+\sigma^2)}}{n}\Big\|\sum_{i=1}^n Z_i Y_i^{\top} \Big\|. 
		$$
		Denote  $A:= [Z_1, \cdots, Z_n]\in\RR^{r\times n}$ and $B: = [Y_1, \cdots, Y_n]^{\top}\in \RR^{n\times (p-r)}$. Since $B\in \RR^{n\times (p-r)}$ has i.i.d. $\mcN(0, 1)$ entries, we are able to write $B = [\tdY_1, \cdots, \tdY_{p-r}]$ with the column $\tdY_i$ i.i.d. to $\mcN(0, I_{p-r})$ for $i\in[p-r]$. We can therefore write $\sum_{i=1}^n Z_iY_i^{\top}=AB$. Conditioned on $A$, the matrix 
		$$
		AB=\big(A\tilde Y_1,\cdots, A\tilde Y_{p-r}\big), 
		$$
		has i.i.d. columns obeying the distribution $\mcN(0, AA^{\top})$. By Lemma~\ref{lem:con-hatSigma}, 
		\begin{align*}
			\PP\bigg(\|AA^{\top}/n-I_r\|\leq C_3\sqrt{\frac{t}{n}}+C_4\sqrt{\frac{r}{n}}\bigg)\geq 1-e^{-t},
		\end{align*}
		for any $t>0$ and  if $n\geq C_4r\log n$. Therefore, there exists an event $\mcF_1$ with $\PP(\mcF_1)\geq 1-e^{-c_1n}$ under which the following bound holds
		$$
		\frac{9n}{10}\leq \lambda_{\min}(AA^{\top})\leq \lambda_{\max}(AA^{\top})\leq \frac{11n}{10}.
		$$
		Under the event $\mcF_1$, we apply Lemma~\ref{lem:con-hatSigma} again and conclude with 
		$$
		\PP\bigg(\Big\|\frac{ABB^{\top}A^{\top}}{p-r}-AA^{\top}\Big\|\leq\frac{n}{100}\bigg)\geq 1-e^{-c_1p},
		$$
		where $c_1>0$ is a small constant and used the condition $p\geq 2r$. Finally, we conclude that 
		$$
		\PP\Big(\|AB\|\leq C_5\sqrt{np}\Big)\geq 1-e^{-c_1(p\wedge n)}.
		$$
		Therefore,
		\begin{equation}\label{eq:proof-hatUDeltaUperp}
			\big\|U^{\top}\Delta U_{\perp}\big\|\leq C_5\sqrt{\frac{p\sigma^2(\lambda+\sigma^2)}{n}}. 
		\end{equation}
		Similar bounds can also be established for $\|U^{\top}\Delta^{(i)}U_{\perp}\|$.  By taking a union bound, we get 
		\begin{align}\label{eq:proof-UDeltaUperp-bd1}
			\PP\bigg(\|U^{\top}\Delta U_{\perp}\|+\max_{i\in[n]} \|U^{\top}\Delta^{(i)}U_{\perp}\|\leq C_5\sqrt{\frac{p\sigma^2(\lambda+\sigma^2)}{n}}\bigg)\geq 1-e^{-c_1(p\wedge n)},
		\end{align}
		where we used the fact $n\geq C_1(r+\log n)$. Applying Lemma~\ref{lem:con-hatSigma}, we can show that the same bound also holds for $\EE \big\|U^{\top}\Delta U_{\perp}\big\|$.

		Similarly, $U^{\top}X_iX_i^{\top}U_{\perp}\stackrel{{\rm d.}}{=}\sigma \Lambda_r^{1/2}Z_iY_i^{\top}$ so that $\|U^{\top}X_iX_i^{\top}U_{\perp}\|\leq C\sqrt{\sigma^2(\lambda+\sigma^2)}\|Z_i\|\|Y_i\|$.  Applying Lemma~\ref{lem:con-X-norm} again, we get 
		\begin{equation}\label{eq:proof-UXXUperp}
			\PP\Big(\max_{i\in [n]} \|U^{\top}X_iX_i^{\top}U_{\perp}\|+\|U^{\top}X_i'X_i'^{\top}U_{\perp}\|\leq C_5\sqrt{\sigma^2(\lambda+\sigma^2)p(r+\log n)}\Big)\geq 1-n^{-100}. 
		\end{equation}
		
		We now consider the last terms in the event $\mcE_1$. Note that $X_i$ is independent of $\sum_{j\neq i}X_iX_j^{\top}$. Conditioned on the latter one, we have 
		$$
		U^{\top} \Big( \frac{1}{n} \sum_{j\ne i} X_j X_j^{\top} \Big) U_{\perp} U_{\perp}^{\top} X_i\sim \mcN\bigg(0, \sigma^2U^{\top} \Big( \frac{1}{n} \sum_{j\ne i} X_j X_j^{\top} \Big) U_{\perp} U_{\perp}^{\top} \Big( \frac{1}{n} \sum_{j\ne i} X_j X_j^{\top} \Big) U\bigg).
		$$
		By the bound (\ref{eq:proof-UDeltaUperp-bd1}), we get
		$$
		\Big\|U^{\top} \Big( \frac{1}{n} \sum_{j\ne i} X_j X_j^{\top} \Big) U_{\perp}  \Big\|\leq C_3\sqrt{\frac{\sigma^2(\lambda+\sigma^2)p}{n}},
		$$
		which holds in event $\mcE_{\Delta}$. Now applying Lemma~\ref{lem:con-X-norm} where $u=C(r+\log n)$, we conclude that in event defined in (\ref{eq:proof-UDeltaUperp-bd1}), 
		\begin{equation}\label{eq:proof-UXXjUperp}
			\PP\bigg(\Big\|U^{\top} \Big( \frac{1}{n} \sum_{j\ne i} X_j X_j^{\top} \Big) U_{\perp} U_{\perp}^{\top} X_i \Big\|\leq \sigma\cdot \sqrt{\frac{\sigma^2(\lambda+\sigma^2)p(r+\log n)}{n}}\bigg)\geq 1-n^{-100}. 
		\end{equation}
		
		By (\ref{eq:proof-UDeltaUperp-bd1}), (\ref{eq:proof-UXXUperp}), and (\ref{eq:proof-UXXjUperp}),  we conclude that $\PP(\mcE_1)\geq 1-e^{-c_1(n\wedge p)}-2n^{-99}$. 

		\subsection{Proof of Lemma~\ref{lem:hatU-bound}}
		
		According to Lemma~\ref{lem:spectral-formula}, 
		$$
		\widehat{U} \widehat{U}^{\top}-U U^{\top}=\sum_{k \geq 1} \mathcal{S}_{\Sigma, k}(\Delta). 
		$$
		For $k \geq 2$, there exists a term either of the form $Q^{-1} \Delta Q^{\perp}$ or $Q^{\perp} \Delta Q^{-1}$ in each summand of $\mathcal{S}_{\Sigma, k}(\Delta)$ and thus 
		$$
		\left\|\mathcal{S}_{\Sigma, k}(\Delta)\right\| \leq \binom{2k}{k}\left(\frac{\|\Delta\|}{\lambda_r}\right)^{k-1}\left\|Q^{-1} \Delta Q^{\perp}\right\| \leq \binom{2k}{k}\left(\frac{\|\Delta\|}{\lambda_r}\right)^{k-1} \frac{\left\|U^{\top} \Delta U_{\perp}\right\|}{\lambda_r}, 
		$$
		where we use the fact $$ \left\|Q^{-1} \Delta Q^{\perp}\right\|=\left\|U \Lambda^{-1} U^{\top} \Delta U_{\perp} U_{\perp}^{\top}\right\| \leq\left(\lambda_r\right)^{-1}\left\|U^{\top} \Delta U_{\perp}\right\|. 
		$$ 
		For all integers $k \geq 1$, $\binom{2(k+1) }{k+1} \binom{2k}{k}^{-1}=2(2 k+1)(k+1)^{-1} \leq 4$ and we have  
		\begin{align*}
			\sum_{k=2}^{\infty}\left\|\mathcal{S}_{\Sigma, k}(\Delta)\right\| & \leq \binom{4}{2} \frac{\|\Delta\|}{\lambda_r} \cdot 	\frac{\left\|U^{\top} \Delta U_{\perp}\right\|}{\lambda_r} \sum_{k=0}^{\infty} 4^k\left(\frac{\|\Delta\|}{\lambda_r}\right)^k \\ 
			& \leq 6 \frac{\|\Delta\|}{\lambda_r} \cdot \frac{\left\|U^{\top} \Delta U_{\perp}\right\|}{\lambda_r} \sum_{k=0}^{\infty}\left(\frac{4}{4+\delta}\right)^k \\
			& = \frac{6(4+\delta)}{\delta} \cdot \frac{\|\Delta\|}{\lambda_r} \cdot \frac{\left\|U^{\top} \Delta 	U_{\perp}\right\|}{\lambda_r}. 
		\end{align*}
		Combining the above result with $\left\|\mathcal{S}_{\Sigma, 1}(\Delta)\right\|=\left\|\Lambda^{-1} U^{\top} \Delta U_{\perp}\right\|$, we have 
		$$
		\left\|\widehat{U} \widehat{U}^{\top}-U U^{\top}\right\| \leq\left\|\mathcal{S}_{\Sigma, 1}(\Delta)\right\|+\sum_{k \geq 2}\left\|\mathcal{S}_{\Sigma, k}(\Delta)\right\| \leq\left\|\Lambda^{-1} U^{\top} \Delta U_{\perp}\right\|+\frac{6(4+\delta)\|\Delta\|\left\|U^{\top} \Delta U_{\perp}\right\|}{\delta \lambda_r^2}. 
		$$

		\subsection{Proof of Lemma~\ref{bridge}} 
		\label{sec:lem-bridge}
		
		Define, for each $k\in[n]$, 
		\begin{equation*}
			\begin{aligned}
				& \eta_k := \frac{\mu_k \wedge\nu_k}{1-t_k}, \\
				& \tilde{\mu}_k:= \frac{\mu_k-\mu_k\wedge \nu_k}{t_k}, \\
				& \tilde{\nu}_k:=  \frac{\nu_k-\mu_k\wedge \nu_k}{t_k}, 
			\end{aligned}
		\end{equation*}
		which are three probability measures on $(\Omega, \mcF)$. These measures provide the following decomposition of $\mu_k$ and $\nu_k$, 
		\begin{equation*}
			\begin{aligned}
				& \mu_k=(1-t_k) \eta_k + t_k \tilde{\mu}_k, \\
				& \nu_k=(1-t_k) \eta_k + t_k \tilde{\nu}_k. 
			\end{aligned}
		\end{equation*} 
		Since $\mu_k$ and $\nu_k$ share $\eta_k$ as the common part, we can make a coupling for 
		$$
		\mu = \mu_1 \times \cdots \times \mu_n \quad {\rm and}\quad \eta := \eta_1 \times \cdots \times \eta_n. 
		$$ 
		Create independent random variables as $ \tilde{X}_k\sim \tilde{\mu}_k, k\in[n]$ and $ Z_k \sim \eta_k, k\in[n]$.  Denote
		$$
		\tilde{X}:= (\tilde{X}_1, \cdots, \tilde{X}_n) \sim \tilde{\mu}_1\times \cdots \times \tilde{\mu}_n, 
		$$ 
		$$
		Z:= (Z_1, \cdots, Z_n)\sim \eta_1\times \cdots \times \eta_n.
		$$ 
		Let $\{W_k: W_k\sim \operatorname{Bern}(t_k)\}_{k\in[n]}$ be a collection of independent Bernoulli random variables that are also independent of $\tilde{X}$ and $Z$. Then, 
		$$
		X_k := W_k \cdot \tilde{X}_k+(1-W_k) \cdot Z_k, 
		$$
		is a random variable with the law $\mu_k$ for all $k\in[n]$. As a result,  the law of random vector $X:= (X_1, \cdots, X_n)^{\top}$ is 
		$$
		\mu = \mu_1 \times \cdots \times \mu_n. 
		$$
		Therefore, $(X, Z)$ is a coupling for $\mu$ and $\eta$.  Note that 
		$$
		X = W \odot \tilde{X}+(\mathbf{1}^{\top}_n-W) \odot Z ,
		$$
		where $\mathbf{1}_n$ is an $n$-dimensional all one vector and $\odot$ denotes the element-wise product. Since $f(t) = e^{2\varepsilon t}-t e^{\varepsilon}+t-1 \geq 0 $ for any $t\in[0, 1]$ and $\varepsilon\in[0, \infty)$, we get
		\begin{equation}
			\label{W}
			\mathbb{E}_{W} \exp \bigg(\varepsilon \sum_{k \in[n]} W_k\bigg) = \prod_{k \in[n]} \left[ (1-t_k) + t_k e^{\varepsilon} \right] \leq \exp \bigg(2\varepsilon \sum_{k \in[n]} t_k\bigg). 
		\end{equation}
		
		Finally, we have 
		\begin{align*}
			\EE_{A, X\sim \mu} \mathbb{I}&\big(  A(X)  \in E\big) =  \EE_{X\sim \mu}  \EE_{A} \mathbb{I}\big(  A(X)  \in E\big) \\
			& = \EE_{(X, Z, W)} \EE_{A} \mathbb{I}\big(  A(X)  \in E\big) \\
			& = \EE_{W}\EE_{\tilde{X}}\EE_{Z} \EE_{A}  \mathbb{I}\Big(  A\big( W \odot \tilde{X}+(\mathbf{1}_n-W) \odot Z\big) \in E \Big) \\
			& = \EE_{\tilde{X}}\EE_{Z}\sum_{W = w}\PP (W = w) \EE_{A}  \mathbb{I}\Big(  A\big( w \odot \tilde{X}+(\mathbf{1}_n-w) \odot Z\big) \in E \Big) \\ 
			& \leq \EE_{\tilde{X}}\EE_{Z} \sum_{W = w}\PP (W = w) \left(  \exp \left(\varepsilon \sum_{k \in[n]} w_k\right)  \EE_{A}  \mathbb{I}\big( A(Z) \in E \big) +\frac{\delta}{e^{\varepsilon}-1}\cdot \exp\left (\varepsilon\sum_{k \in[n]} w_k \right) \right) \\
			& =  \mathbb{E}_{W} \exp \left(\varepsilon \sum_{k \in[n]} W_k \right) \cdot  \left(\EE_{Z}  \EE_{A}  \mathbb{I}\big(  A(Z) \in E \big)+\frac{\delta}{e^{\varepsilon}-1} \right)\\
			& \leq  \exp \bigg(2\varepsilon \sum_{k\in[n]}t_k\bigg) \left(\EE_{A, Z\sim \eta}  \mathbb{I}\big(  A_Z  \in E\big)+\frac{\delta}{e^{\varepsilon}-1}\right),
		\end{align*}
		where the first inequality is due to the $(\varepsilon, \delta)$-differential privacy composition property of $A$ (notice that $Z$ and $\omega\odot\tilde X+({\bf 1}_n-\omega)\odot Z$ differs by ${\bf 1}_n^{\top} w$ entries) and the last inequality is by (\ref{W}). 
		
		In the same fashion, one can also establish a relation:
		$$
		\EE_{A, X\sim \eta} \mathbb{I}\big(  A(X)  \in E\big) \leq \exp\bigg(2 \varepsilon \sum_{k\in[n]}t_k\bigg)	\left(\EE_{A, X\sim \nu}\mathbb{I}\big(  A(X)  \in E\big) + \frac{\delta}{e^{\varepsilon}-1}\right) , 
		$$
		which completes the proof.

		\subsection{Proof of Lemma~\ref{lem:subG-bounds}}
		By definition, we write
		$$
		U^{\top}\Delta U_{\perp}=\frac{1}{n}\sum_{i=1}^n U^{\top}X_iX_i^{\top}U_{\perp}. 
		$$
		Denote $a_i=U^{\top}X_i\in\RR^r$ and $b_i=U_{\perp}^{\top}X_i\in\RR^{p-r}$.  Therefore,  $a_i$ and $b_i$ are sub-Gaussian random vectors with proxy-covariance matrices $U^{\top}\Sigma U$ and $U_{\perp}^{\top}\Sigma U_{\perp}$,  respectively.   However,  $a_i$ and $b_i$ can be dependent for sub-Gaussian distributions.  Nevertheless,  we can write 
		$$
		U^{\top}\Delta U_{\perp}=\frac{1}{n}\sum_{i=1}^n a_ib_i^{\top},
		$$
		which is a sum of independent rank-one random matrices.  To this end,  we apply the matrix Bernstein concentration inequality developed by \cite{koltchinskii2011neumann}.  Observe that 
		\begin{align*}
			\big\|\EE a_ib_i^{\top}b_i a_i \big\|\leq \EE \big\|a_ib_i^{\top}b_i a_i\big\|=\max_{\|v\|\leq 1}\EE \|b_i\|^2 \langle a_i, v\rangle^2\leq& \EE^{1/2}\|b_i\|^4\cdot \sup_{\|v\|\leq 1}\langle a_i, v\rangle^4\\
			\lesssim& p\sigma^2\cdot (\lambda+\sigma^2),
		\end{align*}
		where we used the facts $a_i,  b_i$ are sub-Gaussian random vectors.   Similarly,  we can get 
		$$
		\big\|\EE b_ia_i^{\top}a_i b_i^{\top} \big\|\leq \EE \|b_ia_i^{\top}a_i b_i^{\top}\|\lesssim r(\lambda+\sigma^2)\cdot \sigma^2.  
		$$
		Moreover,  
		$$
		\Big\|\big\|a_ib_i^{\top} \big\|\Big\|_{\psi_1}\lesssim \|a_i\|_{\psi_2}\|b_i\|_{\psi_2}\lesssim \sqrt{(\lambda+\sigma^2)\sigma^2rp}.
		$$
		By \cite[Proposition~2]{koltchinskii2011neumann},  with probability at least $1-(p+n)^{-10}$,
		\begin{align}
			\big\|U^{\top}\Delta U_{\perp} \big\|\leq& C_3\sqrt{\frac{\sigma^2(\lambda+\sigma^2)p\log(p+n)}{n}}+\frac{\sqrt{\sigma^2(\lambda+\sigma^2)rp}\log(p+n)\log r}{n}\nonumber\\
			\leq& C_3\sqrt{\frac{\sigma^2(\lambda+\sigma^2)p\log(p+n)}{n}}, \label{eq:proof-subG-bd1}
		\end{align}
		where the last inequality holds since $n\geq r\log p\log^2 r$.  We can also get 
		$$
		\EE \|U^{\top}\Delta U_{\perp}\|\leq C_3\sqrt{\frac{\sigma^2(\lambda+\sigma^2)p\log p}{n}}. 
		$$
		
		Moreover,  since
		\begin{align*}
			\Big\|U^{\top}(X_iX_i^{\top}/n)U_{\perp} \Big\|\leq \|a_i\|\|b_i\|/n,
		\end{align*}
		we get,  with probability at least $1-(p+n)^{-9}$, that
		\begin{align}\label{eq:proof-subG-bd2}
			\max_{i\in[n]}\ \Big\|U^{\top}(X_iX_i^{\top}/n)U_{\perp} \Big\|+\Big\|U^{\top}(X_i'X_i^{'\top}/n)U_{\perp} \Big\|\leq C_3\frac{\sqrt{\sigma^2(\lambda+\sigma^2)p(r+\log n)}}{n}.
		\end{align}
		 In the event where (\ref{eq:proof-subG-bd1}) and (\ref{eq:proof-subG-bd2}) hod,  we have 
		\begin{align*}
			\max_{i\in[n]}\|U^{\top}\Delta^{(i)}U_{\perp}\|\leq& \|U^{\top}\Delta U_{\perp}\|+\max_{i\in[n]}\ \Big\|U^{\top}(X_iX_i^{\top}/n)U_{\perp} \Big\|+\Big\|U^{\top}(X_i'X_i^{'\top}/n)U_{\perp} \Big\|\\
			\lesssim& \sqrt{\frac{\sigma^2(\lambda+\sigma^2)p\log(p+n)}{n}}.  
		\end{align*}
		
		Since $X_i$ is independent of $\sum_{j\neq i}X_jX_j^{\top}$,  
		\begin{align*}
			\max_{i\in[n]}\bigg\|U^{\top}\bigg(\frac{1}{n}\sum_{j\neq i}X_jX_j^{\top}\bigg)U_{\perp}U_{\perp}^{\top}X_i  \bigg\|\leq& \max_{i\in[n]}\bigg\| U^{\top}\bigg(\frac{1}{n}\sum_{j\neq i}X_jX_j^{\top}\bigg)U_{\perp}\bigg\|\sigma \cdot{\sqrt{(r+\log n)\log n}}\\
			\lesssim& \sigma\cdot  \sqrt{\frac{\sigma^2(\lambda+\sigma^2)p(r+\log n)\log(p+n)\log n}{n}},  
		\end{align*}
		where the first inequality holds with probability at least $1-n^{-10}$.

		\subsection{Proof of Lemma~\ref{long_sequence}}
		
		\emph{Simple case: $k=2$}. It suffices to consider $l=k=2$ and to bound the following terms 
		$$
		\Big\|U^{\top} \Delta U_{\perp} U_{\perp}^{\top} \Big( \frac{1}{n}X_i^{\prime} X_i^{\prime \top} \Big) U_{\perp}\Big\|_{\rm F} \quad {\rm and}\quad  \Big\|U^{\top} \Delta U_{\perp} U_{\perp}^{\top} \Big( \frac{1}{n}X_i X_i^{\top} \Big) U_{\perp}\Big\|_{\rm F}
		$$
		We only show prove the upper bound for $ \|U^{\top} \Delta U_{\perp} U_{\perp}^{\top} \left( X_i X_i^{\top}/n \right) U_{\perp}\|_{\rm F}$ since the other term can be bounded in exactly the same way. We decompose the term into two parts for decoupling purpose:   
		\begin{align*}
			&  \Big\|U^{\top} \Delta U_{\perp} U_{\perp}^{\top} \Big( \frac{1}{n}X_i X_i^{\top} \Big) U_{\perp}\Big\|_{\rm F} \\
			& \leq  \Big\|U^{\top} \Big( \frac{1}{n} X_i X_i^{\top} \Big) U_{\perp} U_{\perp}^{\top} \Big( \frac{1}{n}X_i X_i^{\top} \Big) U_{\perp}\Big\|_{\rm F} + \Big\|U^{\top} \Big( \frac{1}{n} \sum_{j\ne i} X_j X_j^{\top} \Big) U_{\perp} U_{\perp}^{\top} \Big( \frac{1}{n}X_i X_i^{\top} \Big) U_{\perp}\Big\|_{\rm F}. 
		\end{align*}
		By Lemma~\ref{lem:con-bounds} and in event $\mcE^{\ast}$, we have 
		\begin{align}
			\Big\|U^{\top} \Big( \frac{1}{n} X_i X_i^{\top} \Big) U_{\perp} U_{\perp}^{\top} \Big( \frac{1}{n}X_i X_i^{\top} \Big) U_{\perp}\Big\|_{\rm F} \leq& \Big\| U^{\top} \Big(  \frac{1}{n}  X_i X_i^{\top} \Big) U_{\perp}  \Big\|_{\rm F}^2\notag\\
			\lesssim& \bigg(\frac{\sqrt{\sigma^2(\lambda+\sigma^2)p(r+\log n)}}{n}\bigg)^2 \label{eq:proof-UXUperpXUperp-bd1}
		\end{align}
		Meanwhile, the other term can be written as 
		\begin{align*}
			\Big\|U^{\top} \Big( \frac{1}{n} \sum_{j\ne i} X_j X_j^{\top} \Big) U_{\perp} U_{\perp}^{\top} \Big( \frac{1}{n}X_i X_i^{\top} \Big) U_{\perp}\Big\|_{\rm F} 
			& =  \Big\|U^{\top} \Big( \frac{1}{n} \sum_{j\ne i} X_j X_j^{\top} \Big) U_{\perp} U_{\perp}^{\top} \frac{1}{\sqrt{n}}X_i\Big\| \Big\|  \frac{1}{\sqrt{n}}U_{\perp}^{\top}X_i\Big\|\\
			& \lesssim \Big\|U^{\top} \Big( \frac{1}{n} \sum_{j\ne i} X_j X_j^{\top} \Big) U_{\perp} U_{\perp}^{\top} \frac{1}{\sqrt{n}}X_i\Big\|\cdot \sigma\sqrt{\frac{p}{n}}
		\end{align*}
		where the last inequality is due to Lemma~\ref{lem:con-X-norm}. Finally, applying Lemma~\ref{lem:con-bounds}, we get in the event $\mcE^{\ast}$ that 
		\begin{align}\label{eq:proof-UXUperpXUperp-bd2}
			\Big\|U^{\top} \Big( \frac{1}{n} \sum_{j\ne i} X_j X_j^{\top} \Big) U_{\perp} U_{\perp}^{\top} \Big( \frac{1}{n}X_i X_i^{\top} \Big) U_{\perp}\Big\|_{\rm F} \lesssim \sigma^2\sqrt{\frac{p}{n}}\cdot \frac{\sqrt{\sigma^2(\lambda+\sigma^2)(r+\log n)}}{n}.
		\end{align}
		The bound in (\ref{eq:proof-UXUperpXUperp-bd2}) clearly dominates that in (\ref{eq:proof-UXUperpXUperp-bd1}) under the conditions in Lemma~\ref{lem:sense-U}. Finally, we conclude that in event $\mcE^{\ast}$, 
		$$
		\Big\|Q^{-2} \Delta U_{\perp} U_{\perp}^{\top} \Big( \frac{1}{n}X_i X_i^{\top} \Big) Q^{\perp}\Big\|_{\rm F} \lesssim \frac{\sigma^2}{\lambda}\sqrt{\frac{p}{n}}\cdot \frac{\sqrt{\sigma^2(\lambda+\sigma^2)(r+\log n)}}{n\lambda}.
		$$

		\noindent\emph{The cases $k\geq 3$}. 	For any $2\leq l\leq k$, we aim to bound 
		$$
		\| Q^{-k} \Delta^{(i)} Q^{\perp} \cdots Q^{\perp} \frac{1}{n} X_iX_i^{ \top}Q^{\perp} \cdots Q^{\perp} \Delta Q^{\perp}\|_{\rm F}\quad {\rm and}\quad \| Q^{-k} \Delta^{(i)} Q^{\perp} \cdots Q^{\perp} \frac{1}{n} X_i'X_i'^{ \top}Q^{\perp} \cdots Q^{\perp} \Delta Q^{\perp}\|_{\rm F},
		$$
		where there exist $l-1$ $Q_{\perp}$ factors before $X_iX_i^{\top}$ or $X_1'X_i'^{\top}$. We begin with the first term and note the following simple fact
		\begin{align*}
			& \| Q^{-k} \Delta^{(i)} Q^{\perp} \cdots Q^{\perp} \frac{1}{n} X_iX_i^{ \top}Q^{\perp} \cdots Q^{\perp} \Delta Q^{\perp}\|_{\rm F}\\
			& \leq \frac{1}{n}\|U^{\top} \Delta^{(i)}U_{\perp}U_{\perp}^{\top}\cdots U_{\perp}U_{\perp}^{\top}X_i X_i^{\top}U_{\perp}\|_{\rm F} \frac{D_{\max}^{k-l}}{\lambda_r^k}\\
			&\leq \|U^{\top} \Delta^{(i)}U_{\perp}U_{\perp}^{\top}\cdots U_{\perp}U_{\perp}^{\top}X_i \|\cdot \frac{\sigma\sqrt{p}}{n}\frac{D_{\max}^{k-l}}{\lambda_r^k},
		\end{align*}
		where the last inequality is due to Lemma~\ref{lem:con-X-norm} and recall $D_{\max}$ is the upper bound for $\|\Delta\|$ and $\|\Delta^{(i)}\|$ in event $\mcE^{\ast}$. 
		
		Since $\Delta^{(i)}$ is independent of $X_i$, conditioned on $\Delta^{(i)}$, we have 
		$$
		U^{\top} \Delta^{(i)}U_{\perp}U_{\perp}^{\top}\cdots U_{\perp}U_{\perp}^{\top}X_i\sim \mcN(0, \sigma^2U^{\top} \Delta^{(i)}U_{\perp}U_{\perp}^{\top}\cdots U_{\perp}U_{\perp}^{\top}\cdots U_{\perp}U_{\perp}^{\top}\Delta^{(i)}U^{\top} ). 
		$$
		Observe that 
		\begin{align*}
			\big\|U^{\top} \Delta^{(i)}U_{\perp}U_{\perp}^{\top}\cdots U_{\perp}\big\|\leq& \|U^{\top}\Delta^{(i)}U_{\perp}\|\cdot\|\Delta^{(i)}\|^{l-2}\\
			\leq& C_2\sqrt{\frac{(\lambda+\sigma^2)p\sigma^2}{n}}\cdot D_{\max}^{l-2},
		\end{align*}
		where the last inequality is due to Lemma~\ref{lem:con-bounds}. Similarly as the case $k=2$, we apply Lemma~\ref{lem:con-X-norm} with $u=C(r+\log n)$ and conclude, in the event $\mcE^{\ast}$, that 
		$$
		\max_{i\in[n]}\ \|U^{\top} \Delta^{(i)}U_{\perp}U_{\perp}^{\top}\cdots U_{\perp}U_{\perp}^{\top}X_i \|\leq \sigma^2\sqrt{\frac{p}{n}}\cdot D_{\max}^{l-2}\sqrt{(\lambda+\sigma^2)(r+\log n)}. 
		$$
		Therefore, we conclude, in the event $\mcE^{\ast}$, that
		\begin{align*}
			\| Q^{-k} \Delta^{(i)} Q^{\perp} \cdots Q^{\perp} \frac{1}{n} X_iX_i^{\top}Q^{\perp} \cdots Q^{\perp} \Delta Q^{\perp}\|_{\rm F}
			\lesssim \Big(\frac{D_{\max}}{\lambda}\Big)^{k-2}\cdot\frac{\sigma^2}{\lambda}\sqrt{\frac{p}{n}}\cdot \frac{\sqrt{\sigma^2(\lambda+\sigma^2)(r+\log n)}}{n\lambda}. 
		\end{align*}
		
		Next, we prove the upper bound of
		\begin{align*}
			\Big\| Q^{-k} \Delta^{(i)} Q^{\perp} \cdots Q^{\perp} \Big(\frac{1}{n} X_i^{\prime}X_i^{ \prime\top}\Big)Q^{\perp} \cdots Q^{\perp} \Delta Q^{\perp}\Big\|_{\rm F}. 
		\end{align*}
		It is slightly more complicated since now $\Delta^{(i)}$ and $X_i'X_i'^{\top}$ are dependent. However, only one term in the summands of $\Delta^{(i)}$ involves $X_iX_i'^{\top}$. Observe first that
		\begin{align*}
			& \| Q^{-k} \Delta^{(i)} Q^{\perp} \cdots Q^{\perp} \frac{1}{n} X_i^{\prime}X_i^{ \prime\top}Q^{\perp} \cdots Q^{\perp} \Delta Q^{\perp}\|_{\rm F}\\
			& \leq \frac{1}{n}\|U^{\top} \Delta^{(i)}U_{\perp}U_{\perp}^{\top}\cdots U_{\perp}U_{\perp}^{\top}X_i^{\prime} X_i^{\prime\top}U_{\perp}\|_{\rm F} \frac{\|\Delta\|^{k-l}}{\lambda^k}. 
		\end{align*}
		For decoupling purpose, we write $U^{\top} \Delta^{(i)}U_{\perp}U_{\perp}^{\top}\cdots U_{\perp}U_{\perp}^{\top}X_i^{\prime} X_i^{\prime\top}U_{\perp}$ as the summation of $l$ terms. 
		\begin{align*}
			& U^{\top} \Delta^{(i)}U_{\perp}U_{\perp}^{\top}\cdots U_{\perp}U_{\perp}^{\top}X_i^{\prime} X_i^{\prime\top}U_{\perp} \\
			& =  U^{\top}\left( \frac{1}{n}X_i^{\prime} X_i^{\prime\top} \right)U_{\perp}U_{\perp}^{\top}\Delta^{(i)} \cdots U_{\perp}U_{\perp}^{\top} \Delta^{(i)}  U_{\perp}U_{\perp}^{\top}X_i^{\prime} X_i^{\prime\top}U_{\perp} \\
			& \quad + U^{\top}\left(\Delta^{-i} \right)U_{\perp}U_{\perp}^{\top}\left(\frac{1}{n}X_i^{\prime} X_i^{\prime\top} \right)  \cdots U_{\perp}U_{\perp}^{\top}\Delta^{(i)} U_{\perp}U_{\perp}^{\top}X_i^{\prime} X_i^{\prime\top}U_{\perp} \\
			& \quad + \cdots \\
			&  \quad + U^{\top}\left(\Delta^{-i} \right) U_{\perp}U_{\perp}^{\top}\left(\Delta^{-i} \right)  \cdots U_{\perp}U_{\perp}^{\top}  \left(\frac{1}{n}X_i^{\prime} X_i^{\prime\top}\right) U_{\perp}U_{\perp}^{\top}X_i^{\prime} X_i^{\prime\top}U_{\perp} \\
			& \quad + U^{\top}\left(\Delta^{-i} \right) U_{\perp}U_{\perp}^{\top}\left(\Delta^{-i} \right)  \cdots U_{\perp}U_{\perp}^{\top}\left(\Delta^{-i} \right)  U_{\perp}U_{\perp}^{\top}X_i^{\prime} X_i^{\prime\top}U_{\perp}\\
			&=:g_1+g_2+\cdots+g_l,
		\end{align*}
		where $\Delta^{-i}:=n^{-1}\sum_{j\neq i}X_jX_j^{\top}$. 
		Clearly, the term $g_l$ can be bounded in the same fashion as before and we conclude that in event $\mcE^{\ast}$ with
		$$
		\|g_l\|_{\rm F}\lesssim  \sigma^2\sqrt{p}D_{\max}^{l-2}\cdot \sqrt{\frac{\sigma^2(\lambda+\sigma^2)p(r+\log n)}{n}}.
		$$
		On the other hand, for $m=2,\cdots, l-1$, in event $\mcE^{\ast}$,
		\begin{align*}
			\|g_m\|_{\rm F}\leq& \frac{1}{n}\Big\| \underbrace{U^{\top}\left(\Delta^{-i} \right) U_{\perp}U_{\perp}^{\top}\left(\Delta^{-i} \right)  \cdots U_{\perp}U_{\perp}^{\top}}_{\textrm{product of }\Delta^{-i} \textrm{ for } m-1 \textrm{ times}}  X_i'\Big\|  \Big\|\underbrace{X_i'^{\top}Q^{\perp}\Delta^{(i)}\cdots Q^{\perp}X_i'X_i'^{\top}U_{\perp}}_{\textrm{product of }\Delta^{(i)} \textrm{ for } l-m-1 \textrm{ times}}\Big\|\\
			\lesssim& \Big\| \underbrace{U^{\top}\left(\Delta^{-i} \right) U_{\perp}U_{\perp}^{\top}\left(\Delta^{-i} \right)  \cdots U_{\perp}U_{\perp}^{\top}}_{\textrm{product of }\Delta^{-i} \textrm{ for } m-1 \textrm{ times}}  X_i'\Big\| \cdot D_{\max}^{l-m-1}\sigma^3\frac{p^{3/2}}{n}\\
			\lesssim& D_{\max}^{l-3}\sigma^4\frac{p^{3/2}}{n}\cdot \sqrt{\frac{\sigma^2(\lambda+\sigma^2)p(r+\log n)}{n}}\\
			\lesssim& \sigma^2\sqrt{p}D_{\max}^{l-2}\cdot \sqrt{\frac{\sigma^2(\lambda+\sigma^2)p(r+\log n)}{n}},
		\end{align*}
		where the second inequality is due to Lemma~\ref{lem:con-X-norm}, the third inequality is similarly derived as in the case $k=2$, and the last inequality is due to $\lambda/\sigma^2\geq C_1(p/n+\sqrt{p/n})$.  
		
		For the term $g_1$, we write 
		\begin{align*}
			\|g_1\|_{\rm F}\leq& \frac{1}{n}\big\|U^{\top}X_i'X_i'^{\top}U_{\perp}U_{\perp}^{\top}\Delta^{(i)} \cdots U_{\perp}U_{\perp}^{\top} \Delta^{(i)}  U_{\perp}U_{\perp}^{\top}X_i^{\prime}\big\| \big\|U_{\perp}^{\top}X_i'\big\|\\
			&\lesssim \sigma^2\frac{p}{n}\cdot \big\|U^{\top}X_i'X_i'^{\top}U_{\perp} \big\|\cdot D_{\max}^{l-2}\\
			&\lesssim \sigma^2\frac{p}{n}\cdot D_{\max}^{l-2}\sqrt{\sigma^2(\lambda+\sigma^2)p(r+\log n)},
		\end{align*}
		which holds on in event $\mcE^{\ast}$.  Combining all the bounds of $g_1,\cdots,g_m$, we conclude that in event $\mcE^{\ast}$,
		\begin{align*}
			\max_{i\in[n]}\Big\| U^{\top} \Delta^{(i)}U_{\perp}U_{\perp}^{\top}\cdots U_{\perp}U_{\perp}^{\top}X_i^{\prime} X_i^{\prime\top}U_{\perp}\Big\|\lesssim lD_{\max}^{l-2}\sigma^2\Big(\sqrt{\frac{p}{n}}+\frac{p}{n}\Big)\sqrt{\sigma^2(\lambda+\sigma^2)p(r+\log n)}.
		\end{align*}
		Finally, it implies that 
		\begin{align*}
			\Big\| Q^{-k} \Delta^{(i)} Q^{\perp} \cdots Q^{\perp} &\Big(\frac{1}{n} X_i^{\prime}X_i^{ \prime\top}\Big)Q^{\perp} \cdots Q^{\perp} \Delta Q^{\perp}\Big\|_{\rm F}\\
			&\lesssim l\Big(\frac{D_{\max}}{\lambda}\Big)^{l-2}\frac{\sigma^2}{\lambda}\Big(\frac{p}{n}+\sqrt{\frac{p}{n}}\Big)\cdot \frac{\sqrt{\sigma^2(\lambda+\sigma^2)p(r+\log n)}}{\lambda n},
		\end{align*}
		which holds for all $i\in [n]$ in event $\mcE^{\ast}$ and completes the proof.

\end{document}